\documentclass[a4paper,twoside]{article}

\usepackage[english]{babel}
\usepackage[utf8x]{inputenc}
\usepackage{amsmath}
\usepackage{graphicx}
\usepackage[colorinlistoftodos]{todonotes}

\usepackage{amsmath}
\usepackage{epsfig,amsthm}
\usepackage{latexsym}
\usepackage{amsfonts}
\usepackage{amssymb}
\usepackage{amscd}
\usepackage{mathrsfs}
\usepackage[all,cmtip]{xy}
\usepackage{enumerate}
\usepackage{tikz}
\usepackage{arydshln}
\usepackage{bbm}

\usepackage{mathdots}
\usepackage[shortlabels]{enumitem}

\usepackage{color}

\usepackage{lipsum}

\usepackage[colorlinks]{hyperref}
\hypersetup{pdfstartview={FitH},         linkcolor=blue,  citecolor=green }

\numberwithin{equation}{section} 
\usepackage[a4paper,left=1in,right=1in,top=1in,bottom=1in,footskip=.35in]{geometry}

\usepackage{fancyhdr}
\AtEndDocument{%
  \par
  \medskip
  \begin{tabular}{@{}l@{}}%
  \textsc{Radboud University}
  \\
    \textsc{Heyendaalseweg 135
, 6525AJ Nijmegen, the Netherlands
}\\
    \textit{E-mail address}: \texttt{\href{mailto:geotam@science.ru.nl}{geotam@science.ru.nl}}
  \end{tabular}}

\newtheorem{thm}{Theorem}[section]
\newtheorem{cor}[thm]{Corollary}

\newtheorem{prop}[thm]{Proposition}

\theoremstyle{definition}
\newtheorem{dfn}[thm]{Definition}
\newtheorem{rmk}[thm]{Remark}

\newcommand{\Fo}{{F_{\bullet}}}

\newcommand{\qo}{{q_{\bullet}}}

\newcommand{\tr}{{\mathrm{tr}}}

\newcommand{\GL}{{\mathrm{GL}}}

\newcommand{\SO}{{\mathrm{SO}}}
\newcommand{\SP}{{\mathrm{Sp}}}

\newcommand{\Ind}{{\mathrm{Ind}}}
\newcommand{\cInd}{{\mathrm{cInd}}}

\newcommand{\antidiag}{{\text{anti-diag}}}
\newcommand{\diag}{{\mathrm{diag}}}

\newcommand{\Hom}{{\mathrm{Hom}}}
\newcommand{\End}{{\mathrm{End}}}

\newcommand{\disc}{{\mathrm{disc}}}

\newcommand{\Gal}{{\mathrm{Gal}}}

\newcommand{\Ad}{{\mathrm{Ad}}}

\usepackage{tocloft}

\title{Endoscopic liftings of epipelagic representations for classical groups}
\author{
 Geo Kam-Fai Tam (\href{mailto:geotam@science.ru.nl}{geotam@science.ru.nl})
 }
\date{Radboud University Nijmegen
\\[2ex]
\today}

\begin{document}
\maketitle


\begin{abstract}

Let $G$ be a p-adic classical group (orthogonal, symplectic, unitary) and $\pi$ be an epipelagic representation of $G$ defined by Reeder-Yu. Using M{\oe}glin's theory 
of extended cuspidal supports and Bushnell-Kutzko's theory of covering types, we determine the endoscopic lift of $\pi$ into the general linear group whose dual expresses the dual group of $G$ as a complex matrix group, explicitly in terms of the inducing type of $\pi$ that is extended from the character of the first Moy-Prasad filtration subgroup defined by a stable functional. We interpret the inducing type of $\pi$ via Stevens' construction of supercuspidal representations by skew semi-simple strata, and introduce the so-called epipelagic strata, which only require the residual characteristic p to be odd. As an application, we reprove the results of M. Oi on the endoscopic lifts of Gross-Reeder's simple supercuspidal representations of quasi-split classical groups.

\end{abstract}

\setcounter{tocdepth}{2}
\tableofcontents

\section{Introduction}

The local Langlands correspondence (LLC) states that, in layman  terms, the irreducible admissible representations of a connected reductive group $G$ over a p-adic field $F$ may be parameterized by morphisms, known as Langlands parameters, from the absolute Galois group $\Gal(\overline F/F)$ of $F$ into the L-group ${}^LG := \hat{G}\rtimes \Gal(\overline F/F)$ of $G$, where the last group is defined by a related group $\hat{G}$ known as the dual of $G$ over the field of complex numbers.  As a current conjecture, the LLC is proven true when $G$ is a general linear group \cite{Harris-Taylor}\cite{Henniart-simple-proof}\cite{Laumon-Rapoport-Stuhler}\cite{Scholze-LLC-GLn} or a classical group \cite{Arthur-book}\cite{Mok-unitary}\cite{Kaletha-Minguez-Shin-White}\cite{Ishimoto-odd-orthogonal}. \cite{Fargues-Scholze} recently announced a construction of Langlands parameters associated with irreducible representations of  a general reductive group which is compatible with many expected properties of the LLC, using deep algebraic-geometric methods.

One promising approach to prove the LLC, within the scope of the representation theory of p-adic reductive groups, is the theory of endoscopy \cite{Langlands-Shelstad} and its twisted analogue \cite{Kottwitz-Shelstad}. A major idea of this approach, among other important ones, states that, again laymanly speaking, if an L-group ${}^LG$ is contained in another such group ${}^LH$, then viewing a parameter of $G$ as a parameter of $H$ via this inclusion corresponds to a lifting process from the parametrized representations of $G(F)$ to those of $H(F)$, known as the endoscopic lift from $G$ to $H$. As a hypothetical strategy, if the LLC of $H$ is known, one may hope that the knowledge from $H$ can be descent to say something about, and eventually prove, the LLC of $G$. This strategy was shown to be successful in proving the LLC for general linear groups and classical groups in \emph{loc. cit.}

In this paper, we study the endoscopic lifting of epipelagic representations. This kind of representations were first introduced in \cite{Reeder-Yu}. They are irreducible, supercuspidal, and compactly induced from so-called types, a broad term for irreducible representations of compact-mod-center open subgroups, that are built from rather simple data, known as stable functionals, which can be viewed as dual vectors in the quotient of the two shallowest Moy-Prasad filtration subgroups of a parahoric subgroup at the barycenter of a specific facet in the Bruhat-Tits building $\mathcal B(G,F)$ of $G$. They generalize the simple supercuspidal representations introduced in \cite{Gross-Reeder}, in which case the facet is an alcove. We choose to focus on epipelagic representations not only because of their simplicity and popularity \cite{Romano-thesis}\cite{BK-epipelagic}\cite{Kaletha-epipelagic}\cite{epipelagic-unitary}, but also their importance to other subjects in the representation theory of reductive groups \cite{epipelagic-theta-correspondence}\cite{Epipelagic-rigid-local-systems}\cite{epipelagic-Brylinski-Deligne-cover}.

We therefore investigate how the inducing types of epipelagic representations of an endoscopic group $G$ are related to those of a target group $H$. The methodology we use in this paper applies to a classical group $G$ (orthogonal, symplectic, or unitary) with the target group $H$ being the general linear group expressing the dual group of $G$ as a matrix group. According to the theory of endoscopy, since representations of $G(F)$ with the same endoscopic lift comprise an L-packet, i.e., a (finite) set of representations related by the endoscopic character identity (which means they have the same Langlands parameter under the assertion of the LLC), our final results provide an explicit description of the LLC for epipelagic representations of classical groups.

Our methodology combines two theories which are purely local in nature: M{\oe}glin's theory \cite{Moeglin-classification-classical-groups,Moeglin-classification-unitary-groups,Moeglin-twisted-endoscopy-Langlands-parameters} of the reducibility of parabolically induced representations to determine the cuspidal supports of endoscopic liftings a.k.a. extended cuspidal supports, and Bushnell-Kutzko's theory \cite{BK-types} of covering types to translate the above reducibility into an analogous property for modules of (affine) Hecke algebras. This combination was shown to be successful to compute endoscopic liftings in specific cases \cite{BHS}, \cite{Lust-Stevens}, \cite{Tam-unramified-unitary, BT-ramified}, the last two of which were previously studied by the author. In a noteworthy case \cite{BHS} when $G$ is a symplectic group, this methodology is successful to describe the inertial class of the endoscopic lifting of {any} supercuspidal representation of $G(F)$, i.e., the endoscopic lifting is determined up to twists by unramified characters. Our focus on epipelagic representations yields relatively simple descriptions of their inducing types as well as the covering types underlying  representations parabolically induced from epipelagic representations.

Let's explain our methodology more precisely, and describe afterwards the main results, Theorem \ref{main theorem in the intro}. For the simplicity of this introduction, let's assume $G$ to be special odd orthogonal or symplectic, defined by a non-degenerate symmetric or skew-symmetric $F$-bilinear form $h$. Let $K$ be the maximal unramified extension of $F$, and $\mathbb F$ be the residual field of $F$ whose characteristic $p$ is odd. The first two steps of the Moy-Prasad filtration $G(K)_{x}\supset G(K)_{x,0_+}\supset G(K)_{x,0_{++}}$ renders successive quotients $\mathsf G_x := G(K)_{x}/G(K)_{x,0_+}$ and $\mathsf{V}_{x} := G(K)_{x,0_+}/G(K)_{x,0_{++}}$. The conjugation of $G(K)_{x}$ on the filtration subgroups defines a representation $\mathsf G_x\rightarrow \GL(\mathsf{V}_{x})$ over ${\mathbb F}$.

Let $\beta\in \mathsf{V}_{x}^*(\mathbb F)$ be a $\mathbb F$-linear functional on $\mathsf{V}_{x}$. By fixing a non-trivial character $\psi$ of $\mathbb F$, we define a character $\psi_\beta$ of $G(F)_{x,0_+}$ which is trivial on $G(F)_{x,0_{++}}$. With our choices of $G$, the normalizer $N(\psi_\beta)$ of $\psi_\beta$ in $G(F)_{x}$ has the quotient by $G(F)_{x,0_+}$ isomorphic to an elementary abelian 2-group $\{\pm1 \}^{\#I}$, where $I$ is an index set. With the assumption on $p$ being odd, $\psi_\beta$ extends to  a character $\lambda$ of $N(\psi_\beta)$. We now assume:
\begin{equation*}
\text{the functional $\beta\in \mathsf{V}_{x}^*(\mathbb F)$ is stable, in the sense of  geometric invariant theory (GIT) \cite{Mumford-Stabilityofprojectivevarieties}. 
}
\end{equation*}
The stability of $\beta$ then implies that the compactly induced representation $\pi_{}:=\cInd_{N(\psi_\beta)}^{G(F)}\lambda$ is irreducible, and is an epipelagic representation in the sense of \cite{Reeder-Yu}. We note that the character $\lambda$ is uniquely determined by $\beta$ and a tuple of signs $\{\lambda(\omega_i)\}_{i\in I}\in \{\pm1 \}^{\#I}$, where $\omega_i\in G(F)$ corresponds to the element $((1_{j})_{j\neq i},-1_i)\in N(\psi_\beta)/G(F)_{x,0_+}$.

A similar construction of epipelagic representations applies to general linear groups. For our purpose, it is more desirable to express an inducing type of a supercuspidal representation using a maximal simple type defined in \cite[Sec 5]{BK}. Indeed, the stable functional $\beta$ can be viewed as an elliptic regular  semi-simple element in the dual Lie algebra $\mathfrak g^*$ of $G$, and admits a decomposition $\oplus_{i\in I}\beta_i$ where $F[\beta_i]$ generates a field over $F$, except when $G$ is orthogonal and for a unique index $o\in I$, in which case $\beta_i=0$. Assume that $i\neq o$. In the epipelagic case, it turns out that $F[\beta_i]/F$ is totally ramified. Let $\tilde G_i = \GL_{[F[\beta_i]:F]}$ and $\tilde G_{i,\tilde x}$ be an Iwahori subgroup of $\tilde G_i (F)$ (and so $\tilde x$ is the barycenter of an alcove in the Bruhat-Tits building of $\tilde G_i$) whose pro-p unipotent radical $\tilde G_{i,\tilde x,0_+}$ affords the character $\psi_{2\beta_i}$. The normalizer $\tilde N(\psi_{2\beta_i})$ of $\psi_{2\beta_i}$ in $\tilde G_i(F)$ has the quotient by $\tilde G_{i,\tilde x,0_+}$ isomorphic to $F[\beta_i]^\times/(1+\mathfrak p_{F[\beta_i]/F}) \cong \mu_F \times \left<\varpi_i\right>$, where $ \mu_F$ is the subgroup of $F$ of roots of unity of order coprime to $p$, and $\varpi_i$ is a chosen uniformizer of $F[\beta_i]$. By the general construction of maximal simple types (see \cite[Sec 6.1]{BK}), $\psi_{2\beta_i}$ can always be extended to a character $\tilde{\boldsymbol{\lambda}}_i$ of $\tilde N(\psi_{2\beta_i})$, and $\tilde\pi_i:=\cInd_{\tilde N(\psi_{2\beta_i})}^{\tilde G_i(F)}\tilde{\boldsymbol{\lambda}}_i$ is also an epipelagic representation of $\tilde G_i(F)$.

We now view $M:=\tilde G_i\times G$ as a Levi subgroup in a parabolic subgroup $P$ of a classical group $\mathbf G$, of the same type as $G$ but of a higher rank. Define a family of normalized parabolically induced representation $I(s,\tilde\pi_i,\pi):=\iota_{P}^{\mathbf G(F)}(\tilde\pi_i|\det|^s\times \pi)$ parametrized by $s\in \mathbb C$, then M{\oe}glin's theory in \emph{loc. cit.} asserts that $\tilde\pi_i$ belongs to the extended cuspidal support of $\pi$ if and only if $I(s,\tilde\pi_i,\pi)$ is reducible at a half-integer $s\in \tfrac{1}{2}\mathbb Z$ with $s\geq 1$ (with an extra parity condition which is irrelevant in our present setup). In the epipelagic case, we always have $I(s,\tilde\pi_i,\pi)$ reducible at exactly one of $\{0,\tfrac{1}{2},1\}$. Therefore, our objective is to determine those $\tilde\pi_i$ giving $I(1,\tilde\pi_i,\pi)$ reducible.

We can now state our main result on the reducibility points of $I(s,\tilde\pi_i,\pi)$, as a slight simplification of Propositions \ref{first general form of lifting, unitary group} and \ref{first general form of lifting}.

\begin{thm}
\label{main theorem in the intro}
Let $G$ be a special odd orthogonal or symplectic group over $F$, and $\pi_{}:=\cInd_{N(\psi_\beta)}^{G(F)}\lambda$ be an epipelagic representation constructed by the data $(\beta = \oplus_{i\in I}
\beta_i, \{\lambda(\omega_i)\}_{i\in I})$. For each index $i\in I$, with $i\neq o$ when $G$ is orthogonal, we define a character $\tilde{\boldsymbol{\lambda}}_i$ of $\tilde N(\psi_{2\beta_i})$ as follows:
\begin{equation*}
\tilde{\boldsymbol{\lambda}}_i|_{\tilde G_{i,\tilde x,0_+}} = \psi_{2\beta_i},\quad 
\tilde{\boldsymbol{\lambda}}_i|_{\mu_{F}}= \left(\frac{\cdot}{\mu_{F}}\right)^{\epsilon_G},
\quad\text{and}\quad
\tilde{\boldsymbol{\lambda}}_i(\varpi_i)= \tilde{\boldsymbol{\lambda}}_i(-2)\lambda(\omega_i )
\mathfrak{n}_z(\varpi_i,\beta,\psi,h).
\end{equation*}
Here 
\begin{itemize}
\itemsep0em 
\item $\epsilon_G=1$ (resp. $-1$) if $G$ is orthogonal (resp. symplectic), i.e., $h$ is an $\epsilon_G$-bilinear form,
\item $\left(\tfrac{\cdot}{\mu_{F}}\right)$ is the quadratic character of $\mu_{F}$, and 

\item $\mathfrak{n}_z(\varpi_i,\beta,\psi,h)$ is a normalized quadratic Gauss sum (i.e., a fourth root of unity), defined by a quadratic form on a certain $\mathbb F$-space $\mathfrak{W}_z$ related to a covering type for $I(s,\tilde\pi_i,\pi)$ in $\mathbf G(F)$, and which depends on the given data $(\varpi_i,\beta,\psi,h)$.
\end{itemize}
When $G$ is symplectic, we define a character $\tilde\pi_o$ of $F^\times$ by 
$$\tilde\pi_o|_{1+\mathfrak p_F}\equiv 1,\quad  \tilde\pi_o|_{\mu_F} = \left(\tfrac{\cdot}{\mu_F}\right)^{\#I} 
\quad\text{and}\quad
\tilde\pi_o(\varpi_{}) = \prod_{i\in I}\left(\tfrac{\varpi\det\beta_i}{\mu_F}\right),$$
Let $\tilde\pi_i:=\cInd_{\tilde N(\psi_{2\beta_i})}^{\tilde G_i(F)}\tilde{\boldsymbol{\lambda}}_i$  be the associated epipelagic representation of $\tilde G_i(F)$, then  $I(s,\tilde\pi_i,\pi)$ is reducible at $s=1$. 
\qed\end{thm}

The main technicality of the proof is exhibited in Section \ref{section Main calculation}, especially in Section \ref{subsection Expanding the intertwining operator as a sum}. A little more work in Section \ref{subsection Reducibility results for different classical groups} allows us to compute the values of $\mathfrak n_z(\varpi_i,\beta,\psi,h)$. A simple application of M{\oe}glin's theory on the bound of the number of such $\tilde\pi_i$ then implies that $\{\tilde\pi_i\}_{i\in I}$ comprises the extended cuspidal support of $\pi$, i.e., the endoscopic lift of $\pi$ has cuspidal support exactly $\{\tilde\pi_i\}_{i\in I}$.

We remark that analogous, although a bit more complicated, results hold when $G$ is special even orthogonal and unitary. For these two types of classical groups, there is an extra structural complicacy in the inducing types coming from the stability of the functionals. We skip the detail for the sake of this introduction, but refer interested readers to the classification results of stable gradings in \cite[Sec 7.2]{Gross-Levy-Reeder-Yu} and their interpretations in terms of semi-simple strata in Section \ref{subsection Stable functionals for classical groups} of this paper for details.

We hence conclude that, with our methodology, we can compute the endoscopic lifts, or more precisely their inducing types, and hence the L-packets of all epipelagic representations of all classical groups. These results are given in Section \ref{subsection Expanding the intertwining operator as a sum} of this paper.

To compute the reducibility points as stated in Theorem \ref{main theorem in the intro}, we apply Bushnell-Kutzko's theory to functorially identity the Bernstein component of $\tilde\pi_i\times \pi$ in $\mathbf G$ with the category of modules of the intertwining algebra, i.e., the desired Hecke algebra, of a covering type over $\tilde\lambda_i\times \lambda$ in $\mathbf G$, where $\tilde\lambda_i $ is the restriction of $\tilde{\boldsymbol{\lambda}}_i$ to the maximal compact subgroup of $\tilde N(\psi_{2\beta_i})$. The structure of this Hecke algebra is well known due to Lusztig \cite{Lusztig-finite-classical-groups}: it is of the generic type on an infinite dihedral group \cite[Th 1.2]{Stevens-Miyauchi}. We may take from such an algebra two generators, denoted by $T_y$ and $T_z$ in this paper, each satisfies a quadratic equation whose coefficients can be determined by the structure of the covering type (Section \ref{subsection Structures of Hecke algebras}).

The reducibility of $I(s,\tilde\pi_i,\pi)$ is hence converted into that of the corresponding module $X_s$ of the intertwining algebra. By the above functoriality, $X_s$ is induced from the character $D_s$ of the intertwining algebra of $\tilde\lambda_i\times \lambda$ in $M$ corresponding to $\tilde\pi_i|\det|^s\times \pi$. A crucial observation from \cite[(1.13)]{Blasco-Blondel-SP4} implies that, when $X_s$ is reducible, there are two  ways to express the eigenvalues of the product $T_y* T_z$: one by multiplying the respective eigenvalues of the generators from the quadratic equations, and another from the inducing character $D_s$. Equating the two expressions leads to a formula (\ref{main formula of reducibility points}) for computing the points of reducibility, a formula that we will examine in Section \ref{subsection Expanding the intertwining operator as a sum}; see Section \ref{section Reducibility} for the complete details of deducing this formula.

Cuspidal inducing types and the related covering types can be built from arithmetic data known as semi-simple strata, based on the general theories developed in \cite{BK}, \cite{Stevens-supercuspidal}. To apply these theories to epipelagic representations, we will translate the language of stable functionals into the one of semi-simple strata in Section \ref{subsection Stable functionals for classical groups}. For instance, the stability condition on functionals is translated into an elliptic-regularity conditions on the vectors in $\mathfrak g^*$ (see Proposition \ref{equivalent conditions for stability}), and the classification in \cite[Sec 7.2]{Gross-Levy-Reeder-Yu} of the points $x\in \mathcal B(G,F)$ such that $\mathsf{V}^*_{x}$ contains stable functionals is translated into conditions on the semi-simple strata that define the same types of parahoric subgroups admitting characters arising from stable functionals (see Corollary \ref{stability implies totally ramified} and the list that follows). Since these translated conditions can be intrinsically stated within the scope of semi-simple strata, we will eventually define what we will call epipelagic strata, which are then extended to characters inducing epipelagic representations (see Section \ref{subsection Epipelagic inducing types for classical groups}). This translation process is at the advantage that the residual characteristic $p$ of the base field $F$ is required only to be odd, the condition on $p$ for exhaustively constructing 
supercuspidal representations of classical groups in \cite{Stevens-supercuspidal}.

With our methodology explained above, we reprove, and compare with ours, the results of M. Oi \cite{Oi-SO-odd,Oi-Sp-and-SO-even,Oi-U-unram} on  endoscopic liftings of simple supercuspidal representations of $G$ in Section \ref{section Examples for simple supercuspidals}. Oi's methodology of computing the endoscopic liftings mainly uses the endoscopic character identity to compare the character expansions of $\pi$ and $\tilde\pi_i$ at affine generic elements in terms of Gauss and Kloostermann sums. His method applies to quasi-split classical groups except ramified unitary groups, while our method is applicable to all pure inner forms of quasi-split classical groups as long as the cohomological invariant of the Hermitian form defining $G$ is concerned (see Sections \ref{subsection Cohomological classification} and \ref{section Embeddings of lattices}). Moreover, our method applies to a more general kind of representations, the epipelagic ones.  Actually, we will use the inducing types of simple supercuspidal representations, i.e., the affine generic characters, as `building blocks' (as in Section \ref{subsection Liftings of epipelagic representations for classical groups}) to describe the inducing types of epipelagic representations explicitly by signs and quadratic Gauss sums.

Speaking of the endoscopic character identity,  Kaletha's theory \cite{Kaletha-epipelagic} applies to tamely ramified reductive p-adic groups to construct explicitly the L-packets of epipelagic representations (among other important results such as genericity and the formal degree conjecture), with some stricter conditions on the residual characteristic p than just being odd. Since we construct the inducing types of epipelagic representations based on the theory of semi-simple strata in \cite{Stevens-supercuspidal}, our method only requires $p$ to be odd. It would be interesting to compare the resulting endoscopic liftings with the author's.

\subsection{Acknowledgements}

This research is supported by the open competition of NWO under Grant No. OCENW.M20.132, and partially by the Radboud Excellence Initiative. We thank Maarten Solleveld for reading several draft versions of this article with great interest. In the midst of the writing, the author was supported by the encouragements from Julia Gordon, to whom the author is indebted. This article would not exist without the contributions of these two mathematicians.

\subsection{Notations and conventions}

The cardinality of a finite set $X$ is denoted by $\#X$. 

The action of a group $G$ on a set $X$ is denoted by $x\mapsto {}^gx$, for $g\in G$ and $x\in X$.

The real and imaginary parts of a complex number $s\in \mathbb C $ is respectively denoted by $\Re(s)$ and $\Im(s)$.

Given a matrix $A$, we denote the transpose of a matrix $A$ is denoted by ${}^tA$. Suppose that $A$ has entries in a field $K$, and $K/K_\bullet$ is a quadratic extension, we denote by $\overline A$ the entrywise conjugation of $A$.

We denote a diagonal matrix with entries $A_1,\dots,A_n$ from NW to SE by $\diag(A_1,\dots,A_n)$, and an anti-diagonal matrix with entries $A_1,\dots,A_n$ from NE to SW by $\antidiag(A_1,\dots,A_n)$. Each entry $A_i$ is allowed to be also a matrix.

A representation $\pi$ of a topological, locally profinite group $G$ is assumed to be smooth. We write $(G,\pi)$ if we want to emphasize the underlying group of the representation $\pi$. In many cases, $G$ is the subgroup of rational points of a reductive group over a non-Archimedean local field, and $\pi$ is a smooth representation of $G$, usually irreducible and supercuspidal. To simplify the writing, all supercuspidal representations of a reductive group is presumed to be irreducible unless otherwise specified.

Given a non-Archimedean local field $F$, denote $\mathfrak o_F$ the ring of integers with maximal ideal $\mathfrak p_F$. The residual field $\mathbb F:=\mathfrak o_F/\mathfrak p_F$ has cardinality $q = q_F$ which is a power of a prime number $p$. By fixing a uniformizer $\varpi$ of $F$, the multiplicative subgroup $F^\times$ of $F$ then decomposes as $\left<\varpi\right>\times \mathfrak o_F^\times$, and also $ \mathfrak o_F^\times =  \mu_F\times \mathcal{U}^1(F)$, where $\mu_F$ is the subgroup of root of unity of order coprime to $p$, and $\mathcal{U}^1(F)= 1+\mathfrak p_F$.

If $H$ is a finite cyclic group of even order, denote by $\left(\frac{\cdot }{H}\right):H\rightarrow \{\pm 1\}$ the quadratic character of $H$. Let $\psi:\mathbb F\rightarrow \mathbb C^\times$ be a non-trivial additive character. We denote the normalized quadratic Gauss sum by $\mathfrak{n}_\psi := q^{-1/2}\sum_{x\in \mathbb F^\times }\psi(x)\left(\frac{x}{\mathbb F^\times}\right)$. It is known that $\mathfrak{n}_\psi^2  = \left(\frac{-1}{\mathbb F^\times}\right)$, so that $\mathfrak{n}_\psi$ is a 4th root of unity.

The extended real line \cite[(6.4.1)]{Bruhat-Tits-reductive-group-1} is $\tilde{\mathbb R} = \{r,r_+:r\in \mathbb R\}$, equipped with an order $r>s$ if and only if $r>r_+\geq  s$. We will work with filtrations of groups compatibly parametrized by $\tilde{\mathbb R}$. To match with the enumeration of Moy-Prasad filtrations, unless otherwise specified, we parametrize any lattice sequence $\Lambda$ using the normalized valuation, i.e., $\varpi\Lambda(r) = \Lambda(r+1) $ for all $r\in \mathbb R$.

 \section{Epipelagic Representations}

 Let $F$ be a non-Archimedean local field with residual characteristic $p$, and $G$ be a connected reductive group over $F$. Let $K$ be the maximal unramified extension of $F$ in its separable closure $F^{\text{sep}}$. Fix a maximal $F$-split torus $S$ of $G$, and take a maximal $K$-split $F$-torus $T$ containing $S$.

The torus $T$ determines an apartment $\mathcal{A} =\mathcal{A}(T,K)$ in the Bruhat-Tits building $\mathcal{B}(G) = \mathcal{B}(G,K)$ of $G$, as well as an affine root system $\Phi_{\mathrm{aff}}$ containing the underlying (relative) root system $\Phi = \Phi(G,T)$. Fix a special point $o$ and an alcove $\mathcal{C}$ in $\mathcal A$ containing $o$. The boundary hyperplanes of $\mathcal C$ determine a set of simple roots $\Delta_{\mathrm{aff}}$, consisting of the roots in a simple root system $\Delta\subset \Phi$ together with an affine root $\alpha_0:=1-\alpha_l$, where $\alpha_l$ is the longest root in $\Phi$.

We take a point $x\in \mathcal{B}(G,F) = \mathcal B(G)^{\Gal(K/F)}$, and assume it lies in $\mathcal{A}$ by $G$-translation. For $r\in \tilde{\mathbb R}_{\geq 0}$, denote by $G(K)_{x,r}$ the Moy-Prasad filtration subgroup in $G(K)$, and put $G_{x,r}=G(F)_{x,r}:=G(K)_{x,r}^{\Gal(K/F)}$. The stabilizer group $G(F)_{x}$ of $x$ in $G(F)$ contains the parahoric subgroup $G_{x,0}$ with finite index modulo center. Now take $ r(x)$ to be the minimal positive value in $\{\psi(x):\psi\in \Phi_{\mathrm{aff}}\}$. We call an irreducible representation of $G$ {\bf epipelagic} if it has depth $r(x)$ and contains a non-zero vector fixed by $G_{x,r(x)_+}$.

 \subsection{Stable functionals}
 \label{subsection Stable functionals}

The maximal compact subgroup $T_0$ of $T(K)$ acts on the affine root subgroup $U_{\psi} = U_{\psi}(\mathfrak o_K)$ for each $\psi\in \Phi_{\mathrm{aff}}$, and on the quotient $\mathfrak g_{\dot\psi}:=U_{\psi}/U_{\psi_+}$, where $U_{\psi_+}$ is the next filtration subgroup of $U_{\psi}$ in $U_{\dot\psi}$ and ${\dot\psi}\in \Phi$ is the direction of $\psi$. Denote by $T_{0_+}$ the kernel of the action of $T_0$ on all $\mathfrak g_{\alpha}$ for all $\alpha\in \Delta$, so that $\mathfrak g_{\alpha}$ is the $\alpha$-weight subspace of the Lie algebra $\mathfrak g$ of $G(\overline{\mathbb F})$ under the action of $T_0/T_{0_+}$.

Put
$\Phi_{x,r} = \{ \alpha\in \Phi:
  \alpha =   \dot\psi\text{ for some $\psi \in \Phi_{\mathrm{aff}}$}
    \text{ such that $\psi(x)= r$}\}$ for $r\in \tilde{\mathbb R}$, 
     then we have a decomposition on the quotient
    $$
\mathsf{V}_{x,r} := G(K)_{x,r}/G(K)_{x,r_+}\cong \bigoplus_{\alpha\in \Phi_{x,r}
}\mathfrak g_{\alpha}.$$ 
A functional $\beta\in {\mathsf{V}}^*_{x,r}$ the linear dual space of ${\mathsf{V}_{x,r}}$ is called {\bf stable} for the action of $\mathsf G_x := G(K)_{x,0}/G(K)_{x,0_+}$ if the $\mathsf G_x $-orbit of $\beta$ in ${\mathsf{V}}^*_{x,r}$ is closed and the stabilizer of $\beta$ in $\mathsf G_x $ is finite (i.e., a finite algebraic group).

 Fix an additive character $\psi$ of $F$ throughout the paper, which is trivial on $\mathfrak p_F$ but non-trivial on $\mathfrak o_F$, and denote also by $\psi$ the induced character of $\mathbb F = \mathfrak o_F/\mathfrak p_F$. Given $x\in \mathcal B(G,F)$ and $\beta\in {\mathsf{V}}^*_{x,r}(\mathbb F) := ({\mathsf{V}}^*_{x,r})^{\Gal(\bar{\mathbb F}/\mathbb F)}$, we define a character on the compact subgroup $G(F)_{x,r}$ by 
 $$\psi_\beta: G(F)_{x,r} \xrightarrow{} (G_{x,r}/G_{x,r_+} )^{\Gal(\bar{\mathbb F}/\mathbb F)} \cong \mathsf{V}_{x,r}(\mathbb F) \xrightarrow{\beta} \mathbb F\xrightarrow{\psi}\mathbb C^\times.$$ 
Denote by $G(F)_{x,\beta}$ the stabilizer of $\beta$ in $G(F)_x$. We now take $r=r(x)$. If  $\beta$ is moreover stable, we take an irreducible representation $(G(F)_{x,\beta},\lambda)$ containing $( G_{x,r}, \psi_\beta)$, i.e., $\lambda$ an irreducible constituent in $\Ind_{G(F)_{x,r}}^{G(F)_{x,\beta}}\psi_\beta$. Then $\cInd_{G(F)_{x,\beta}}^{G(F)} \lambda$ is an irreducible supercuspidal representation of $G(F)$, and is moreover epipelagic \cite[Prop 2.4]{Reeder-Yu}.

\begin{rmk} The existence of rational stable functions in ${\mathsf{V}}^*_{x,r}$ is first shown in \cite[Sec 5 and 6]{Reeder-Yu} for large $p$, and is then guaranteed by \cite{Fintzen-Romano-stable-vectors} for arbitrary $p$ if we extend from $F$ to a large enough unramified extension of $F$.
 \qed\end{rmk}

\subsection{Example: Simple supercuspidal representations}
\label{subsection Example: Simple supercuspidal representations}

Simple supercuspidal representations are examples of epipelagic representations. These representations are first constructed in \cite{Gross-Reeder} for simply connected simple $F$-split groups, using affine generic characters. For quasi-split reductive groups, we recall the statements from \cite{Oi-SO-odd}.

With the setup of the previous section, we take $x\in \mathcal C$ to be the barycenter, i.e., the point where the affine roots in $\Delta_{\mathrm{aff}}$ attain a common value, which is $r = 1/h$, where $h$ is the twisted Coxeter number \cite{Reeder-torsion-autom}. Put $\mathcal{I}^+ = G_{x,r}$ and $\mathcal{I}^{++} = G_{x,r_+}$, and denote by $Z$ the center of $G(F)$. We call a character $\chi$ of $\mathcal{I}:= Z\mathcal{I}^+$ {\bf affine generic} if it is trivial on $\mathcal{I}^{++}$ and is non-trivial on all $U_{\alpha}/U_{\alpha_+}$ for $\alpha\in \Delta_{\text{aff}}$. In this paper, we also call the functional $\beta$ that gives rise to $\psi_\beta = \chi|_{\mathcal I^+}$ {\bf affine generic}.

Denote by $N(\chi)$ the normalizer of $\chi$ in $G(F)$. Given an affine generic character $\chi$ of $Z\mathcal{I}^+$, the compact induction $\cInd_{Z\mathcal{I}^+}^{G(F)}\chi$ admits a decomposition
$$\cInd_{\mathcal{I}}^{G(F)}\chi = \bigoplus_{\chi'} (\dim {\chi'} )\pi_{\chi'},$$
where $\chi'$ ranging over irreducible constituents of $\cInd_{Z\mathcal{I}^+}^{N(\chi)}\chi$ and $\pi_{\chi'} = \cInd_{N(\chi)}^{G(F)} \chi'$. Each $\pi_{\chi'}$ is an irreducible supercuspidal representation of $G(F)$. \cite[Prop 2.8]{Oi-SO-odd}.

Recall the decomposition $N_G(T)/T_0 = W_{\mathrm{aff}}\rtimes \bar\Omega$, where $\bar \Omega$ is the (finite) group $\Gamma_{\mathcal C}$ defined in \cite[p.189]{Bourbaki-Lie-4-6}, and let $\Omega$ be a set of representatives of $\bar\Omega$ in $N_G(T)$. If $\pi_{\chi'}$ and $\pi_{\xi'}$ are two such representations arising respectively from affine generic characters $\chi$ and $\xi$ of $\mathcal{I}$, then
\begin{equation}
\label{criterion of isomorphic simple supercuspidal}
   \begin{split}
    & \text{$\pi_{\chi'}$ and $\pi_{\xi'}$ are isomorphic if and only if }
     \\
     &\text{there exists $t\in T_0\Omega$ such that ${}^t(\chi,\chi')=  (\xi,\xi')$.}
   \end{split}
\end{equation}
 This gives an upper bound 
$$\mathcal{I}\subseteq N(\chi) \subseteq \mathcal{I}\Omega = N_{G(F)}(\mathcal{I}^+)$$
 for the stabilizer $N(\chi)$. Note that in \cite[Sec 8]{Gross-Reeder}, the group $G$ is assumed to be simply-connected, so that $\Omega$ is trivial.

\section{Covering types}
\label{section Covers}

Let $P= MU$ be a parabolic subgroup of a connected reductive group $\mathbf G$ over $F$, where $M$ is a Levi subgroup in $P$ and $U$ the unipotent radical of $P$. Let $P^- = M U^-$ be the opposite of $P$.

A compact subgroup $\mathcal J_P$ of $\mathbf G(F)$ is called {\bf decomposed with respect to} $(M,P)$, or just $(M,P)$-decomposed, if 
$\mathcal{J}_P = \mathcal{J}_P^- \mathcal{J}_M \mathcal{J}_P^+$, where 
$$\mathcal{J}_P^- = \mathcal{J}_P\cap U^-, \quad \mathcal{J}_M = \mathcal{J}_P\cap M,\quad \text{and}\quad \mathcal{J}_P^+ = \mathcal{J}_P\cap U^+.$$
Let $\mathcal{J}_P$ be $(M,P)$-decomposed. We call an element ${\mathbbm z}\in Z_M(F)$ in the center $Z_M$ of $ M$ {\bf strongly positive with respect to }$(P,\mathcal{J}_P)$ if 
$${\mathbbm z} \mathcal{J}_P^+{\mathbbm z}^{-1} \subset \mathcal{J}_P^+, \quad {\mathbbm z}^{-1} \mathcal{J}_P^-{\mathbbm z} \subset \mathcal{J}_P^-,$$
and for any compact subgroups $\mathcal{J}_1,\mathcal{J}_2\subset U^+$ and $\mathcal{J}_3,\mathcal{J}_4\subset U^-$, there exists an integer $N\geq 0$ such that 
$${\mathbbm z}^N \mathcal{J}_1{\mathbbm z}^{-N} \subset \mathcal{J}_2, \quad {\mathbbm z}^{-N} \mathcal{J}_3{\mathbbm z}^N \subset \mathcal{J}_4.$$

An irreducible representation $\lambda_P$ of $\mathcal{J}_P$ is called {\bf decomposed with respect to} $(M,P)$ if $\mathcal{J}_P$ is $(M,P)$-decomposed and both $\mathcal{J}_P^-$ and $\mathcal{J}_P^+$ are contained in $\ker \lambda_P$. We also call $(\mathcal{J}_P,\lambda_P)$ an $(M,P)$-decomposed pair.

\begin{dfn}
\label{definition of cover}
Suppose that $\lambda_M = \lambda_P|_{\mathcal{J}_M}$ is irreducible. We call $(\mathcal{J}_P,\lambda_P)$ a {\bf covering type}, or simply a {\bf cover}, of $(\mathcal{J}_M,\lambda_M)$ in $\mathbf G(F)$ if 
\begin{enumerate}[(a)]
\item $(\mathcal{J}_P,\lambda_P)$ is a decomposed pair,  and 

\label{existence of decomposed pair}

\item  there is a strongly $(P,\mathcal{J}_P)$-positive element ${\mathbbm z}\in Z_M$, and an invertible element in the Hecke algebra $ \mathcal H(\mathbf G(F),\lambda_P)$ supported on $\mathcal{J}_P{\mathbbm z} \mathcal{J}_P$.
\label{existence of invertible supported on positive}
\end{enumerate}
\end{dfn}

We will use the strongly positive element in \ref{existence of invertible supported on positive} above for the computations in Section \ref{section Main calculation}. There is an equivalent criterion in \cite[4.2]{Kim-Yu-tame-types} which is more convenient to confirm directly that a constructed $(\mathcal{J}_P,\lambda_P)$ is a cover of $(\mathcal{J}_M,\lambda_M)$ and will be recalled in the next section. We remark that that criterion is based on the following proposition.

\begin{prop}
\cite{Blondel_Injectivity-of-jacquet} $(\mathcal{J}_P,\lambda_P)$ a cover of $(\mathcal{J}_M,\lambda_M)$ if and only if \ref{existence of decomposed pair} in Definition \ref{definition of cover} and (b') below are satisfied:
\begin{equation*}
   \begin{split}
     \text{(b')}\quad &\text{for any smooth representation $V$ of $\mathbf G(F)$, the Jacquet map $V \rightarrow V_U$} \\
        &\text{induces an injection on the $(\mathcal{J}_P,\lambda_P)$-isotypic subspace of $V$.}
   \end{split}
\end{equation*}
\end{prop}

\subsection{Preliminaries on the construction of covering types}
\label{subsection Preliminaries on the construction of covering types}

Consider the following setup: suppose that $\mathcal{J}_P$ is a compact subgroup in $\mathbf G(F)$ containing a parahoric subgroup $\mathcal{J}_{P,0} $ (whose associating facet is unimportant here and hence ignored), and the pro-unipotent radical $\mathcal{J}_{P,0_+} $ is a normal subgroup of $\mathcal{J}_P$. Assume that all $\mathcal{J}_P, \mathcal{J}_{P,0}$, and $\mathcal{J}_{P,0_+}$ are decomposed with respect to $(M,P)$, and moreover that $\mathcal{J}_P/\mathcal{J}_{P,0_+} \cong \mathcal{J}_M /\mathcal{J}_{M,0_+}  $, where $\mathcal{J}_M = \mathcal{J}_P\cap M $ and $\mathcal{J}_{M,0_+}  = \mathcal{J}_{P,0_+} \cap M$.

Let $(\mathcal{J}_{P,0_+},\theta_P)$ be a character, trivial on both $\mathcal{J}_P^{+ }$ and $\mathcal{J}_P^{- }$. Denote by $(\mathcal{J}_{M,0_+},\theta_M)$ its restriction. Take an irreducible representation $(\mathcal{J}_{M,0},\lambda_M)$ which is $\theta_M$-isotypic, and an $(M,P)$-decomposed pair $(\mathcal{J}_P,\lambda_P)$ with $\lambda_P|_{\mathcal{J}_M}\cong \lambda_M$. 

\begin{prop}
\label{results in Kim-Yu-tame-types}
Under the above setup, we have the following.
\begin{enumerate}[(i)]

\item \cite[Th 6.3]{Kim-Yu-tame-types} The pair $(\mathcal{J}_{P,0_+}, \theta_P)$ is a cover of $(\mathcal{J}_{M,0_+}, \theta_M)$ in $\mathbf G(F)$. 
\label{extension of simple character}

\item \cite[Cor 6.4]{Kim-Yu-tame-types} If $(\mathcal{J}_P,\lambda_P)$ satisfies criterion (\ref{existence of decomposed pair}) of Definition \ref{definition of cover} as well as $\lambda_P|_{\mathcal{J}_{P,+}}$ is $\theta_P$-isotypic, then it is a cover of $(\mathcal{J}_M,\lambda_M)$ in $\mathbf G(F)$. (This is indeed a consequence of (\ref{extension of simple character})).
\label{cover of simple character extends to cover of cuspidal type}

\end{enumerate}

\end{prop}

When constructing covers for classical groups in the next section, we will implicitly identify the notions of lattice filtrations and Moy-Prasad filtrations. Interested readers may consult \cite{Broussous-Stevens-buildings} or \cite{Lemaire-compare-lattice-Moy-Prasad} for details. In particular, we will take a stable functional $\beta$ over $F$ in $M$ and lift it to a functional in $\mathbf G$, such that the building of the centralizer $\mathbf G(F)_\beta = \mathbf G_\beta(F)$ is indeed a tree. We will take two facets  $\mathcal F_y,\mathcal F_z\subset \mathcal B(\mathbf G,F)$ such that each
$$w := \mathcal F_w\cap \mathcal B( \mathbf G_\beta,F) ,\quad w\in \{y,z\},$$ 
becomes a vertex, and moreover $y$ and $z$ are adjacent of each other. Let $\mathcal F_{\mathfrak m} \subset \mathcal B(\mathbf G,F)$ be a facet such that $\mathfrak m:= \mathcal F_{\mathfrak m}\cap \mathcal B( \mathbf G_\beta,F)  $ is the edge connecting $y$ and $z$.

The parahoric subgroups $\mathbf G(F)_{\mathcal F_w,0}$, with $w\in \{y,z\}$, and also $\mathbf G(F)_{\mathcal F_{\mathfrak m},0}$ will be decomposed with respect to a given $(M,P)$. Upon restricting to $ \mathbf G_\beta$, the parahoric subgroup $\mathbf G_\beta(F)_w$ will also be decomposed with respect to $(M\cap \mathbf G_\beta,P\cap \mathbf G_\beta)$, which is also a Levi-parabolic pair in $\mathbf G_\beta$. The indices $\{y,z\}$ will be labelled to match the pair of parabolics $\{P,P^-\}$, in the sense that
 $$\mathbf G_\beta(F)_{\mathfrak m} = (\mathbf G_\beta(F)_y\cap P)\mathbf G_\beta(F)_{y,0+} = (\mathbf G_\beta(F)_z\cap P^-)\mathbf G_\beta(F)_{z,0+}.$$
For both $w\in \{y,z\}$, the quotient $\boldsymbol{\mathsf P}_{\beta,\mathfrak m,w}:= \mathbf G_\beta(F)_{\mathfrak m}/\mathbf G_\beta(F)_{w,0_+}$ will be a maximal parabolic subgroup of $\boldsymbol{\mathsf G}_{\beta,w} :=\mathbf G_\beta(F)_{w}/\mathbf G_\beta(F)_{w,0_+}$.

Later in Section \ref{subsection Lattices and covers for classical groups} when we study reducibilities of certain induced repesentations of classical groups, we will define a character on the pro-p-subgroup $\mathbf G(F)_{\mathcal F_{\mathfrak m},0_+}$. If we can extend this character to an irreducible representation of its normalizer $\mathbf G(F)_{\mathcal F_\mathfrak m}$, then according to the conditions in Proposition \ref{results in Kim-Yu-tame-types} it will be automatically a covering type. In fact, our constructed covering type will also be a character.

\section{Classical groups}
\label{section Classical groups}

This section is basically to translate the general language of Moy-Prasad filtrations and characters into the one involving lattices and strata, after providing the preliminaries of classical groups.

Let $\Fo$ be a non-Archimedean local field with residual characteristic $p\neq 2$. For an extension $F/\Fo$ which is either trivial or quadratic, we denote by $(x\mapsto \bar x)\in \Gal(F/\Fo)$ the trivial or involutive automorphism respectively. We choose representatives $\{1,\zeta,\varpi,\zeta\varpi\}$ of $F^\times/F^{\times2}$, where $\zeta$ is a generator of the group $ \mu_F$ of roots of unity of order coprime to $p$, and $\varpi$ is a uniformizer of $F$ fixed throughout the paper.

Given a finite dimensional vector space $V$ over $F$, we denote by $\tilde{G}=\tilde{G}(V)$ the $F$-algebraic group of linear automorphisms of $V$, which is isomorphic to the general linear group $\GL_{n}$ over $F$, where $n=\dim _FV$. Suppose now $V$ is equipped with a non-degenerate $\epsilon$-Hermitian form $h=h_V$, relative to the extension $F/\Fo$ and $\epsilon=\epsilon_G\in \{\pm 1\}$, which means that
$$h(av,bw)={}\bar ab\cdot h(v,w)\quad\text{ and }\quad
h(w,v)=\epsilon\overline{ h(v,w)},
\quad \text{for all }a,b\in F\text{ and }v,w\in V.$$
The group $G^\sharp = G^\sharp(V,h_V)$ of isometries is the $\Fo$-algebraic group of a classical group, which is either 
\begin{equation*}
\begin{split}
&\text{orthogonal ($\mathrm{O_{odd}}$, 
$\mathrm{O_{even}}$) }
\\
&\text{symplectic ($\mathrm{Sp}$)}
\\
&\text{unitary ($\mathrm{U_{N}}$)}
\end{split}
\qquad\text{when}\qquad
\begin{split}
&\text{$[F/\Fo]=1$, $\epsilon=1$,}
\\
&\text{$[F/\Fo]=1$, $\epsilon=-1$,}
\\
&\text{$[F/\Fo]=2$, $\epsilon=(-1)^{N-1}$}.
\end{split}
\end{equation*}
Denote by $\sigma $ the  involutive automorphism on $\tilde{{G}}$ such that $\tilde{{G}}^\sigma = {{G}}^\sharp$. Let ${G}$ be the connected component of ${G}^\sharp$, so that ${G}={G}^\sharp$ except when ${G}^\sharp$ is an orthogonal group, then ${G}$ is the underlying special orthogonal group ($\mathrm{SO_{odd}}$, $\mathrm{SO_{even}}$).

The involution $\sigma$ on $\tilde G$ induces an adjoint action $\alpha$ on its Lie algebra $\tilde{\mathfrak g}=\tilde{\mathfrak g}(V):=\mathrm{End}_F(V)$, so that $\mathfrak g:=\tilde{\mathfrak g} ^\alpha$ is the Lie algebra of ${G}$.

Suppose now that $(V,h_V)$ is fixed. Let $\tilde V = \tilde V_-$ be another finite dimensional vector space over $F$, and $\tilde V_+$ be its $F$-linear dual. Put  $\tilde V_\pm =\tilde V_-\oplus \tilde V_+$, equipped with the structure of a hyperbolic space with respect to the Hermitian form
\begin{equation*}
h_{\tilde V}:\left<(z-,z+),(w_-,w_+)\right> \mapsto (z_-,w_+) +\epsilon\cdot      {}\overline{(w_-,z_+)},
\quad 
z_-,w_-\in \tilde V_-\text{ and }z_+,w_+\in \tilde V_+.
\end{equation*} 
We define the Hermitian space 
$W := V \perp \tilde V_\pm $, equipped with the form $h_W = h_{\tilde V_\pm } \perp h_V$. The isometry group $\mathbf G^\sharp = G^\sharp(W,h_W)$ is a classical group of the same type as $G^\sharp$ but of a higher rank. We denote the involutive and adjoint automorphisms respectively associated to $\mathbf G^\sharp$ and its Lie algebra $\boldsymbol{\mathfrak g}$ also by $\sigma$ and $\alpha$.

Fix a chain $\tilde V_-\subset V \oplus \tilde V_- \subset W$ and denote by $P$ the corresponding parabolic subgroup, with its unipotent radical by $U$ and the opposite of $U$ by $U^-$. Denote by $M$ the Levi component of $P$, and fix an embedding 
\begin{equation*}
\label{embedding of Levi subgroup}
i_M:\tilde{G}({\tilde V}_-)\times G^\sharp\rightarrow \mathbf G^\sharp, \qquad (g,h)\mapsto (g,h,^\sigma g).
\end{equation*}

For the computations in Section \ref{section Main calculation}, it is more desirable to write down some matrix presentations for elements in $U$ and $U^-$. By choosing a suitable basis, let $H$ be the $\epsilon$-Hermitian matrix, i.e., $H = \epsilon{}^{t}\overline H$, that defines the Hermitian form $h_V$. Then  ${}^\sigma g = H^{-1}{}^t\overline g^{-1} H$ and ${}^\alpha X = -H^{-1}{}^t\overline X H$. Let $\tilde H$ be the matrix that identifies the dual space $\tilde V_+$ with $\tilde V_-$ (relative to implicitly chosen bases), such that the matrix $H_W$ that defines $h_W$ is given by 
\begin{equation}
\label{The big matrix JW}
H_W = \text{anti-diag}(\tilde H, H ,  \epsilon {}^{t}\overline{\tilde H}).\end{equation}
If we denote a typical blocked unipotent matrix by 
\begin{equation*}
(X,Y,Z)^+ := \begin{bmatrix}
I_{\tilde V_-} & X & Y
\\
& I_V& Z
\\
&& I_{\tilde V_+}
\end{bmatrix},
\end{equation*}
then $\alpha:(X,Y,Z)^+ \mapsto ({}^\alpha Z,{}^\alpha Y,{}^\alpha X)^+$, where
\begin{equation}
\label{alpha on X and Y}
{}^\alpha X :=  -H^{-1} {}^{t} \overline X \tilde H,\qquad
{}^\alpha Y := -\epsilon {}^{t}\overline{\tilde H}^{-1}   {}^{t}\overline Y  \tilde H,
\qquad\text{and}\qquad
{}^\alpha Z := -\epsilon {}^{t}\overline{\tilde H}^{-1}  {}^{t} \overline Z  H.
\end{equation}
We have $(X,Y,Z)^+\in U$ if and only if both relations $Z={}^\alpha X$ and \begin{equation}
\label{X-alpha-X-equals-Y-minus-alpha-Y}
X{}^\alpha X = Y-{}^\alpha Y
\end{equation} hold, in which case we simply denote $(X,Y,{}^\alpha X)^+$ by $(X,Y)^+$. Similar results hold for $
(X,Y,Z)^- := \begin{bmatrix}
I_{\tilde V_-} & &
\\
Z& I_V& 
\\
Y&X& I_{\tilde V_+}
\end{bmatrix}\in U^-
$, and we put $(X,Y)^- = (X,Y,{}^\alpha X)^-\in U^-$.

\subsection{Semi-simple strata}

We extract the definitions from \cite{BK, Stevens-supercuspidal} involving lattices and strata. Given an $\mathfrak o(F)$-lattice sequence $\Lambda$ in $V$ and $r\in \mathbb R$, we put  
$$\mathfrak P^{r}(\Lambda):=\{X\in \End_{F}(V):X\Lambda(t)\subseteq \Lambda(t+r)\text{ for all }t\in \mathbb R\}.$$
We emphasize that the filtration of any lattice sequence is normalized, i.e., $\varpi\Lambda(t) = \Lambda(t+1) $ for all $t\in \mathbb R$.

A stratum in an $F$-vector space $V$ is a triple $\mathbf{s}=[\Lambda,t,\beta]$ consisting of an $\mathfrak o_F$-lattice sequence $\Lambda$ in $V$, an element $\beta\in \tilde {\mathfrak g}(V)$, and a real number $t\leq -v_\Lambda(\beta)$. With $\Lambda$ and $t$ fixed, two strata $\mathbf{s}_i=[\Lambda,t,\beta_i]$, with $i=1,2$, are called equivalent if $\beta_1-\beta_2\in \mathfrak P^{-t}(\Lambda)$. We call $\mathbf{s}$ {\bf simple} if $F[\beta]$ is a field and $\Lambda$ is $\mathfrak o_{F[\beta]}$-invariant.

Suppose that $V$ admits a decomposition $\oplus_{i\in I}V_i$, with $\mathbf{1}_i: V\rightarrow V_i$ the projection onto $V_i$ with kernel $\oplus_{j\neq i}V_j$, and $\Lambda_i = \Lambda\cap V_i$ a lattice sequence in $V_i$. A stratum $\mathbf{s}=[\Lambda,t,\beta]$ is now called {\bf semi-simple} if $\beta = \sum_{i\in I}\beta_i$, with each $\beta_i =\mathbf{1}_i\circ \beta\circ \mathbf{1}_i $, each $\mathbf s_i = [\Lambda_i,t,\beta_i]$ is simple or null, i.e., $\beta_i=0$, and $[\Lambda_i\oplus \Lambda_j,t,\beta_i\oplus \beta_j]$ is not equivalent to a simple stratum for all $i,j\in I$ with $i\neq j$.

We can put the definitions above in the self-dual setting. Suppose that $V$ is now a Hermitian space. We call a stratum $\mathbf s = \oplus_{i\in I}\mathbf s_i$ self-dual if $\Lambda$ is self-dual and $\beta\in \tilde {\mathfrak g}^\alpha$, and is called {\bf skew} if moreover the decomposition $\oplus_{i\in I}V^i$ is orthogonal and each $\beta_i\in \tilde {\mathfrak g}(V_i)^\alpha$. Note that $\mathbf s_i$ is then a skew simple stratum for each $i\in I$.

The idea of strata is that, if further $t\geq -v_{\Lambda}(\beta)/2$, the equivalence class of a stratum $\mathbf s$ corresponds to an additive character 
$$\psi_\beta:X\mapsto \psi\circ\tr_{\tilde {\mathfrak g}/F}(\beta X), \quad X\in \mathfrak P^{t_+}(\Lambda),$$
which is trivial on $\mathfrak P^{-v_{\Lambda}(\beta)_+}(\Lambda)$. For all $t\in \mathbb R_{>0}$, denote $\mathcal{U}^{t}(\Lambda): = I+\mathfrak P^{t}(\Lambda)$. We use the same symbol for the character $I+X\mapsto \psi_\beta(X)$ on $\mathcal{U}^{t_+}(\Lambda)$, which is trivial on $\mathcal{U}^{-v_{\Lambda}(\beta)_+}(\Lambda)$. If $\mathbf s$ is furthermore self-dual, then $\psi_\beta$ respectively corresponds to a character on $\mathfrak P^{t_+}(\Lambda)^\alpha$ and on $\mathcal{U}^{t_+}(\Lambda)^\sigma$ by restrictions.

We remark that if in contrast $t<-v_{\Lambda}(\beta)/2$, then we need to go through an approximation process to create the suitable characters for constructing supercuspidal representations, \emph{c.f.} \cite[Sec 2.4]{BK}, \cite[Sec 3.1]{Stevens-semi-simple-char}. Moreover, during the construction we need to extend each character into an irreducible representation of Heisenberg type, which is of dimension $>1$ in general. These processes do not concern us when we later consider epipelagic representations which have shallow depths, and so we will try to keep the related discussions minimal.

\subsection{Stable functionals and epipelagic strata}
\label{subsection Stable functionals for classical groups}

Fix a classical group $G =G(V,h)$ and let $\mathbf{s}=[\Lambda,0,\beta]$ be a skew semi-simple stratum in $V$. By definition, the lattice sequence $\Lambda$ admits an orthogonal decomposition $\Lambda = \oplus_{i\in I}\Lambda_i$, as well as $\beta = \oplus_{i\in I}\beta_i$. Put $E_i = F[\beta_i]$, then each $\Lambda^i$ is $\mathfrak o_{E_i}$-invariant.

Starting from such a stratum $\mathbf s$, we define a subgroup called $H^1(\Lambda,\beta)$ containing $\mathcal{U}^{(-v_\Lambda(\beta)/2)_+}(\Lambda)$ as in \cite[Sec 3.2]{Stevens-semi-simple-char}, and on which a semi-simple character $\theta$ whose restriction on $\mathcal{U}^{(-v_\Lambda(\beta)/2)_+}(\Lambda)$ is $\psi_\beta$. In general $H^1(\Lambda,\beta)$ is a product of groups of the form $\prod_{j}\mathcal{U}^{t_j/2_+}(\Lambda)\cap Z_{\tilde G(F)}({\gamma_j})$, where $\mathbf s_j = [\Lambda,t_j,\gamma_j]$ is a finite sequence of semi-simple strata which approximate $\mathbf s$.

Let's impose a condition on $\mathbf s$ which imitates the epipelagic condition on representations of $G(F)$. To simply put, we want the (positive) depth $-v_{\Lambda_i}(\beta_i)$, for each $\beta_i\neq 0$, to be as small as possible, i.e., if $e_i:=e(\Lambda_i/\mathfrak o_F)$ is the $\mathfrak o_F$-period of $\Lambda_i$, then  $-v_{\Lambda_i}(\beta_i) = 1/e_i$.

 We now state an equivalent condition for $H^1(\Lambda,\beta)$ to have only one factor in its product form, i.e., it is equal to $\mathcal{U}^{(-v_\Lambda(\beta)/2)_+}(\Lambda)$. 
 
 \begin{prop}
 \label{epipelagic and parahoric forcing stability}
Under the above conditions on $\mathbf s$, the equality $H^1(\Lambda,\beta) = \mathcal{U}^{(-v_\Lambda(\beta)/2)_+}(\Lambda)$ holds if and only if $e_i$ are equal for all $i\in I$ with $\beta_i\neq 0$.
 \end{prop}
 \proof
To prove the sufficiency, we just have to recall the definitions from \cite[Sec 3.1 and 3.2]{Stevens-semi-simple-char} and see that $\mathbf s$ has only one step of approximation, namely itself. Hence by definition the product form of $H^1(\Lambda,\beta) $ has only one factor. Conversely, suppose that  we have ordered the indices in $I$ such that $e_i\leq e_{i+1}$ for all $I$, and there exists $i\in I$ with both $\beta_{i}$ and $\beta_{i+1}$ non-zero and $e_i\neq e_{i+1}$, then there are at least two steps of approximations of the form $\gamma_j = \oplus_{k\leq i+1}\beta_k$ and $\gamma_{j+1} = \oplus_{k\leq i}\beta_k$, which render at least two distinct factors in the product form of $H^1(\Lambda,\beta)$. This proves the necessity.
 \qed

Continue from the conditions in Proposition \ref{epipelagic and parahoric forcing stability}, if $x$ is the barycenter of the facet in $\mathcal B(G,F)$ corresponding to $\Lambda$, \emph{c.f.} \cite{Broussous-Stevens-buildings} or \cite{Lemaire-compare-lattice-Moy-Prasad}, then 
$\mathcal{U}^{(-v_\Lambda(\beta)/2)_+}(\Lambda) = G(F)_{x,r}$ where $r = -v_\Lambda(\beta)$. This depth is just $r(x)$ defined in Section \ref{subsection Stable functionals}.

In the remainder of this section, we will show that the definition of stable functionals and the epipelagic condition in \cite{Reeder-Yu} also derive the same forms of semi-simple strata as above.

We first notice that the coset $\beta +\mathfrak P^{(-r)_+}(\Lambda)$ defines a character $\psi_\beta$ on ${\mathsf{V}}_{x,r}$ via the identifications $\mathfrak P^{t}(\Lambda) = \tilde{\mathfrak g}_{x,t}(F)$ and $\mathfrak P^{t}(\Lambda)^\alpha = \mathfrak g_{x,t}(F)^\alpha$ for all $t= -r$ and $(-r)_+$ respectively, and the Moy-Prasad isomorphism $\mathfrak P^{-r}(\Lambda)^\alpha/\mathfrak P^{(-r)_+}(\Lambda)^\alpha\cong {\mathsf{V}}^*_{x,r}(\mathbb F)$.

Let $K_\beta$ be the field extension of $K$ that splits $\beta$, i.e., the field generated by all $\{\beta_i\}_{i\in I}$ over $K$. Put $e = [K_\beta:K]$. Let $\mathfrak g$ be the Lie algebra of $G$, and put ${\mathfrak g}^*_{K_\beta} = {\mathfrak g}^*\times_F K_\beta$. If $x\in \mathcal B(G,F)$ is a rational point of order dividing $e$, then it becomes a hyperspecial vertex in $\mathcal B^{red}(G,K_\beta)$, and the dual space ${\mathsf{V}}^*_{x,r}$ can be regarded as a subspace of $\overline{ {\mathfrak g}^*_{K_{\beta}}} := ({\mathfrak g}^*_{K_{\beta}})_{x,0}/({\mathfrak g}^*_{K_{\beta}})_{x,0_+}$. We now assume that $\beta$ lies in a Cartan subspace, in the sense of \cite[Sec 5]{Reeder-Yu}, of ${\mathsf{V}}^*_{x,r}$, so that we may regard $\beta$ as an semi-simple element in $\overline{ {\mathfrak g}^*_{K_{\beta}}}$. In particular, we have the usual definition of regularity of $\beta$.

The base-changed space  ${\mathfrak g}^*_{K_\beta} $ is equipped with the action of $\Gal(K_\beta/K)$. Let $\mathfrak s_\beta$ be the centralizer of $\beta$ in ${\mathfrak g}^*_{K_{\beta}}$, which is $\Gal(K_\beta/K)$-invariant.

 The following proposition is the first step towards defining what we will call an epipelagic stratum (see Definition \ref{definition of Epipelagic stratum} below). 
 
 \begin{prop}
 \label{equivalent conditions for stability}
(\cite[Lemma 13]{Gross-Levy-Reeder-Yu} or \cite[Prop 3.1]{Kaletha-epipelagic}) Under the above setup, the functional $\beta\in {\mathsf{V}}^*_{x,r}$ is stable if and only if 
\begin{enumerate}[(i)]
\item $\beta$ is regular as an element in $\overline{ {\mathfrak g}^*_{K_{\beta}}}$, and

\label{equivalent conditions for stability - regularity}
\item the $\Gal(K_\beta/K)$-fixed point subspace of $\mathfrak s_\beta$ is in the center of ${\mathfrak g}^*_{K_{\beta}}$.

\label{equivalent conditions for stability - ellipticity}
\end{enumerate}
\end{prop}

From the viewpoint of semi-simple strata, the regularity condition (\ref{equivalent conditions for stability - regularity}) in Proposition \ref{equivalent conditions for stability} of $\beta$ implies that it satisfies the regularity condition of semi-simple strata in \cite[Def. 2.4]{Stevens-supercuspidal}, while condition (\ref{equivalent conditions for stability - ellipticity}) implies that each component $\mathbf s_i:=[\Lambda^i,0,\beta_i]$ is minimal, in the sense of \cite[1.4.15]{BK}. In particular, $\mathbf s$ is a skew semi-simple stratum \cite[Def. 2.5]{Stevens-supercuspidal}.

\begin{cor}
\label{stability implies totally ramified}
Continue with the above setup, if $\beta$ is a stable functional, then each $E_i/F$ is totally ramified.
\end{cor}
\proof Suppose the contrary, i.e., if there exists some residual degree $f_{E_i/F}>1$, then $\beta_i$ is non-regular over $K$ in $K_\beta$. \qed

We can hence translate the tables in \cite[Sec 7.2]{Gross-Levy-Reeder-Yu} (for groups of type A,B,C, or D) into a list of structures of a stable $\beta$.  
  \begin{enumerate}[(i)]
\item For general linear groups and unramified unitary groups, the only allowable order $e$ is the Coxeter number which is just $n = \dim _FV$, and $F[ \beta]/F$ is a totally ramified extension of degree $n$. 
\label{Table GL and unram-U}
 \end{enumerate}
 
For other classical groups, all $E_i/F$ are totally ramified as stated in Corollary \ref{stability implies totally ramified}. There is possibly an index $o\in I$ such that $\beta_o=0$, and is unique by regularity if it exists. There are addition conditions listed as follows.

   \begin{enumerate}[resume*]
\item  For odd orthogonal or symplectic groups,  all extensions $E_i/F$ have a constant even degree for all $i\neq o$. This index does not exist if $G$ is symplectic. 
\label{Table SO-odd}

\item For ramified unitary groups, there exists at most one index $o\in I$ such that $E_o = F$, and all other $E_i/F$ have a constant odd degree. (The oddness condition agrees with the result in \cite[Sec 3.1]{BT-ramified}.)
\label{Table ram-U}

\item  For even orthogonal groups, there exists at most one index $o\in I$ such that $E_o/F$ is ramified quadratic, all other $E_i/F$ have a constant even degree, and furthermore:

   \begin{itemize}
\item if $G$ is split or unramified, then $\#I$ is even;

\item if $G$ is ramified, then $\#I$ is odd.

   \end{itemize}

\label{Table SO-even}

\end{enumerate}

\begin{rmk}
Indeed, the conclusions in Proposition \ref{equivalent conditions for stability} and Corollary \ref{stability implies totally ramified} already imply the shape of $\beta$ in the above table. This implication is clear for general linear groups. The same implication applies to unramified unitary groups, as they are isomorphic to general linear groups over $K$. For other classical groups, the parities of the degrees $e_{E_i/F}$ are already known from the classification in \cite{Stevens-supercuspidal} (and see also \cite{BT-ramified} for ramified unitary groups). For even orthogonal groups, let $d_{E_i/F}$ be the discriminant of $E_i/F$, then we must have $\prod_{i\in I}d_{E_i/F} = \mathrm{disc}(h) \bmod (F^\times)^2$, where $\mathrm{disc}(h)$ is the discriminant of the hermitian form $h$ defining $G$. Represent the possible discriminants, i.e., $F^\times/F^{\times 2}$, by the quartet $\{1,\zeta,\varpi,\zeta\varpi\}$. Since every $E_i/F$ is totally ramified, each $d_{E_i/F} $ can only be $\varpi$ or $\zeta\varpi$. This implies that if $\mathrm{disc}(h)  = 1$ or $\zeta$ (i.e., $G$ is split or unramified orthogonal), then $\#I$ must be even, and if $\mathrm{disc}(h)  = \varpi$ or $\zeta\varpi$ (i.e., $G$ is ramified orthogonal), then $\#I$ must be odd.
 \qed\end{rmk}

We therefore see that the classification of stable functionals in \cite[Sec 7.2]{Gross-Levy-Reeder-Yu} also results in semi-simple strata of the same form as in Proposition \ref{epipelagic and parahoric forcing stability}. Let's summarize this in the following definition.

\begin{dfn}
\label{definition of Epipelagic stratum}
We call a skew semi-simple stratum $\mathbf s = [\Lambda,0,\beta]$ {\bf epipelagic} if  $v_\Lambda(\beta)=-r(x)$ and its components $\beta_i$ satisfy the conditions in the above table \ref{Table GL and unram-U}-\ref{Table SO-even}.
\qed\end{dfn}

As a final remark of this section, the stability of the functional $\beta$ was investigated under the setting in \cite{Gross-Levy-Reeder-Yu}, which requires $F[\beta]/F$ to be tamely ramified. Our definition of epipelagic strata only requires $p$ to be odd, as a particular instance of semi-simple strata in general defined in \cite{Stevens-semi-simple-char}.

\subsection{Epipelagic inducing types}
\label{subsection Epipelagic inducing types for classical groups}

Let's fix a character $(G(F)_{x,r},\psi_\beta)$ arising from a stable functional $\beta$. The stabilizer group $G(F)_{x,0,\beta}$ of $\beta$ in $G(F)_{x,0}$ has quotient $\mathsf A_{x,\beta}: = G(F)_{x,0,\beta} / G(F)_{x,r}$ isomorphic to 
\begin{align*}
& \GL_1(\mathbb F)\cong \mu_F 
 && G = \GL_m,
\\
& \mathrm U_1(\mathbb F/\mathbb F_\bullet)\cong \mu_F^\sigma   
&\text{if}\qquad&\text{$G$ is unramified $\mathrm U_n(F/\Fo)$, or }
\\
&\text{a product of $\mathrm O_1(\mathbb F)\cong \{\pm 1\}$} 
&&\text{$G$ is of other types of classical groups}.
\end{align*}
Since all these quotients are abelian and have order coprime to $p$, we can extend $\psi_\beta$ to a character of $G(F)_{x,0,\beta}$, and we have a bijection
$$\{\lambda\in [G(F)_{x,0,\beta}]^\wedge\text{ where }\lambda|_{G(F)_{x,r}}=\psi_\beta\}\xrightarrow{\sim } \hat {\mathsf A}_{x,\beta}.$$
In the last case above, the precise numbers of components can be read from the tables in \cite[Sec 7.2]{Gross-Levy-Reeder-Yu}, with extra considerations on the simply-connectedness of the groups. The number of components is $\#I-2$ for odd orthogonal groups (note that the index $o\in I$ does not contribute to a $\{\pm 1\}$-factor), $\#I-1$ for even split or unramified $\SO$ and even ramified unitary groups, and $\# I$ for other types of classical groups. The element
\begin{equation}
\label{definition of omega_i element}
\omega_i = \diag((I_j)_{i\in I\smallsetminus i}, -I_i)\in \left(\textstyle\prod_{i\in I}G(V_i)(F)\right)\cap G(F)_{x,0}\subset G(F)_{x,0,\beta}
\end{equation}
corresponds to $((1)_{j\neq i},-1_i)\in \mathsf A_{x,\beta}$, and an extension $\lambda$ of a fixed $\psi_\beta$ is determined by $\{\lambda(\omega_i)\}_{i\in I}$.

In the last two cases, the full stabilizer group ${\mathcal{J}} = G(F)_{x,\beta}$ is isomorphic to $G(F)_{x,0,\beta}$ for unramified unitary groups, to $\{\pm 1\}^{\# I-1}$ for odd SO, and to $\{\pm 1\}^{\# I}$ for other types of classical groups. Note that $\mathcal J$ is the normalizer denoted by $N(\psi_\beta)$ in the Intro.. Again $\lambda$ can be extended to ${\mathcal{J}}$, so that $\pi_\lambda = \cInd_{\mathcal J}^{G(F)}\lambda$ is a supercuspidal representation, and is epipelagic.

In the first case $G = \GL_m$, the full stabilizer group $G(F)_{x,\beta}$ is $\tilde{\boldsymbol{\mathcal{J}}} = F[\beta]^\times G(F)_{x,r}$ and contains $\tilde{{\mathcal{J}}} = G(F)_{x,0,\beta}$ as the maximal compact subgroup. Note that $\tilde{\boldsymbol{\mathcal{J}}}$ is the normalizer denoted by $\tilde N(\psi_\beta)$ in the Intro.. By \cite[6.2.2]{BK}, $\tilde{\boldsymbol{\mathcal{J}}} $ is indeed equal to the intertwining set $I_{\tilde{G}(F)}(\tilde{\lambda})$. As in \cite[6.1.2]{BK}, we can extend $\tilde\lambda$ to a character $\tilde{\boldsymbol{\lambda}}$, so that $\tilde\pi_{\tilde{\boldsymbol{\lambda}}} = \cInd_{\tilde{\boldsymbol{\mathcal{J}}}}^{\tilde G(F)}\tilde{\boldsymbol{\lambda}}$ is a supercuspidal representation, and is epipelagic.

Continue from above, if $\beta$ is moreover self-dual, then we can extend $\psi_\beta$ to a self-dual character $(\tilde{\mathcal{J}}, \tilde\lambda)$ by requiring $\tilde\lambda|_{\mu_F}$ to be or order $\leq 2$. Since $F[\beta]/F$ is totally ramified, we can choose a uniformizer $\varpi_\beta$ of $F[\beta]$ such that ${}^\sigma\varpi_\beta = -\varpi_\beta^{-1}$. This implies that we can further extend $\tilde\lambda$ to a character $(G(F)_{x,\beta},\tilde{\boldsymbol{\lambda}})$ such that $\tilde{\boldsymbol{\lambda}}(\varpi_\beta)^2  = \tilde\lambda(-1)$.

 \section{Reducibility}
\label{section Reducibility}

Recall the setup at the beginning of Section \ref{section Classical groups}: let $G$ be a classical group defined by a Hermitian space of dimension $n$, and $\tilde G$ be a general linear group $\GL_m$ for some $m$. We denote by $\mathbf G $ a classical group, of the same type as $G$ but of a higher rank, and containing $M=\tilde G\times G$ as a Levi subgroup.

 Let $\pi_M = \tilde\pi\times \pi$ be an irreducible supercuspidal representation of $M$. Take a parabolic subgroup $P$ containing $M$. We are interested in the points $s\in \mathbb {C}$ where the normalized  parabolic induction 
 \begin{equation}
\label{reducibility of parabolic induction with complex point}
I(s,\tilde\pi,\pi):=\iota_{P}^{\mathbf G}(\tilde\pi|\det|^s\times \pi)
\end{equation}
is reducible. Due to the obvious reason, we confine ourselves in the domain 
\begin{equation}
\label{complex domain}
\{s\in\mathbb C: 0 \leq\Im(s) < \tfrac{2\pi }{f_{\tilde\pi} \log q}\},
\end{equation}
where $f_{\tilde\pi}$ is the order of the subgroups of unramified characters $\chi$ of $F^\times$ such that $\chi\tilde\pi\cong \tilde\pi$. Since we only concerns epipelagic representations in this paper, the field datum of the underlying stratum of $\tilde\pi$ is totally ramified over $F$, and so $f_{\tilde\pi}=1$.

If $I(s,\tilde\pi,\pi)$ is reducible for some $s\in \mathbb C$, then up to the twisting by an unramified character, we  assume that $\tilde{\pi}$ is self-dual. There is then a unique \emph{real} $s_{\tilde{\pi},{\pi}}\geq 0$, indeed a half integer \cite[Sec 4]{Moeglin-classification-classical-groups},\cite[Sec 3]{Moeglin-classification-unitary-groups}, such that $I(s,\tilde\pi,\pi)$ is reducible at $s=\pm s_{\tilde{\pi},{\pi}}$. By  \cite[6.2.5]{BK}, there are exactly two self-dual representations in the inertial class of $\tilde{\pi}$, namely $\tilde\pi $ and its twist  $\tilde\pi' $ by the unramified character $|\det|^{ \frac{\pi \sqrt{-1}}{\log q}}$. The {complex} points of reducibility of $I(s,\tilde\pi,\pi)$ are therefore of the form \begin{equation*}
\label{4 reducibility points}
\mathrm{Red}(\tilde{\pi},{\pi}) := \left\{\pm s_1,\,\pm s_2+\frac{\pi \sqrt{-1}}{\log q}\right\}\qquad\text{ for some }s_1,\,s_2\in \frac{1}{2}\mathbb{Z}_{\geq 0},
\end{equation*}
and those in $\mathrm{Red}(\tilde{\pi}',{\pi})$ take   the same form, with $s_1$ and $s_2$ exchanged.

 The reason we are interested in the reducibility of (\ref{reducibility of parabolic induction with complex point}) is due to a theory of M{\oe}glin \cite{Moeglin-classification-classical-groups,Moeglin-classification-unitary-groups,Moeglin-twisted-endoscopy-Langlands-parameters}; summarizing it in one sentence, it asserts that if $s_{\tilde{\pi},{\pi}}\geq 1$, then \begin{equation*}
\text{ the Langlands parameter of $\tilde\pi$ is a component of the Langlands parameter of $\pi$.}
\end{equation*}
To be more precise, let $\mathcal E (\pi)$ be the set 
 $$\left\{(\tilde G,\tilde\pi): 
     \begin{matrix}   \text{$\tilde G$ a general linear group, $\tilde\pi$ an (isomorphism class of)}
     \\
     \text{irreducible supercuspidal representation of $\tilde G(F)$,}
     \\
    \text{such that  (\ref{reducibility of parabolic induction with complex point}) is reducible for some $s\in  \tfrac{1}{2}\mathbb N$ with $s\geq 1$}     \end{matrix}
\right\}.$$
Suppose that $s=1$ for all $(\tilde G,\tilde\pi)\in \mathcal E(\pi)$, which is the only case we concern for epipelagic representations. (We refer to the general situation to \cite[Sec 4]{Moeglin-twisted-endoscopy-Langlands-parameters}, in which certain parity conditions on $2s-1$ are required.) This set $\mathcal E(\pi)$ is the extended cuspidal support of $\pi$, in the sense that the endoscopic lifting of $\pi$ is the parabolically induced representation 
$$\text{$\Pi:=\prod_{\tilde\pi\in \mathcal E(\pi)}\tilde\pi$ of $\GL_{\hat N}(F)$,}
\quad
\text{where 
 $\hat N = \sum_{(\GL_m(F),\tilde\pi)\in \mathcal E(\pi)}m$.}$$ 
 (In fact, $\hat N$ is known to be equal to the rank of the matrices in the Langlands dual group of $G$, and is equal to $N$ or $N\pm 1$ depending on the type of $G$.) The lifting $\Pi$ is $\sigma$-elliptic in the sense of \cite[Sec 2.8]{Moeglin-twisted-endoscopy-Langlands-parameters}, and is also irreducible. Hence understanding the set $\mathcal E(\pi)$ suffices for understanding the lifting $\Pi$ of $\pi$.

 \subsection{The Hecke algebra approach}
 
 The harmonic analysis for the parbolically induced representation (\ref{reducibility of parabolic induction with complex point}), especially concerning character expansions, could be eventually difficult. One algebraic approach, stemming from \cite{BK-types}, is to view this representation in its Bernstein component of the category of smooth representations of $\mathbf G$ and relate this component with the category of the associated Hecke algebra, which is explained below.

 Let $\mathfrak s_M$ be the inertial class of $\pi_M := \tilde\pi\times \pi$. Suppose that $\pi = \cInd_{\mathcal J}^{G(F)}\lambda$ and $\tilde\pi = \cInd_{\tilde{\boldsymbol{\mathcal{J}}}}^{\tilde G(F)}\tilde{\boldsymbol{\lambda}}$ are the epipelagic representations constructed as in Section \ref{subsection Epipelagic inducing types for classical groups}, with $\tilde{\boldsymbol\lambda}$ being extended from a character $\tilde\lambda$ of the maximal compact subgroup $\tilde {\mathcal J}$ of $\tilde{\boldsymbol{\mathcal{J}}}$. We put $\mathcal J_M = \tilde{\mathcal{J}}\times {\mathcal J}$ and $\lambda_M = \tilde\lambda\times \lambda$, and let $\mathcal{H}(M,\lambda_M)$ be the associated Hecke algebra \cite[Sec 4.1]{BK}. Since $\pi_M$ is supercuspidal, it is known \cite{BK-types} that $({\mathcal J}_M,\lambda_M)$ is an $\mathfrak{s}_M$-type, which means that there is an equivalence of categories
\begin{equation*}
\mathcal{M}_M:\mathcal{R}^{\mathfrak{s}_M}(M)\rightarrow \text{Mod-}\mathcal{H}(M,\lambda_M),\,\qquad\tau\mapsto \Hom_{{\mathcal J}_M}(\lambda_M,\tau).
\end{equation*}
If $({\mathcal J}_P,\lambda_P)$ is a cover of $({\mathcal J}_M,\lambda_M)$ in $\mathbf G(F)$, from \cite{BK-types} we know that $({\mathcal{J}}_P,\lambda_P)$ is an $\mathfrak s$-type, where $\mathfrak s = \Ind_M^{\mathbf G}\mathfrak s_M$, and there is an analogous equivalence of categories
$
\mathcal{M}_{\mathbf G}:\mathcal{R}^{\mathfrak{s}}(\mathbf G)\rightarrow \text{Mod-}\mathcal{H}(\mathbf G,\lambda_P)
$. 
Let
\begin{equation}
\label{Morphism t_P}
t_P:\mathcal{H}(M,\lambda_M)\rightarrow
\mathcal{H}(G_W,\lambda_P)
\end{equation}
be the injective morphism of algebras defined in \cite[(8.3, 8.4)]{BK-types}, and denote by $(t_P)_*:\text{Mod-}\mathcal{H}(M,\lambda_M)\rightarrow \text{Mod-}\mathcal{H}(\mathbf G,\lambda_P)$ the co-induced functor between the module categories. The above constructions render the following commutative diagram:
\begin{equation}\label{commutative diagram}
  \xymatrixcolsep{5pc}\xymatrix{
\mathcal{R}^{\mathfrak s}(\mathbf G)
\ar[r]^{\mathcal{M}_{\mathbf G}}    
&\text{Mod-}\mathcal{H}(\mathbf G,\lambda_P)
\\
\mathcal{R}^{\mathfrak s_M}(M) \ar[u]_{\iota_P^{\mathbf G}}
\ar[r]^{\mathcal{M}_M }
&\text{Mod-}\mathcal{H}(M,\lambda_M).  
\ar[u]_{(t_P)_*}
}
\end{equation}
Hence we can study the reducibility of (\ref{reducibility of parabolic induction with complex point}) by understanding the structures of $\mathcal{H}(M,\lambda_M)$ and $
\mathcal{H}(\mathbf G,\lambda_P)$, as well as their module categories.

\subsection{Structures of Hecke algebras}
\label{subsection Structures of Hecke algebras}

 We refer to \cite[Ch 5]{BK} and \cite[Sec 6 and 7]{Stevens-supercuspidal} the fine structures of the inducing types of supercuspidal representations of general linear groups and classical groups.

Continue with the setup from the previous section. Recall that $\tilde\lambda$ is constructed from the epipelagic simple stratum 
$[\tilde\Lambda,0,\tilde\beta]$. We decompose $\tilde \lambda = \tilde \rho \otimes \tilde \kappa$, where $\tilde \kappa$ has order a power of $p$, and $\tilde \rho $ has order coprime to $p$. Hence $\tilde \rho $ is a depth zero character, i.e., a character of $\mathcal{U}(\tilde\Lambda)_{F[\tilde\beta]} $ trivial on $\mathcal{U}^{0_+}(\tilde\Lambda)_{F[\tilde\beta]}$. Similarly, $\lambda$ is constructed from the epipelagic semi-simple stratum 
$[\Lambda,0,\beta]$, and we decompose $ \lambda =  \rho \otimes  \kappa$, where $\kappa$ has order a power of $p$ and $ \rho $ is a depth zero character of $\mathcal{U}(\Lambda)_{F[\beta]}^\sigma $ trivial on $\mathcal{U}^{0_+}(\Lambda)_{F[\beta]}^\sigma$. The character $\tilde{\rho}\times \rho$ hence descends to 
\begin{equation}
\label{product of finite reductive groups}
\mathcal{U}(\tilde\Lambda)_{F[\tilde\beta]} /\mathcal{U}^{0_+}(\tilde\Lambda)_{F[\tilde\beta]} \times \mathcal{U}(\Lambda)_{F[\beta]}^\sigma /\mathcal{U}^{0_+}(\Lambda)_{F[\beta]}^\sigma,
\end{equation}
 which is a product of finite groups, each of whose factors is the subgroup of rational points of a reductive (but not necessarily connected) group over the finite field $\mathbb F_{F[\tilde \beta]}$ and $\mathbb F_{F[\beta_i]}$ with $i\in I$ respectively, or all over the base finite field $\mathbb F$ by restriction of scalars.

We now describe the structures of the related Hecke algebras. Recall that $\tilde{\boldsymbol{\mathcal J}}$ is equal to the intertwining set $I_{\tilde{G}(F)}(\tilde{\lambda})$. We have $\tilde{\boldsymbol{\mathcal J}} = F[\tilde\beta]^\times\tilde{\mathcal J}$, and so $
\tilde{\boldsymbol{\mathcal J}}/\tilde{\mathcal J} \cong F[\tilde\beta]^\times/\mathfrak o^\times_{F[\tilde\beta]} \cong \left<\varpi_{\tilde \beta}\right>$, where $\varpi_{\tilde \beta}$ is a uniformizer of $F[\tilde\beta]$. The Hecke algebra $\mathcal{H}(\tilde{G},\tilde{\lambda})$ is isomorphic to $\mathbb{C}[Z,Z^{-1}]$, where $Z$ is supported on the single coset $\varpi_{\tilde\beta} \tilde{\mathcal J}$. As for $\lambda$, we have $I_{{G}(F)}({\lambda})={\mathcal J}$, and so $\mathcal{H}({G},{\lambda})\cong \mathbb{C}$. Therefore, $\mathcal{H}(M,\lambda_M)\cong \mathbb{C}[Z,Z^{-1}]$.

We now describe the structure of $\mathcal{H}(\mathbf G,\lambda_P)$. If $({\mathcal J}_P,\lambda_P)$ is a cover of $({\mathcal J}_M,\lambda_M)$ but $\tilde{\lambda}$ is not self-dual, then 
$\mathcal{H}(\mathbf G,\lambda_P)\cong \mathcal{H}(M,{\lambda}_M)\cong \mathbb{C}[Z,Z^{-1}]$. This case is not of our interest. In the interesting contrast, when $\tilde{\lambda}$ is self-dual, it is known \cite[Cor 6.16]{Stevens-supercuspidal},  \cite[Prop 3.3]{Blondel-Weil} that     
\begin{equation*}
\mathrm{rank}_{\mathcal{H}(M,\lambda_M)}
(\mathcal{H}(\mathbf G,\lambda_P))=
\#(N_{\mathbf G(F)}(\mathfrak s_M)/M(F))=2,
\end{equation*}
so that the Hecke algebra $\mathcal{H}(\mathbf G,\lambda_P)$ is a rank-2 module over $ \mathcal{H}(M,{\lambda}_M)$ \cite[(11.4)]{BK-types} via the injective morphism $t_P$ in (\ref{Morphism t_P}). We may describe this algebra by choosing two special generators $T_w$, for $w\in \{y,z\}$, as follows. First pick two elements $s_y$ and $s_z$ in $\mathbf G(F)$ (for example, we may choose $s_1$ and $s^\varpi_1$ in \cite{Stevens-supercuspidal} or \cite{Blondel-Weil}), each of which is a generator for the normalizer group $N_{\mathbf G(F)}(\mathfrak s_M)$ mod $M(F)$, satisfying the following conditions.
\begin{enumerate}[(i)]
\item $s_y{\mathcal J}_P^-s_y^{-1}\subset {\mathcal J}_P^+$\text{ and }$s_z{\mathcal J}_P^+s_z^{-1}\subset {\mathcal J}_P^-$.

\item If we put $\tilde{\mathbbm z} :=s_ys_z=i_M(\varpi_{\tilde\beta}I_{\tilde V_-},I_V)$, then ${\mathcal J}_P \tilde{\mathbbm z}^{e(F[\tilde\beta]/F)}{\mathcal J}_P = {\mathcal J}_P {\mathbbm z}{\mathcal J}_P  $ for a strongly $(P,{\mathcal J}_P)$-positive element ${\mathbbm z}$ in the center of $M$.

\item The generator $T_w$, for each $w\in \{y,z\}$, is supported on the double coset ${\mathcal J}_Ps_w{\mathcal J}_P$.

\item  Each $T_w$ satisfies a quadratic relation of the form
\begin{equation}
\label{quadratic relation general b and c}
T_w*T_w=b_wT_w+c_w     \mathbbm{1}  
\end{equation}
for certain real numbers $b_w$ and $c_w$ (here $\mathbbm{1}$ is the unit in $\mathcal{H}(\mathbf  G,\lambda_P)$, which is the function supported  on ${\mathcal J}_P$ with $\mathbbm{1}(1)=I_{\lambda_P}$ the identity operator on the representation space of $\lambda_P$).

\item The relation $T_y*T_z=t_P(Z)$ holds.\end{enumerate}
We have in particular $t_P(Z) ( {\mathbbm z}) =  T_y(s_y) \circ T_z (s_z) $. Moreover, $t_P(Z) ( {\mathbbm z})^{* e(F[\tilde\beta]/F)}$ is supported on the strongly positive element ${\mathbbm z}$, endorsing our assumption that $({\mathcal J}_P,\lambda_P)$ is a cover of $({\mathcal J}_M,\lambda_M)$ (see Definition \ref{definition of cover}).

There are two directions to study the real coefficients $b_w$ and $c_w$. 
In the one direction, we follow \cite[Sec 1]{Blasco-Blondel-SP4} to compute these coefficients directly: $c_y =[{\mathcal J}_P^+:s_y{\mathcal J}_P^- s_y^{-1}]$, 
\begin{equation}
\label{formula for b and c for y}
b_y=\sum_{u\in \mathcal S_y} {T_y(u)}, \quad\text{where \quad $\mathcal S_y:= { \frac{s_y{\mathcal J}_P^+s_y^{-1}\cap {\mathcal J}_Ps_y{\mathcal J}_P}{{\mathcal J}_P^- }}$},
\end{equation}
and similarly $c_z  =[s_z {\mathcal J}_P^-s_z^{-1}:{\mathcal J}_P^+] $,
\begin{equation}
\label{formula for b and c for z}
b_z =\sum_{u\in \mathcal S_z} { T_z(u)}{},\quad 
\text{
where \quad $ \mathcal S_z:=\frac{s_z^{-1}{\mathcal J}_P^-s_z\cap {\mathcal J}_Ps_z{\mathcal J}_P}{{\mathcal J}_P^+}$}.
\end{equation}
Under the condition $\dim\lambda_P=1$, the summands in $b_w$ are just scalars. (In general, $b_w$ is a sum of traces of operators in $\End_{\mathbb C}( V_{\lambda_P})$.) We have implicitly chosen $s_w$, for $w\in \{y,z\}$, such that $s_w^2\in {\mathcal J}_M$, so that we may and do normalize each $T_w$, up to a sign,
 such that $
T_w(s_w)^2=\lambda_P(s_w^2)$, giving the positive values for both $c_w$. The choices of $s_w$ will become clear in (\ref{representative-sy-and-sz}).


In the another direction, we view $[\Lambda_W,0,\beta_W]:=[\tilde\Lambda\oplus \Lambda,0,\tilde\beta \oplus \beta]$ as a stratum in $\mathbf G(F)$, and $\tilde\Lambda\oplus \Lambda$ as an $(\mathfrak o_{\tilde\beta}\oplus \mathfrak o_{\beta})$-lattice sequence. There are two minimal lattice sequences properly contained in $\tilde\Lambda\oplus \Lambda$, denoted by $\mathfrak M^w$ with $w\in \{y,z\}$. Their precise definitions will be given in Section \ref{subsection Lattices and covers for classical groups}; in fact, they correspond to the edge $\mathfrak m$ and the vertices $\{y,z\}$ of the building $\mathcal B(\mathbf G_\beta,F)$ mentioned in Section \ref{subsection Preliminaries on the construction of covering types}. Put 
$$\boldsymbol{\mathsf{G}}_{\beta,w} = \mathbf G_\beta(F)_{w}/\mathbf G_\beta(F)_{w,0_+} = \mathcal{U}(\mathfrak{M}^w)_{F[\beta]} / \mathcal{U}^{0_+}(\mathfrak{M}^w)_{F[\beta]},$$ which is the subgroup of rational points of a product of reductive (but not necessarily connected) groups, over $\mathbb F_{{F[\tilde\beta]} }$ and  $\mathbb F_{{F[\beta]} }$ respectively, and containing the product group (\ref{product of finite reductive groups}) supporting the character $\tilde\rho\times\rho$. By \cite[Section 7.1]{Stevens-supercuspidal}, there exist injections 
$$\mathcal{H}(\boldsymbol{\mathsf{G}}_{\beta,w}  ,\tilde{\rho}\times \rho) \cong \mathcal{H}(\mathcal{U}(\mathfrak{M}^w)_{F[\beta]},\tilde{\rho}\times \rho)\hookrightarrow  \mathcal{H}(\mathbf G,\lambda_P), \quad w\in \{y,z\},$$ 
with both images together generate the whole algebra $ \mathcal{H}(\mathbf G,\lambda_P)$. Hence we reduce to considering Hecke algebras for finite reductive groups. The structure of these algebras is well-known \cite{Lusztig-finite-classical-groups}: under suitable normalizations of $T_w$, that we will denote by $\mathcal T_w$  to avoid confusions, the quadratic relation can be written as 
\begin{equation}
\label{normalized quadratic relation general b and c}
(\mathcal T_w+\mathbbm{1})*(\mathcal T_w -q_{F[\tilde\beta]}^{r_w}\mathbbm{1})=0
\end{equation}
for certain integers $r_w\geq 0$. The values of $r_w$ can be read from \cite[Table II]{Lusztig-Chevalley-groups}. In our present situation, the extension $F[\tilde\beta]/F$ is totally ramified, and so $q_{F[\tilde\beta]}= q$.

Therefore, using (\ref{formula for b and c for y}) and (\ref{formula for b and c for z}), and comparing (\ref{quadratic relation general b and c}) and (\ref{normalized quadratic relation general b and c}), we can determine the value of $b_w$, for $w\in \{y,z\}$, up to a sign. We will analyze this sign in the next section, and relate it with the reducibility of $I(s,\tilde\pi,\pi)$ as well.

\subsection{Eigenvalues of Hecke algebra  elements}


For $s\in \mathbb C$, let $D_s = \mathcal M_M(\tilde\lambda|\det|^s\times \lambda)$ and  $X_s = \mathcal {M}_{\mathbf G}(I(s,\tilde\pi,\pi))$ be respectively the $\mathcal  H(M,\lambda_M)$- and $\mathcal  H(\mathbf G,\lambda_P)$-modules under the categorial equivalences in the diagram (\ref{commutative diagram}). Note that $D_s$ is just a line over $\mathbb C$, and $X_s$ is two dimensional over $\mathbb C$.

The following proposition is a crucial relation between the reducibility of parabolically induced representations and the eigenvalues for modules of Hecke algebras.
\begin{prop} 
\label{Blondel-Blasco-reducibility-criterion} \cite[(1.13)]{Blasco-Blondel-SP4}
The module $X_s$ is reducible if and only if the product of eigenvalues of $T_y$ and of $T_z$ on $X_s$ is equal to the scalar action of $Z$ on $D_s$.
\end{prop}

We continue to assume that $\dim \lambda_P=1$. On the one hand, since $\tilde\pi$ is induced from a self-dual extension $(\tilde{\boldsymbol{\mathcal J} }, \tilde{\boldsymbol\lambda})$ of $( \tilde{\mathcal  J}, \tilde{\lambda})$, the scalar action of $Z$ on $D_s$ is \cite[(1.15)]{Blasco-Blondel-SP4}
$$
q^{s} \Delta_P({\mathbbm z})^{1/2} {\tilde{\boldsymbol{\lambda}}}(\varpi_{\tilde\beta})^{-1} T_y(s_y)T_z(s_z).$$
Here $\Delta_P$ is the modular character of $P$ coming from normalizing the parabolic induction, and we have indeed
 \begin{equation*}
 \Delta_P({\mathbbm z})=[{\mathcal J}_P^+: {\mathbbm z} {\mathcal J}_P^+{\mathbbm z}^{-1}]=[s_y {\mathcal J}_P^+s_y:{\mathcal J}_P^-][s_z {\mathcal J}_P^-s_z^{-1}:{\mathcal J}_P^+] = c_yc_z. 
 \end{equation*}
On the other hand, by comparing (\ref{quadratic relation general b and c}) and (\ref{normalized quadratic relation general b and c}), the eigenvalues of $T_w$ are 
\begin{equation*}
     \epsilon_w c_w^{1/2} q^{r_w/2  }
     \qquad\text{and}\qquad
      -\epsilon_w  c_w^{1/2} q^{-r_w/2  }
\end{equation*}
and the possible products of eigenvalues of $T_y$ and   $T_z$ are  
\begin{equation*}
\label{product of eigenvalues of Ty and Tz}
   \epsilon_y \epsilon_z  (c_y c_z)^{1/2}  q^{\pm (r_y+r_z)/2  }
        \qquad\text{and}\qquad
         -\epsilon_y \epsilon_z  (c_yc_z)^{1/2}  q^{\pm (r_y-r_z)/2  }.
\end{equation*}
When $X_s$, and hence $I(s,\tilde{\pi},\pi)$, is reducible,  Proposition \ref{Blondel-Blasco-reducibility-criterion} implies that
 \begin{equation*}
\label{values of real parts}
   \Re \mathrm{Red}(\tilde{\pi},{\pi}) = \{\pm \frac{r_y+r_z}{2},\,\pm \frac{r_y-r_z}{2}\}.
\end{equation*}
We will compute the signs $\epsilon_w$, for $w\in \{y,z\}$, in Section \ref{section Main calculation}. The result gives the exact four values of $s\in \mathbb C$ (possibly with multiplicities) in the domain (\ref{complex domain}) at which $I(s,\tilde\pi,\pi)$ is reducible: 
\begin{equation}
\label{main formula of reducibility points}
q^s = \delta\epsilon_y T_y(s_y)\epsilon_z T_z(s_z)^{-1} {\tilde{\boldsymbol{\lambda}}}(\varpi_{\tilde\beta})  q^{\epsilon (r_y+\delta r_z)/2},\quad \epsilon, \delta \in \{\pm 1\}.
\end{equation}
When $F[\tilde\beta]/F$ is a maximal subfield in $\tilde {\mathfrak g}(\tilde V)$ and is totally ramified, $\dim \tilde\lambda=1$, and so the intertwining operators $T_w(s_w)$, for $w\in \{y,z\}$, are just scalars. 

\begin{prop}
\label{main prop of reducibility points}
In the situation above, if ${\tilde{\boldsymbol{\lambda}}}(\varpi_{\tilde\beta}) = \delta\epsilon_y T_y(s_y)\epsilon_z T_z(s_z)$ for some $\delta \in \{\pm 1\}$, then 
$$\mathrm{Red}(\tilde\pi,\pi) = \left\{\pm \frac{r_y+\delta r_z}{2},\pm \frac{r_y-\delta r_z}{2}+\frac{\pi \sqrt{-1}}{\log q} \right\}.$$
\end{prop}

\subsection{Examples of liftings of characters}
\label{Example: depth zero characters}

We provide two simple examples to demonstrate our methodology explained in this section.
The results will be also useful in the calculation in Section \ref{subsection Expanding the intertwining operator as a sum}. Readers who are unfamiliar with the idea behind may first 
consult with the next section.

\subsubsection{Example 1: U(1)}

Consider $\tilde G = \GL_1$ and $G = \mathrm U_{1,F/\Fo}$. Take characters $(F^\times, {\tilde{\boldsymbol{\lambda}}})$ and $(\mathrm U_1(F/\Fo), {\lambda})$, and put ${\tilde{{\lambda}}} = 
{\tilde{\boldsymbol{\lambda}}}|_{\mathfrak o^\times_F}$. Put $\mathcal{U}^{0_+}(F)_1 = \mathcal{U}^{0_+}(F) \cap \mathrm U_1(F/\Fo)$ and $\mu(F)_1 = \mu_F \cap \mathrm U_1(F/\Fo)$.

Denote by $(1-c)$ the map $x\mapsto x \bar x^{-1}:F^\times\rightarrow \mathrm U_1(F/\Fo)$. Assume that $\tilde{\lambda}|_{\mathcal{U}^{0_+}(F)} = {\lambda}|_{\mathcal{U}^{0_+}(F)_1} \circ(1-c)$, and abbreviate $\tilde{\lambda} \leftrightarrow {\lambda}$ when the following situation occurs: 
$$\tilde{\lambda}|_{\mu_F } = {\lambda}|_{\mu(F)_1} \circ(1-c).$$
First consider ${F}/{\Fo}$ being unramified. The self-dual condition on $\tilde{\lambda}$ implies that $\tilde{\lambda}|\mu_{{\Fo}}\equiv 1$. If $\zeta\in \mu_{{F}}\smallsetminus\mu_{{\Fo}}^\times$ but $\zeta^2\in \mu_{\Fo}$, then $\tilde{\lambda}(\zeta)  = \lambda(-1) \in \{ \pm 1\}$ when $\tilde\lambda\leftrightarrow \lambda$. We also pinpoint a character $\boldsymbol{\tilde{\lambda}}$  by taking $\boldsymbol{\tilde{\lambda}}(\varpi)=1$, and denote by $\boldsymbol{\tilde{\lambda}}'$ the self-dual unramified twist of $\boldsymbol{\tilde{\lambda}}$, i.e., $\boldsymbol{\tilde{\lambda}}(\varpi)=-1$.

We use (\ref{formula for b and c for y}) to calculate  $c_y = q_{}^{3/2}$ and
$$b_yT_y(s_y) = \sum_{\begin{smallmatrix}
X\in \mathfrak{o}_{F}/\mathfrak{p}_{F}
\\
Y \in \mu_{F}
\\
Y+{}\overline Y = -X\overline X
\end{smallmatrix}}\tilde{\lambda}(Y)\lambda(1+\overline XY^{-1} X) 
 = \begin{cases}
 \lambda(-1) (\qo^3 -1) & \text{if }\tilde{\lambda}\leftrightarrow \lambda,
 \\
 - \lambda(-1) \qo(\qo-1) & \text{if }\tilde{\lambda}\not\leftrightarrow \lambda \text{ and $\tilde{\lambda}$ is self-dual}.
 \end{cases} $$
The calculation for $z$ is easier; using (\ref{formula for b and c for z}), we have
$c_z = q_{}^{1/2} $ 
and 
$$b_z T_z(s_z) = \sum_{\begin{smallmatrix}
Y \in \mu_{F}
\\
Y= -{}\overline Y \end{smallmatrix}}\tilde{\lambda}(Y)
 =  \lambda(-1) (\qo-1) \quad \text{ if $\tilde{\lambda}$ is self-dual}. $$
 As a remark, one can compute that $b_wT_w(s_w)=0$ for both $w\in \{y,z\}$ when $\tilde\lambda$ is not self-dual, but we do not need this result later on.

 When $\tilde{\lambda}$ is self-dual, we have $r_z = 1/2$, and $r_y = 3/2$ or $ 1/2$ depending on whether $\tilde{\lambda}\leftrightarrow \lambda$ or not. Proposition \ref{main prop of reducibility points} implies that 
$$\mathrm{Red}(\boldsymbol{\tilde{\lambda}},\lambda)= 
   \begin{cases}     \{\pm1,\pm\frac{1}{2} +\frac{\pi \sqrt{-1}}{\log q}\}
   & \text{$\tilde{\lambda}\leftrightarrow \lambda$}, \\
   \{0,\pm \frac{1}{2}+\frac{\pi \sqrt{-1}}{\log q}\} & \tilde{\lambda}\not\leftrightarrow \lambda.
   \end{cases}
$$
In the first case $\tilde{\lambda}\leftrightarrow \lambda$, we also have $\mathrm{Red}(\boldsymbol{\tilde{\lambda}}',\lambda)= 
   \{\pm\frac{1}{2} ,\pm1+\frac{\pi \sqrt{-1}}{\log q}\}$.

The calculation for ramified quadratic ${F}/{\Fo}$ is similar but a bit simpler, so we just highlight some key points. If $\tilde{\lambda}$ is self-dual, then $\tilde{\lambda}|_{\mu_{\Fo}}$ has order at most 2 and, if we assume that $\varpi^2 = \varpi_{\bullet}$, then $\boldsymbol{\tilde{\lambda}}(-\varpi_{}^2) = 1$. We compute that 
$$b_yT_y(s_y)= \tilde\lambda(-2)\lambda(-1)(q-1)
\quad\text{and}\quad 
b_z T_z(s_z) = \begin{cases}
q-1 & \tilde{\lambda}|_{\mu_{\Fo}} \equiv 1,
\\
0 & \tilde{\lambda}|_{\mu_{\Fo}}\text{ is quadratic.}
\end{cases}$$  
Since $c_y = c_z=q$, Proposition \ref{main prop of reducibility points} implies that 
$$r_y = 1,\quad r_z=\begin{cases}
1& \tilde{\lambda}|_{\mu_{\Fo}} \equiv 1,
\\
0& \tilde{\lambda}|_{\mu_{\Fo}}\text{ is quadratic,}
\end{cases} 
$$
and hence
$$
\mathrm{Red}(\boldsymbol{\tilde{\lambda}},\lambda)
 = 
 \begin{cases}
\{\pm 1,\frac{\pi \sqrt{-1}}{\log q}\} & \tilde{\lambda}|_{\mu_{\Fo}} \equiv 1\text{ and ${\tilde{\boldsymbol{\lambda}}}(\varpi_{})=\lambda(-1)$},
\\
\{0,\pm 1+\frac{\pi \sqrt{-1}}{\log q}\} & \tilde{\lambda}|_{\mu_{\Fo}} \equiv 1\text{ and ${\tilde{\boldsymbol{\lambda}}}(\varpi_{})=-\lambda(-1)$},
\\
\{\pm \frac{1}{2},\pm \frac{1}{2}+\frac{\pi \sqrt{-1}}{\log q}\}
 & \tilde{\lambda}|_{\mu_{\Fo}}\text{ is quadratic.}
\end{cases}$$
The result agrees with the calculation in \cite[(3.5) and Cor 3.6]{BT-ramified}.

\subsubsection{Example 2: ramified SO(2)}
\label{subsubsection Example 2: ramified SO(2)}

Fix a uniformizer $\varpi$ of $F$ and let $G = \SO_{2,F[\sqrt\varpi]/F}$, thr absolute rank 1 special orthogonal group that splits over $F[\sqrt\varpi]$ but not over $F$. Then $G(F)$ is just isomorphic to $\mathrm U_1(F[\sqrt\varpi]/F)$. Take a depth zero character $(G(F),\lambda)$, i.e. just a character of $F[\sqrt\varpi]^\times_1=\ker N_{F[\sqrt\varpi]/F}$. Viewing $G$ as an $F$-group, the theory of endoscopy asserts that there is an irreducible representation $\tilde\pi$ of $\GL_2(F)$ lifted from $\lambda$. It turns out that $\tilde\pi$ is non-supercuspidal, and we determine its cuspidal support, which is necessarily a set of two characters of $F^\times$.

The group $G$ is the connected isometry group defined by the matrix $H = \diag(-\varpi,1)$ and consists of elements represented by $$a+b\sqrt{\varpi} = \begin{bmatrix}
a&b \\ b\varpi & a
\end{bmatrix}, \quad a,b\in F\text{ and }a^2-b^2\varpi = 1.$$
 The element $\mathbbm p : = \diag(1,-1)\in G^\sharp(F)\smallsetminus G(F)$ conjugates $a+b\sqrt{\varpi}\mapsto a-b\sqrt{\varpi}$. We take a self-dual lattice $\Lambda$ in $V$ determined by 
\begin{equation*}
     \Lambda(0) = \mathfrak{o}_F\oplus  \mathfrak{o}_F, \quad \Lambda(1) = \mathfrak{o}_F\oplus  \mathfrak{p}_F.
\end{equation*}
As we will test the reducibility using characters, we have $\dim \tilde V_{-} = \dim \tilde V_{+}=1$, then the lattice 
$$\mathfrak m (k) = \Lambda_0(\tfrac{k-1}{3})\oplus\Lambda(\tfrac{k}{3})\oplus\Lambda_0(\tfrac{k+1}{3}) $$ 
renders the compact groups in $\mathbf G(F) = \SO_{4,F[\sqrt\varpi]/F}(F)$:
\begin{equation*}
{\mathcal J}_P = {\begin{bmatrix}
\mathfrak o^\times&\mathfrak p&\mathfrak o&\mathfrak o
\\
\mathfrak o&\mathfrak o^\times&\mathfrak o&\mathfrak o
\\
\mathfrak p&\mathfrak p&\mathfrak o^\times&\mathfrak o
\\
\mathfrak p&\mathfrak p&\mathfrak p&\mathfrak o^\times
\end{bmatrix}}\cap \mathbf G(F) 
\quad \text{and}\quad 
{\mathcal J}_{P,0_+} =I+\begin{bmatrix}
\mathfrak p&\mathfrak p&\mathfrak o&\mathfrak o
\\
\mathfrak o&\mathfrak p&\mathfrak o&\mathfrak o
\\
\mathfrak p&\mathfrak p&\mathfrak p&\mathfrak o
\\
\mathfrak p&\mathfrak p&\mathfrak p&\mathfrak p
\end{bmatrix}\cap \mathbf G(F) .
\end{equation*}
Let $\tilde G = \GL_1$ and $(\tilde G(F), \tilde{\boldsymbol \lambda})$ be a depth zero character, whose values on $\mu_F\times \left<\varpi\right>$ is yet to be determined.

First take $w=y$. For $(X,Y)^-\in \mathcal S_y$, we represent $X $ by $[0,b]$ with $b\in \mathfrak o\bmod \mathfrak p$ and $Y\in \mathfrak o^\times \bmod \mathfrak p$, then we have $2Y = -b^2$. One can compute that $I - {}^\alpha X Y^{-1}X $ is always $ \diag(1,-1) = \mathbbm p$. Hence 
$$b_y T_y(s_y) = \sum_{-b^2=2Y}\tilde\lambda (Y) \lambda( \mathbbm p^2) = \tilde\lambda(-2) (q-1).$$
Now take $w=z$. For $(X,Y\varpi^{-1})^+\in \mathcal S_z$, with representatives $X = [a,0]$ where $a\in \mathfrak o\bmod \mathfrak p$ and $Y\in \mathfrak o^\times \bmod \mathfrak p$, we have $2Y\varpi^{-1} = a^2\varpi^{-1}$ and $I - {}^\alpha X Y^{-1}X = \diag(-1,1) =- \mathbbm p$. Hence 
$$b_z T_z(s_z) = \tilde\lambda(-1)\sum_{a^2=2Y}\tilde\lambda (Y) \lambda(- \mathbbm p^2) = \tilde\lambda(-2)\lambda(-1) (q-1).$$
We therefore have two possibilities for the character $\tilde\lambda$:
$$\tilde \lambda_1 |_{\mu_F}\equiv \mathbf 1_{\mu_F}, \tilde{\boldsymbol \lambda}_1(\varpi) = \lambda(-1)
\quad\text{and}\quad
\tilde \lambda_2  |_{\mu_F}\equiv \left(\frac{\cdot}{\mu_F}\right), \tilde{\boldsymbol \lambda}_2(\varpi) = \lambda(-1).
$$
Both $\mathrm{Red}\{\tilde {\boldsymbol \lambda}_1, \lambda\}$ and $\mathrm{Red}\{\tilde {\boldsymbol \lambda}_2, \lambda\}$ contains 1. The representation $\iota^{\GL_2(F)}_{\GL_1(F)\times \GL_1(F)}(\tilde{\boldsymbol \lambda}_1\times \tilde{\boldsymbol \lambda}_2)$ is clearly irreducible, and is the lifting of $\lambda$.

\section{Main calculations}
\label{section Main calculation}

We will dive into the calculations of the signs $\epsilon_w T_w(s_w)$, with $w\in \{y,z\}$, appearing in Proposition \ref{main prop of reducibility points} and eventually deduce the relation between $\tilde\pi$ and $\pi$ such that $I(s,\tilde\pi,\pi)$ is reducible at $s=1$. The calculations will be provided in Section \ref{subsection Expanding the intertwining operator as a sum}, after the preparatory knowledge in Sections \ref{subsection Cohomological classification}
 - \ref{subsection Lattices and covers for classical groups}. The main results are given in Propositions \ref{first general form of lifting, unitary group} and \ref{first general form of lifting}. Some properties of normalized quadratic Gauss sums are necessary; they will be provided in Section \ref{subsection Non-degeneracy of a quadratic form}.

\subsection{Cohomological classification}
\label{subsection Cohomological classification}

We provide a classification of classical groups by their pure inner forms. Let $G$ be a connected classical group over $\Fo$ determined by a Hermitian form $h$, and view the first cohomology group $H^1(\Fo,{G})$ as a pointed set. Let $\hat G$ be the Langlands dual group of $G$, and $Z\hat G$ be its center. We use Kottwitz isomorphism $$H^1(\Fo,{G}) \cong \pi_0(Z\hat{{G}})^{\Gal(F^{\text{sep}}/\Fo)}$$
to enumerate the set $H^1(\Fo,G)$ explicitly.

\begin{itemize}

\item If $G$ is a unitary group, then $H^1(\Fo,G) \cong \{\pm 1\} $ corresponds bijectively to the two isomorphism classes of Hermitian spaces with the same dimension. These two classes can be distinguished by the discriminant
\begin{equation}
\label{discriminant of Hermitian forms}
\mathrm{disc}:\{\text{non-degenerate }F/\Fo\text{-Hermitian forms on }V\}\rightarrow \Fo^\times /N_{F/\Fo}(F^\times).
\end{equation}

\item For $G=\mathrm{SO}_{N}$, we view $G$ as defined by the quadratic space $(V,q)$ with $\dim V = N$ and $q(x) = h(x,x)$. Then $H^1(F,G) $ is a singleton when $N\leq 2$; otherwise, $H^1(F,G) \cong \{\pm 1\} $ corresponds bijectively to the isomorphism classes of orthogonal spaces with the same dimension and discriminant (defined similarly as in (\ref{discriminant of Hermitian forms}) but with image in $F^\times / F^\times {}^2$). The two classes of orthogonal spaces can be distinguished by the {Hasse-Witt invariant} \cite[IV, Sec 2.1]{Serre-A-course-in-arithmetic} defined as follows. If $\dim V=1$, we define $e(q)=1$; if $\dim V\geq 2$, by choosing an ordered orthogonal basis $\{v_i\}$ of $V$ which is eventually irrelevant, we define
$$e(q) = \prod_{i<j} (q(v_i), q(v_j))$$
where $(\cdot,\cdot)$ is the Hilbert symbol. Its value lies in $\{\pm 1\}$ for non-Archimedean $F$.

\item For $G=\mathrm{Sp}$, the set $H^1(F,G)  $ is trivial, so there is no non-trivial pure inner form of $G$.

\end{itemize}
From now on, we choose a quasi-split form of $G$ to be parametrized by $+1\in \{\pm 1\}$, and denote this form by $G_+$; the another form is hence $G_-$. With fixed $(\epsilon_G,F/\Fo) $ and $N = \dim_FV$, as well as a discriminant $d\in F^\times/F^{\times 2}$ when $\epsilon_G=1$ and $F=\Fo$, we temporarily define a map
$$H^1:\{\text{Hermitian forms with fixed }(\epsilon, F/\Fo, N, d)\}\rightarrow \{\pm 1\}$$
to unify the classification above. Hence $H^1$ has image $\{1\}$ in the symplectic case, $\Fo^\times / N_{F/\Fo}(F^\times)\cong \{\pm 1\}$ by the discriminant in the unitary case, and $\{\pm 1\}$ by the Hasse-Witt invariant in the orthogonal case.

For computational convenience, we will choose the following Hermitian matrices (under some choices of bases that are eventually unimportant) 
to define the forms that represent their isometry classes.

\begin{itemize}
\item If $G$ is odd ramified unitary, then the two classes of forms can be represented by
\begin{equation*}
   \begin{split}
     H_+ = \antidiag(1,-1,1,\dots,1)
     \quad\text{and}\quad H_- = \zeta H_+ 
     \end{split}
  \end{equation*}
We see that $H^1(H_+) = 1$ and $H^1( H_-)=-1$.

\item If $G$ is ramified $\SO_{2n}$ with discriminant $-\varpi'$ (with $\varpi'$ not necessarily equal to the fixed uniformizer $\varpi$), then two classes of forms can be represented by
\begin{equation}
\label{ramified even SO H plus and minus}
   \begin{split}
       H_+ & = \antidiag(1,\dots,1,\diag(-(-1)^{n-1}\varpi',1),1,\dots,1),       \\
         H_- & = \antidiag(1,\dots,1,\zeta^{}\diag(-(-1)^{n-1}\varpi',1),1,\dots,1).   \end{split}
  \end{equation}
It is easy to show that $H^1(H_+) =  - H^1( H_-) $, and both forms are quasi-split over $F$.
\end{itemize}

\subsection{Embeddings}
\label{section Embeddings of lattices}

Fix a classical group $G = G(V,h)$. Let $\mathbf s$ be an epipelagic stratum in $V$ as in Definition \ref{definition of Epipelagic stratum}. In particular, $\mathbf s$ is skew semi-simple, and $F[\beta] = \oplus_{i\in I}F[\beta_i]$ is a maximal elliptic Cartan subspace in $\tilde {\mathfrak g}(V)(F)$. Each summand $E_i:=F[\beta_i]$ is $\alpha$-invariant, and the restriction of $-\alpha$ on $E_i$ defines a Galois involution simply denoted by $c$, whose fixed field is denoted by $E_{i\bullet}$.

Put $(E_i^\times)_1 = \mathrm U_1(E_i/E_{i\bullet})$ and $E_1^\times = \prod_{i\in I}(E_i^\times)_1$. To describe the $\Fo$-embeddings of $E_1^\times$ into $G(\Fo)$, up to $G(\Fo)$-conjugacy, we first impose a Hermitian structure on $E_i$ (as an $F$-space) for each $i\in I$: for all $x_i,y_i\in E_i$,
\begin{equation}
\label{hermitian forms on E_i quadratic}
\begin{split}
h_{E_i}(x_i,y_i) &=\tr_{E_i/F}({}^cx_iy_i),
\\
 &= \tr_{E_i/F}(\beta_i{}^cx_iy_i),
\end{split}
\qquad\text{when}\qquad
\begin{split}
\epsilon_G &= 1,
\\
\epsilon_G &= -1.
\end{split}
\end{equation}
We then put $h_{F[\beta]} = \oplus_{i\in I}h_{E_i}$ on $F[\beta]$.

We now return to $(V,h)$ and consider different isometry classes of Hermitian forms. For $\delta\in \{1,\zeta\}$, let $h_i^\delta$ be the Hermitian form defined by the representative $H^\delta_i$, where $H^1_i$ is of the form $H_+$ and $H^\zeta_i$ is of the form $H_-$ in Section \ref{subsection Cohomological classification}. Given a partition on the index set 
$$I = I_1\sqcup I_\zeta$$
(and we sometimes just say $I_\zeta$ is a partition of $I$), denote by 
$$(V,h_{I_\zeta}), \quad h_{I_\zeta} = \sum_{i\in I_1}h_i^1 +\sum_{i\in I_\zeta}h_i^\zeta$$ the resulting Hermitian space defined by orthogonal sum.

In our epipelagic setting, we have $\dim _FV = [F[\beta]:F] = \sum_{i\in I}[E_i:F]$. By the classification using $H^1$ in Section \ref{subsection Cohomological classification}, 
$(F[\beta],h_{F[\beta]})$ is $F/\Fo$-isometric to any $(V,h_{I_\zeta})$ with partition $I_\zeta$ such that $H^1(h_{I_\zeta})=H^1(h_{F[\beta]})$ . This isometry map $(F[\beta],h_{F[\beta]})\rightarrow (V,h_{I_\zeta})$ then induces an embedding of $F[\beta]^\times$ into $\tilde G(V,h_{I_\zeta})(F)$, and by restriction an embedding 
$$m_{I_\zeta}: E_1^\times\hookrightarrow  G(\Fo)=G(V,h_{I_\zeta})(\Fo).$$
It is routine to show that the $G(\Fo)$-isometry class of $m_{I_\zeta}$ is independent of the choice of the isometry map, and so depends only on the partition $I = I_1\sqcup I_\zeta$. Given two partitions $I_1\sqcup I_\zeta$ and $I'_1\sqcup I'_\zeta$ of $I$, the embeddings $m_{I_\zeta}$ and $m_{I'_\zeta}$ are $G(\Fo)$-conjugate if and only if ${I_\zeta}={I'_\zeta}$, i.e., the partitions are equal. (Similar results were shown in \cite[Rem 3,26(ii) and Def 9.15]{KSS-endo-parameters}, using the language of concordance.)

The classification results of $G$ using $H^1(\Fo,G)$ in Section \ref{subsection Cohomological classification} imply that, since these invariants are just $\pm 1$, switching an index from one of the sets $I_1$ and $I_\zeta$ to another results in changing the pure inner form from one to another, and switching two indices results in a form equivalent to the original one. Therefore, given a fixed pure inner form of $G$, we can parametrize the $\Fo$-embeddings $E^\times_1\rightarrow G(\Fo)$ by partitions of $I$ such that one of $\#I_1$ and $ \#I_\zeta$ has a fixed parity. For instance, we will parametrize the embeddings into $G=G(V,h)$ with $H^1(h)=1$ by the partitions where $\#I_\zeta$ is even.

\subsection{Lattices and covers for classical groups}
\label{subsection Lattices and covers for classical groups}

We check how our previous results change our choices of lattices. 

For brevity, we identify $V=F[\beta]$ and equip on which the form $h_{I_\zeta}$ defined in Section \ref{section Embeddings of lattices}, for a choice of partition $I_\zeta$ of $I$. For each $j\in I$ and $\delta\in \{1,\zeta\}$, let $H_j = H^\delta_j\in \tilde{\mathfrak g}(V_j)(F)$ be the Hermitian matrix (i.e., ${}^{t}\overline {H^\delta_j}= \epsilon_G H^\delta_j$) corresponding to the form $h_j^\delta$ on $E_j$, then the form $h_{I_\zeta}$ can be presented as 
\begin{equation}
\label{the specific choice of hermitian form}
H =H_{I_\zeta} =  \diag(\underbrace{\dots, H^1_j,\dots}_{j\in I_1},\underbrace{\dots, H^\zeta_j,\dots}_{j\in I_\zeta}).
\end{equation}
Note that if $X\in \tilde{\mathfrak g}(V_j)(F)$, then $X\mapsto H_j^{-1}{}^{t} \overline X H_j$ defines the conjugation on $E_j$.

We now fix an index $i\in I$ and put $\tilde V_- = E_i$. Take $\tilde H = H_i$ and, with $H$ defined above, form the matrix $H_W$ in (\ref{The big matrix JW}), defining a Hermitian form $h_W$ on $W =  (\tilde V_- \oplus \tilde V_+)\perp V$. The $\alpha$-operators on the entries $(X,Y)$ in (\ref{alpha on X and Y}) become
 \begin{equation*}
     {}^\alpha X  =
-\sum_{j\in I}  H_j^{-1}{}^{t} \overline X_j H_i, \quad 
{}^\alpha Y = -\epsilon {}^{t}\overline{\tilde H}_i^{-1}  {}^{t}\overline  Y \tilde H_i =   -  H_i^{-1}{}^{t} \overline Y H_i .
\end{equation*}
Given an embedding $m_{I_\zeta}: E_1^\times \hookrightarrow G(\Fo)$, defined in Section \ref{section Embeddings of lattices} and corresponding to the partition $I_\zeta$, we define another embedding
$$\mathbf m_{I_\zeta}: E_1^\times \hookrightarrow \mathbf G(\Fo) = G(W,h_W)(\Fo), \quad g = (g_j )_{j\in I}\mapsto (m_{I_\zeta}(g_i),m_{I_\zeta}(g),m_{I_\zeta}(g_i)).
$$
When ${I_\zeta}$ is fixed, we simply denote the image of $\mathbf m_{I_\zeta}$ by $S$.

We now construct some self-dual lattices in $W$. Let $\Lambda$ and $\Lambda_i$ be the self-dual lattice sequences appearing in the strata $\mathbf s$ and its component $\mathbf s_i$ respectively. We then follow \cite[Sec 6.2]{Stevens-supercuspidal}: in the space $W_i = (\tilde V_- \oplus \tilde V_+)\perp V_i$, which is isomorphic to $E_i^{\oplus 3}$ as an $F$-space, we define two minimal self-dual $\mathfrak o_{E_i}$-lattice sequences $\mathfrak M_i^w{}$, with $w\in \{y,z\}$ such that $q^y_+=q^y_-=0$ and $q^z_+=-q^z_-=-1/2$ are the unique numbers in the interval $[0,1)$ giving
$$\mathfrak M_i^w(r)\cap \tilde V_\delta \supsetneq \mathfrak M_i^w(r_+)\cap \tilde V_\delta\quad\Leftrightarrow \quad r = q^w_\delta, \quad \text{ for }\delta\in \{+,-\}\text{ and }w\in \{y,z\}.$$
as well as $\mathfrak M_i^w(0)\cap V_i\supsetneq \mathfrak M_i^w(0)\cap  V_i$. They are contained in the maximal self-dual $\mathfrak o_{E_i}$-lattice sequence $\mathfrak m_i$ such that 
$q^\mathfrak m_+=-q^\mathfrak m_-=-1/3$ are the unique numbers in the interval $[0,1)$ giving
$$\mathfrak m_i(r)\cap \tilde V_\delta \supsetneq \mathfrak m_i(r_+)\cap \tilde V_\delta\quad\Leftrightarrow \quad r = q^\mathfrak m_\delta, \quad \text{ for }\delta\in \{+,-\},$$
and also $\mathfrak m_i(0)\cap V_i\supsetneq \mathfrak m_i(0_+)\cap  V_i$. Finally, for $\mathfrak L \in \{\mathfrak M^y,\mathfrak M^z,\mathfrak m\}$, we define
$$\mathfrak L = \mathfrak L_i \oplus \bigoplus_{j\in I\smallsetminus \{i\}} \Lambda_j.$$
These are the lattice sequences we mentioned in Section \ref{subsection Structures of Hecke algebras}, with their corresponding facets in the appropriate Bruhat-Tits buildings mentioned in Section \ref{subsection Preliminaries on the construction of covering types}.

Define for $j,k\in I\cup\{+,-\}$ and $r\in \mathbb R$,
$$\mathfrak P^r_{(j,k)}=\{X\in \Hom_F(V_k,V_j): X\Lambda_k(s)\subset \Lambda_j(s+r)\text{ for all }s\in \mathbb R\},$$
so that $\mathfrak P^r_{(j,j)} = \mathfrak P^r(\Lambda_{j})$. In the epipelagic case, since $E_j$ is a maximal subfield in $\tilde{\mathfrak g}(V_j)(F)$, we simply have $\mathfrak p^r({E_j}) = \mathfrak P^r_{(j,j)}\cap E_j$.

The compact subgroup $ {\mathcal J}_P$ is defined as 
   $$ {\mathcal J}_P :=  (\mathfrak P^0(\mathfrak m) \cap Z_{\mathbf G(\Fo)}(S))(\mathfrak P^{0_+}(\mathfrak m) \cap \mathbf G(\Fo)).$$  
To present it as a matrix group, we recall the form from \cite[1.3. Prop 1]{Blondel-cover-propag} with block size $t=3$ and obtain a presentation of the lattice \begin{equation*}
 {\mathfrak J}_P 
   = \begin{bmatrix}
    \mathfrak o_{E_i}+ \mathfrak{P}^{0_+}(\Lambda_-) & \mathfrak o_{E_i}+\mathfrak P^{0_+}_{(-,I)}
     &\mathfrak P^0_{(-,+)}
     \\
\mathfrak P^{0_+}_{(I,-)} &  \mathfrak o_{E}  +\mathfrak{P}^{0_+}(\Lambda) &  \mathfrak o_{E_i}+\mathfrak P^{0_+}_{(I,+)}
      \\
  \mathfrak  p_{E_i}+ \mathfrak P^{0_{++}}_{(+,-)} & \mathfrak P^{0_+}_{(+,I)}  &    \mathfrak o_{E_i} + \mathfrak{P}^{0_+}(\Lambda_+)
    \end{bmatrix}{} \subset \mathfrak P^0(\mathfrak m),
     \end{equation*}
where at the $(+,-)$-corner, if $\mathfrak P^{0_{+}}_{(+,-)} = \mathfrak P^{1/e_i}_{(+,-)}$, then $\mathfrak P^{0_{++}}_{(+,-)} = \mathfrak P^{2/e_i}_{(+,-)}$. We then put $\mathcal J_P = \mathfrak J_P \cap \mathbf G(\Fo)  $.

We now define the elements $s_y$ and $s_z$ appeared in Section \ref{subsection Structures of Hecke algebras}. Let $I_{(+,-)}$ be the matrix that maps the basis of $\tilde V_-$ defining the Hermitian matrix $H$ of $h$ into its dual basis in $\tilde V_+$, hence is represented by the identity matrix with the above bases. Let $I_{(-,+)}$ be defined similarly. Fix a uniformizer $\varpi_i\in E_i$ and define, for $w\in \{y,z\}$,
\begin{equation}
\label{representative-sy-and-sz}
s_y  = \text{anti-diag}(I_{(-,+)},{\mathbbm p},I_{(+,-)}),
\quad s_z =\text{anti-diag}(\varpi_{i}^{-1}I_{(-,+)},  {\mathbbm p},  - \varpi_{i}^{}I_{(+,-)}),
\end{equation}
where $\mathbbm p \in G^\sharp(\Fo) $ is just the identity unless $G$ is orthogonal, in which case $\mathbbm p$ will be a specifically chosen element in $ G^\sharp(\Fo)$ such that $\mathbbm p^2=1$ and $\det s_y = \det s_z = 1$. The elements in (\ref{representative-sy-and-sz}) satisfy $s_y^2=1$ and  $s_z^2=\diag(-I_{\tilde V},I_{V},-I_{\tilde V})$.

We put $T_w(s_w) = \tilde T_w(s_w)\times \lambda(\mathbbm p)$  for some intertwining operator $\tilde T_w(s_w)\in \End(V_{\tilde\lambda})$ and $\lambda(\mathbbm p)\in\End(V_{\lambda}) $. Again, in the epipelagic case, these operators are just scalars. We normalize $T_w$ such that $\tilde T_y(s_y)^2 = 1$ and $\tilde T_z(s_z)^2 = \tilde\lambda(-1)$.

\subsection{Expanding the intertwining operator as a sum}
\label{subsection Expanding the intertwining operator as a sum}


In this section, we compute $b_w$ for $w\in \{y,z\}$ from (\ref{formula for b and c for y}) and (\ref{formula for b and c for z}) in the case $i\in I\smallsetminus \{o\}$, while the case $i=o$ is postponed to Section \ref{section Examples for simple supercuspidals}.

Recall that the epipelagic strata $[\tilde\Lambda,0,\tilde \beta]$ and $[\Lambda,0,\beta]$ associated to the inducing types of the epipelagic representations $\pi = \cInd_{\mathcal J}^{G(F)}\lambda$ and $\tilde\pi = \cInd_{\tilde{\boldsymbol{\mathcal{J}}}}^{\tilde G(F)}\tilde{\boldsymbol{\lambda}}$ respectively. We will now assume that $\tilde\beta = 2\beta_i$ for some $i\in I\smallsetminus \{o\}$, and determine the relations between the extensions $\tilde{\boldsymbol{\lambda}}$ and ${\lambda}$ such that $I(s,\tilde\pi,\pi)$ is reducible at $s\in \mathbb R_{\geq 1}$. The results will be given in Propositions \ref{first general form of lifting, unitary group} and \ref{first general form of lifting}, which will be used to determine the endoscopic lift of $\pi$ in Section \ref{subsection Liftings of epipelagic representations for classical groups}.  

\subsubsection{Computing $b_y$}

We first compute $b_y$. By applying the calculation in \cite[Sec 4.1]{BT-ramified}, we can expend (\ref{formula for b and c for y}) into
\begin{equation}
\label{main calculation: by begins}
b_y\tilde T_y(s_y)=\sum_{(X,Y)\in \mathcal S_y}\tilde{\lambda}(Y)\lambda (I-{}^\alpha X Y^{-1}X),
\end{equation}
where $\mathcal S_y = 
      ( s_y{\mathcal J}_P^+s_y^{-1}\cap {\mathcal J}_Ps_y{\mathcal J}_P)/{\mathcal J}_P^-$. The quotient $s_y{\mathcal J}_P^+s_y^{-1}/{\mathcal J}_P^-$ consists of elements of the form 
\begin{equation*}
   \begin{split}
     (X,Y)&\in  \frac{\mathfrak o_{E_i}+\mathfrak P^{0_+}_{(-,I_y)}}{\mathfrak P^{0_+}_{(-,I_y)}}
\oplus 
\frac{\mathfrak P^{0}_{(-,+)}}{\mathfrak  p_{E_i}+ \mathfrak P^{0_{++}}_{(-,+)}}
   \end{split}
\end{equation*}
satisfying the relation (\ref{X-alpha-X-equals-Y-minus-alpha-Y}). Write $X = (X^j)_{j\in I}$, then $X^j=0$ for all $j\neq i$, and $X^i\in \mathfrak o_{E_i}/\mathfrak p_{E_i} = \mathfrak o_{F}/\mathfrak p_{F}$. The condition $(X,Y)\in {\mathcal J}_Ps_y{\mathcal J}_P$ forces $Y\in \tilde{\mathcal J} = \mathfrak o_{E_i}^\times+ \mathfrak P^{0_{+}}({\Lambda_i})$, so that we write $Y = Y_0(1+Y_1) \in \mathfrak o_{E_i}^\times/ \mathfrak p_{E_i}+ \mathfrak P^{0_+}({\Lambda_i})/(\mathfrak p_{E_i}+\mathfrak P^{0_{++}}({\Lambda_i}))$. 
Hence (\ref{X-alpha-X-equals-Y-minus-alpha-Y}) implies that 
\begin{subequations}
\begin{align}
     & Y_0+  {}^c Y_0 = - X^i {}^c X^i,
     \label{relation of X and Y at level 0 and w=y}
     \\
     & Y_0Y_1- {}^\alpha Y_1{}^cY_0=0\quad \text{i.e., }{}^\alpha Y_1 =  Y_0Y_1{}^cY_0^{-1}.
     \label{relation of X and Y at level 1 and w=y}
\end{align}
\end{subequations}
We first obtain
\begin{equation}
\label{main calculation: tilde-lambda-Y in by}
 \tilde{\lambda}(Y) = \tilde\lambda(Y_0)\psi_{\tilde\beta}(Y_1).
\end{equation}
We then compute $\lambda(I-{}^\alpha X Y^{-1}X)$. We represent $\beta = \oplus_{i\in I}\beta_i$ by diagonal blocks, and look at the corresponding blocks of $I-{}^\alpha X Y^{-1}X$, which are $I_j$ for all $j\neq i$, and is 
\begin{equation}
\label{value W_0(1+W_1)}
   I_i-{}^\alpha X^i Y^{-1}X^i =  I_i - {}^\alpha X^i(I+Y_1)^{-1}Y_0^{-1}X^i \equiv I_i - {}^\alpha X^i(I-Y_1)Y_0^{-1}X^i \mod \mathcal{U}^{0_{++}}(\Lambda_i)
\end{equation}
when $j=i$ the $i$-th diagonal block. From here we branch into two cases, depending on the type of $G$.

If $G$ is orthogonal, symplectic, or ramified unitary, we follow similar arguments in \cite[Lem 4.2]{BT-ramified}. Indeed, (\ref{relation of X and Y at level 0 and w=y}) becomes $2Y_0 = -(X^i)^2$, and so $X^i\in \mathfrak o_{E_i}^\times$. Hence indeed 
$I_i-{}^\alpha X^i Y_0^{-1}X^i =-I_i$, and we obtain
\begin{equation}
\label{main calculation: lambda-p-W in by, not unramified unitary}
\lambda(I-{}^\alpha X Y^{-1}X) =\lambda(\omega_i)\psi_{\beta_i}(2Y_1)^{-1},
\end{equation}
where $\omega_i=\diag((I_j)_{j\neq i},-I_i)$ as in (\ref{definition of omega_i element}). Since $\tilde\beta = 2\beta_i$, putting (\ref{main calculation: tilde-lambda-Y in by}) and (\ref{main calculation: lambda-p-W in by, not unramified unitary}) into (\ref{main calculation: by begins}) yields
\begin{equation*}
 c_y = q^2, \quad b_y T_y(s_y) = \tilde\lambda(-2)\lambda(\omega_i)(q-1) = \tilde\lambda(-2)\lambda(\omega_i) (q^{1/2}-q^{-1/2})(c_y/q)^{1/2},
\end{equation*}
which implies that
\begin{equation}
\label{results of ry and eyTy}
 r_y =1\quad\text{and} \quad \epsilon_y T_y(s_y)= \tilde\lambda(-2)\lambda(\omega_i ).
\end{equation}
This result is regardless of whether $\tilde\lambda$ is trivial or quadratic.

If $G$ is unramified unitary, then $I$ is just a singleton. We express the value in (\ref{value W_0(1+W_1)}) as $W_0(I+W_1)$ where $W_0\in \mu_F$ and $I+W_1\in \mathcal{U}^{0_+}(\Lambda)$, then we can take 
\begin{equation*}
W_0 = I - {}^\alpha XY_0^{-1}X   = 
   \begin{cases}
     1& \text{if $X=0$,}\\
        -Y_0^{-1}{}^cY_0& \text{if $X\neq 0$},
   \end{cases}
   \end{equation*}
and $ W_1 =  W_0^{-1}{}^\alpha XY_1Y_0^{-1}X$. We continue separating into cases $X\neq 0$ and $X=0$. In the former case, we have $
 W_1 =  -(1+ {}^cY_0^{-1}Y_0)Y_1$. Using (\ref{relation of X and Y at level 1 and w=y}) and also ${}^\alpha \beta = \beta$, we have $\psi_{\beta}(W_1) =\psi_{\beta}(-2Y_1)$ and so
$$\lambda(I-{}^\alpha X Y^{-1}X) =\lambda(W_0)\psi_{\beta}(W_1)=\lambda(-{}^cY_0 Y_0^{-1})\psi_{\beta}(2Y_1)^{-1}.$$
Since $\tilde\beta = 2\beta$, the summand with $X\neq 0$ is 
\begin{equation}
\label{by unramified unitary X not 0 summand}
\sum_{
   \begin{smallmatrix} X\neq 0 \\ Y_1   \end{smallmatrix}} \tilde{\lambda}(Y) \lambda(I-{}^\alpha X Y^{-1}X) = \lambda(-1)(\#Y_1)^{1/2}\sum_{X\neq 0}\tilde\lambda(Y_0)\lambda({}^cY_0 Y_0^{-1})
\end{equation}
(the factor being $(\#Y_1)^{1/2}$ instead of $(\#Y_1)$ because of (\ref{relation of X and Y at level 1 and w=y})). In the latter case, we have $W_0=1$ and $W_1=0$. The relation (\ref{relation of X and Y at level 0 and w=y}) becomes ${}^cY_0 = -Y_0$, and so $\psi_{\tilde\beta}(Y_1) = \psi_{\tilde\beta}(^\alpha Y_1) =  \psi_{\tilde\beta}(- Y_0Y_1Y_0^{-1})  = \psi_{\tilde\beta}(- Y_1) $ which implies that $\psi_{\tilde\beta}(Y_1) =1$ if $p\neq 2$. Hence the summand with $X=0$ is 
\begin{equation}
\label{by unramified unitary X equal 0 summand}
\sum_{ \begin{smallmatrix} X= 0 \\ Y_1   \end{smallmatrix}} \tilde{\lambda}(Y) \lambda(I-{}^\alpha X Y^{-1}X) = \lambda(-1 )(\#Y_1)^{1/2}\sum_{X= 0}\tilde\lambda(Y_0).
\end{equation}
The total sum (\ref{by unramified unitary X not 0 summand})+(\ref{by unramified unitary X equal 0 summand}) is $\lambda(-1)(\#Y_1)^{1/2}$ times the sum considered in the first example of Section \ref{Example: depth zero characters}, and the results there implies that, when $G$ is unramified unitary: 
\begin{equation*}
c_y = \qo^3(\#Y_1), \quad 
b_y \tilde T_y(s_y)= \begin{cases}
\lambda(-1) (\qo^3-1)(c_y/\qo^3)^{1/2} & \tilde\lambda \leftrightarrow \lambda,
\\
-\lambda(-1) \qo(\qo-1)(c_y/\qo^3)^{1/2}
 & \tilde\lambda \not\leftrightarrow \lambda\text{ but is self-dual},
\end{cases}
\end{equation*}
which implies that
\begin{equation}
\label{results of ry and eyTy for unramified unitary group}
r_y = \begin{cases}
3/2
\\
1/2
\end{cases}
\quad\text{and}\quad 
\epsilon_y T_y(s_y)= \begin{cases}
\lambda(-1)
\\
-\lambda(-1).
\end{cases}
\end{equation}
Note that the results in (\ref{results of ry and eyTy}) and (\ref{results of ry and eyTy for unramified unitary group}) are independent of whether $i\in I_1$ or $I_\zeta$.

\subsubsection{Computing $b_z$}

We then compute $b_z$. By applying the calculation in \cite[Sec 4.2]{BT-ramified}, we expand (\ref{formula for b and c for z}) into
\begin{equation}
\label{main calculation: bz begins}
b_z \tilde T_z(s_z)=\sum_{(X,Y)\in \mathcal S_z}\tilde{\lambda}(Y\varpi_i) \lambda(I-{}^\alpha X Y^{-1}X),
\end{equation}
where $
  \mathcal S_z = 
    (s_z{\mathcal J}_P^-s_z^{-1}\cap {\mathcal J}_Ps_z{\mathcal J}_P)/{\mathcal J}_P^+ $. The quotient $s_z{\mathcal J}_P^-s_z^{-1}/{\mathcal J}_P^+$ consists of elements 
\begin{equation*}
   \begin{split}
      (X,Y)& \in     \frac{\mathfrak P^0_{(+,I_y)}}{\mathfrak o_{E_i}+\mathfrak P^{0_{+}}_{(+,I_y)}}
\oplus
\frac{\mathfrak  p_{E_i}^{-1}+\mathfrak P^{0_{}}_{(+,-)}}{ \mathfrak P^{0_{}}_{(+,-)}},
   \end{split}
\end{equation*}
satisfying the relation (\ref{X-alpha-X-equals-Y-minus-alpha-Y}). The entries $X^j$ for $j\in I\smallsetminus\{i\}$ lies in $\mathfrak W^{}_{z,j}:={\mathfrak P^{0_{}}_{(+,j)}}/{\mathfrak P^{0_{+}}_{(+,j)}}$ while $X^i\in \mathfrak W^{}_{z,i}:={\mathfrak P^{0_{}}(\Lambda_{i})}/({\mathfrak o_{E_i}+ \mathfrak P^{0_{+}}(\Lambda_{i}}))$. The condition $(X,Y)\in {\mathcal J}_Ps_z{\mathcal J}_P$ forces $Y\in \varpi_i^{-1}\tilde {\mathcal J}$, so that we write $Y = \varpi_i^{-1}Y_0(I+Y_1)$ where $Y_1$ is just auxiliary, and choose $Y_0\in \mu_{E_i} = \mu_F$ such that 
$$
Y-{}^\alpha Y = \varpi_i^{-1}
Y_0(I+Y_1)- (I-{}^\alpha Y_1 ){}^cY_0\varpi_i^{-1} = 
X{}^\alpha X.$$
Comparing valuations, the above relation implies that
\begin{subequations}
\begin{align}
  & Y_0 = {}^cY_0\quad \text{ (i.e., $Y_0\in \mu_{\Fo}$, since ${}^\alpha\varpi_i = \varpi_i$)},
  \label{relation of X and Y at level 0 and w=z}
   \\
&\varpi_i^{-1}
Y_0Y_1+ {}^\alpha Y_1 {}^cY_0\varpi_i^{-1} = 
X{}^\alpha X . 
     \label{relation of X and Y at level 1 and w=z}
\end{align}
\end{subequations}
We have $\tilde{\lambda}(Y\varpi_i)=\tilde \lambda(Y_0)\psi_{\tilde \beta}(1+Y_1) $, and using (\ref{relation of X and Y at level 1 and w=z}), the last factor is
\begin{equation}
\label{bz X sum 1}
\psi_{\tilde \beta}(1+Y_1) =  \psi(\beta_i(Y_1+{}^\alpha Y_1)) = \psi( \beta_i\varpi_i Y_0^{-1}
\sum_{j\in I}X^j{}^\alpha X^j ).
\end{equation}
For  $\lambda(I-{}^\alpha X Y^{-1}X)$, it is easy to see that $I-{}^\alpha X Y^{-1}X = I-{}^\alpha X Y_0^{-1}\varpi_i X \mod \mathcal{U}^{0_{++}}(\Lambda)$. The calculation again reduces to looking at diagonal blocks, so that  
\begin{equation}
\label{bz X sum 2}
\lambda(I-{}^\alpha X Y^{-1}X)= \prod_{j\in I_{}}\psi(- \beta_j{}^\alpha X^j\varpi_i Y_0^{-1} X^j).
\end{equation}
Putting (\ref{bz X sum 1}) and (\ref{bz X sum 2}) into (\ref{main calculation: bz begins}), we obtain
\begin{equation}
\label{bz expanded as quadratic form}
   \begin{split}
b_z \tilde T_z(s_z) = \sum_{Y_0}\tilde \lambda(Y_0)\sum_X 
\prod_{j\in I_{}}\psi( \varpi_iY_0^{-1}( \beta_i X^j - X^j \beta_j) {}^\alpha X^j  ).
   \end{split}
\end{equation}
We will show in Section \ref{subsection Non-degeneracy of a quadratic form} that the last sum $\sum_X$ is a quadratic Gauss sum on the space 
\begin{equation}
\label{direct sum of mathfrak W}
\mathfrak W^{i}_z :=\mathfrak W^{}_{z,i}\oplus
\bigoplus_{j\in I\smallsetminus\{i\}} 
\mathfrak W^{}_{z,j}  = {\mathfrak P^{0_{}}_{i}}/({\mathfrak o_{E_i}+ \mathfrak P^{0_{+}}(\Lambda_{i}})) \oplus
\bigoplus_{j\in I\smallsetminus\{i\}}
{\mathfrak P^{0_{}}_{(+,j)}}/{\mathfrak P^{0_{+}}_{(+,j)}}
\end{equation}
containing $X$, equipped with a non-degenerate quadratic form. We will also compute this Gauss sum, which is then of the form
\begin{equation}
\label{Gauss sum of the form}
\left(\frac{Y_0}{\mu_\Fo}\right)^{\dim_{\mathbb F_\bullet}\mathfrak{W}^{i}_z}
q^{\dim_{\mathbb F_\bullet} \mathfrak{W}^{i}_z /2}\mathfrak{n}_z(\varpi_i,\mathbf{s},\psi,h) , 
\end{equation}
 where 
$\mathfrak{n}_z(\varpi_i,\mathbf{s},\psi,h)$ is a 4-th root of unity.

We put the form (\ref{Gauss sum of the form}) of the Gauss sum into (\ref{bz expanded as quadratic form}): 
\begin{equation}
\label{b_z T_z(s_z) general form}
b_z \tilde T_z(s_z)=\mathfrak{n}_z(\varpi_i,\mathbf{s},\psi,h)
q^{\dim_{\mathbb F_\bullet} \mathfrak W^i_z /2}
\sum_{Y_0\in \mu_{\Fo}}\tilde \lambda(Y_0)\left(\frac{Y_0}{\mu_\Fo}\right)^{\dim_{\mathbb F_\bullet}\mathfrak{W}^i_z},
\end{equation}
and obtain
\begin{equation*}
\text{$b_z=0$ \quad if and only if \quad $\tilde \lambda|_{\mu_{\Fo}} \neq \left(\frac{\cdot}{\mu_{F_\bullet}}\right)^{\dim_{\mathbb F_\bullet}\mathfrak{W}^i_z}$.}
\end{equation*}
We hence assume otherwise, $b_z\neq 0$, in the subsequent discussions. In this case, we have
\begin{equation*}
b_z \tilde T_z(s_z)= \mathfrak{n}_z(\varpi_i,\mathbf{s},\psi,h)
q^{\dim_{\mathbb F_\bullet} \mathfrak W^i_z /2}
(\qo-1).
\end{equation*}

\subsubsection{Preliminary results on endoscopic liftings}

When $G$ is unramified unitary, the index set $
I$ is a singleton. As $\mathfrak{W}_z = \mathfrak{W}^i_z$ is an $\mathbb F$-space, $\dim_{\mathbb F_\bullet}\mathfrak{W}_z$ is always even, and so the sum $\sum_{Y_0}$ in  (\ref{b_z T_z(s_z) general form}) is just $\qo-1$. We 
hence have 
\begin{equation}
\label{}
c_z = q^{(\dim_{\mathbb F_\bullet} \mathfrak W_z -1)/2}, \quad r_z = 1/2 
\quad\text{and}\quad 
\epsilon_z T_z(s_z)= \mathfrak{n}_z(\varpi_i,\mathbf{s},\psi,h_{})
\end{equation}
which gives the following proposition when combined with (\ref{results of ry and eyTy for unramified unitary group}).

\begin{prop}
\label{first general form of lifting, unitary group}
Let $G=G(V,h)$ be an unramified unitary group, and $\pi = \cInd_{\mathcal J}^{G(F)}\lambda$ be an epipelagic representation constructed from an epipelagic simple stratum $[\Lambda,0,\beta]$, i.e., $\lambda|_{\mathcal{U}^{0_+}(\Lambda)^\sigma} = \psi_\beta$. We construct a character $( F[\beta]^\times \mathcal{U}^{0_+}(\Lambda),\tilde{\boldsymbol{\lambda}})$ as follows
\begin{equation*}
 \lambda|_{\mathcal{U}^{0_+}(\Lambda)} = \psi_{2\beta},\quad 
 \tilde \lambda|_{\mu_F}\leftrightarrow \lambda|_{\mu_1(F)},
\quad\text{and}\quad
\tilde{\boldsymbol{\lambda}}( \varpi_i)= \lambda(-1)
\mathfrak{n}_z(\varpi_\beta,\mathbf{s},\psi,h).
\end{equation*}
Let $\tilde \pi = \tilde \pi_{\tilde{\boldsymbol{\lambda}}}$ be the associated epipelagic representation of $\tilde G(F)$, then we have $$\mathrm{Red}(\tilde\pi,\pi) = \{\pm 1, \pm \tfrac{1}{2}+\tfrac{\pi \sqrt{-1}}{\log q} \}.$$
 In particular, $I(s,\tilde\pi,\pi)$ is reducible at $s=1$. 
\end{prop}

When $G$ is not unramified unitary, we interpret the results in (\ref{b_z T_z(s_z) general form}) directly and obtain 
\begin{equation}
\label{}
c_z = q^{(\dim_{\mathbb F} \mathfrak W^i_z -1)/2}, \quad r_z = 1
\quad\text{and}\quad 
\epsilon_z T_z(s_z)= \mathfrak{n}_z(\varpi_i,\mathbf{s},\psi,h).
\end{equation}
Combining it with (\ref{results of ry and eyTy}), we get the following proposition.

\begin{prop}
\label{first general form of lifting}
Let $G=G(V,h)$ be a connected classical group not of unramified unitary type, and $\pi = \cInd_{\mathcal J}^{G(F)}\lambda$ be an epipelagic representation constructed from a epipelagic semi-simple stratum $[\Lambda,0,\beta]$, i.e., $\lambda|_{\mathcal{U}^{0_+}(\Lambda)^\sigma} = \psi_\beta$. Write $\beta = \sum_{i\in I}\beta_i$. Fix $i\in I\smallsetminus \{o\}$ and construct a character $( F[\beta_i]^\times \mathcal{U}^{0_+}(\Lambda_i), \tilde{\boldsymbol{\lambda}}_i)$ as follows
\begin{equation*}
\tilde\lambda|_{\mathcal{U}^{0_+}(\Lambda_i)} = \psi_{2\beta_i},\quad 
\tilde \lambda_i|_{\mu_{F}}= \left(\frac{\cdot}{\mu_{F}}\right)^{\dim_{\mathbb F}\mathfrak{W}^i_z},
\quad\text{and}\quad
\tilde{\boldsymbol{\lambda}}_i(\varpi_i)= \tilde\lambda(-2)\lambda(\omega_i )
\mathfrak{n}_z(\varpi_i,\mathbf{s},\psi,h).
\end{equation*}
Let $\tilde \pi_i = \tilde \pi_{\tilde{\boldsymbol{\lambda}}_i}$ be the associated epipelagic representation of $\tilde G(V_i)(F)$, then we have $$\mathrm{Red}(\tilde\pi_i,\pi) = \left\{\pm 1,\pm \tfrac{1}{2}+\tfrac{\pi \sqrt{-1}}{\log q} \right\}.$$ In particular, $I(s,\tilde\pi_i,\pi)$ is reducible at $s=1$. 
\end{prop}

In Sections \ref{section Examples for simple supercuspidals} and \ref{subsection Reducibility results for different classical groups} below, we will further analyze the value $\tilde{\boldsymbol{\lambda}}(\varpi_i)$ for different types of $G$. Note that the choice of the uniformizer $\varpi_i$ of $E_i$ is arbitrary, and our result is indeed independent of this choice.

\begin{prop}
\label{prop of independence: uniformizer}
The results in  Propositions \ref{first general form of lifting, unitary group} and \ref{first general form of lifting} relating $\tilde{\boldsymbol{\lambda}}(\varpi_i)$ and $\mathfrak{n}_z(\varpi_i,\mathbf{s},\psi,h_{})$ is independent of the choice of $\varpi_i$. 
\end{prop}
\proof
The proof is exactly the same as in \cite[Prop 4.8(ii)]{BT-ramified}. It is essentially because, to obtain the value in (\ref{b_z T_z(s_z) general form}), we assigned $\tilde\lambda|_{\mu_{\Fo}}$ to be the character defined on $Y_0\in \mu_{\Fo}$ brought out from the Gauss sum in (\ref{bz expanded as quadratic form}).  
\qed

\subsection{Properties of a quadratic form}
\label{subsection Non-degeneracy of a quadratic form}

We first recall a summary of general properties of quadratic Gauss sums from \cite[Sec 4.5]{BH-ET2}. Let $(\mathbf V,\mathbf q)$ be a quadratic form on an $\mathbb F$-vector space. Assume that $\mathbf q$ is non-degenerate. With the fixed non-trivial additive character $\psi$ of $\mathbb F$, we define the normalized quadratic Gauss sum
$$\mathfrak n(\mathbf q) = q^{-\dim_{\mathbb F}\mathbf V/2}\sum_{X\in \mathbf V}\psi(\mathbf q(X)).$$
Its value can be expressed as follows. Define the symmetric bilinear form $\mathbf h$ associated to $\mathbf q$ as 
\begin{equation}
\label{definition of bilinear form}
\mathbf h(X,Y) = \tfrac{1}{2}(\mathbf q(X+Y) - \mathbf q(X)-\mathbf q(Y)),\quad X,Y\in \mathbf V.
\end{equation}
If $\mathbf H$ is the symmetric matrix associated to $\mathbf h$, then $\det \mathbf H\neq 0$, and we simply denote $\det \mathbf q := \det \mathbf H $. Put $\mathfrak n_\psi = \mathfrak n(\mathbf q_0)$ where $\mathbf q_0:x\mapsto x^2$ on $\mathbf V = \mathbb F$, then we have
$$\mathfrak n(\mathbf q) = \left(\frac{\det \mathbf q}{\mathbb  F^\times}\right)\mathfrak n_\psi^{\dim_{\mathbb F}\mathbf V}.$$
It is well-known that $\mathfrak n_\psi$ is a 4th root of unity, and hence so is $\mathfrak n(\mathbf q)$.

We now analyze the normalized Gauss sum $\mathfrak{n}_z(\varpi_i,\mathbf{s},\psi,h_{})$ appeared in the previous section. With the setup from (\ref{bz expanded as quadratic form}), in particular a fixed index $i\in I\smallsetminus\{o\}$, we have a bilinear form on $\mathfrak{W}_{z,j}$, for $j\in I$,  by
\begin{equation*}
\mathbf q_{z,\mathbf s,j}(X^j) =\tr_{\tilde{ \mathfrak g}(V_i)/\mathbb F}(\varpi_i
(\beta_iX^{j}{}-  X^{j}\beta_j ) {}^\alpha X^j) \mod \mathfrak{p}_F,\quad \text{$X^j\in \mathfrak{W}_{z,j}$},
\end{equation*}
(note that $\mathfrak{W}_{z,i}$ is structurally different from the other $\mathfrak{W}_{z,j}$, see (\ref{direct sum of mathfrak W})), and the orthogonal sum
$
\mathbf q_{z,\mathbf s} = \perp_{j\in I}\mathbf q_{z,\mathbf s,j},
$
a bilinear form equipped on $\mathfrak{W}^i_{z}=\bigoplus_{j\in I}\mathfrak{W}_{z,j}$. Using the independence in Proposition \ref{prop of independence: uniformizer}, we pick $\varpi_i = \beta_i^{-1}$ in the case $i\neq o$, and so
$$\mathbf q_{z,\mathbf s}:X\mapsto \sum_{j\in I}\tr_{\tilde{ \mathfrak g}(V_i)/\mathbb F}(
( X_{j}{}- \beta_i^{-1} X_{j}\beta_j ) {}^\alpha X_j).$$

\begin{prop}
\label{properties of quadratic forms}
\begin{enumerate}[(i)]
\item The quadratic form $\mathbf q_{z,\mathbf s}$ is non-degenerate. \label{non-degenerate quad form}

\item 
The form $\mathbf q_{z,\mathbf s}$ is Hermitian-symmetric, i.e., ${}^\alpha\mathbf q_{z,\mathbf s} = -\mathbf q_{z,\mathbf s}$.
\label{Hermitian symmetry of the quadratic form}

\end{enumerate}

\end{prop}
\proof
To show \ref{non-degenerate quad form}, it suffices to show that each $\mathbf q_{z,\mathbf s,j}$ is non-degenerate. If $\beta_i$ is non-null, then it is equivalent to show that 
$$\beta_i X_{j}{}- X_{j}\beta_j \in \mathfrak{P}^{-s}_{(i,j)}
\quad\Rightarrow \quad
X\in 
\begin{cases}
\mathfrak{p}^{s}({E_i})+
\mathfrak{P}^{s_+}_{(i,j)}  & j=i,
\\
\mathfrak{P}^{s_+}_{(i,j)} & j\neq i.
\end{cases} 
$$
In the epipelagic case, $s=0$. The first case ($j=i$) is proved in \cite[Prop 4.7]{BT-ramified}, and indeed the second case ($j\neq i$) is just analogous if we apply \cite[Lem 3.7(i)]{Stevens-semi-simple-char}. If $\beta_i=0$, then we just need to show that 
$$X_{j}\beta_j \in \mathfrak P^0_{(i,j)}
\quad\Rightarrow \quad
X_j\in 
\mathfrak{P}^{0_+}_{(i,j)},$$
which is obvious since $v_{\Lambda_j}(\beta_j)<0$. \ref{Hermitian symmetry of the quadratic form} comes from a simple calculation that ${}^\alpha(
( X_{j}{}- \beta_i^{-1} X_{j}\beta_j ) {}^\alpha X_j) = -X({}^\alpha X -\beta_j{}^\alpha X \beta_i^{-1} )$, using the $\alpha$-invariance of $\beta_j$ for all $j\in I\smallsetminus\{o\}$.
\qed

Expanding by definition (\ref{definition of bilinear form}) and using the symmetry in Proposition \ref{properties of quadratic forms}\ref{Hermitian symmetry of the quadratic form}, the associated symmetric bilinear form is hence
\begin{equation*}
\mathbf h_{z,\mathbf s,j}(X^j,Y^j) =\tr_{\tilde{ \mathfrak g}(V_i)/\mathbb F}((X^{j}{}- \beta_i^{-1} X^{j}\beta_j ) {}^\alpha Y^j), \quad \text{$X^j,Y^j\in \mathfrak{W}_{z,j}$}.
\end{equation*}
When $j\neq i$, the discriminants of $\mathbf h_{z,\mathbf s,j}$ is equal to $\gamma_j\theta_j$, where 
\begin{equation}
\label{discriminants for j neq i}
\gamma_j = 1-\det \beta_j \det\beta_i^{-1},\quad \theta_j = \det H_j \det H_i^{-1}, \quad j\neq i.
\end{equation}
When $j=i$, we will compute in the next section the discriminants for different types of groups, as we have reduced to the case of simple supercuspidal representations in the sense of \cite{Gross-Reeder}. Note that $\theta_j$ depends on the Hermitian form defining $G=G(V,h)$; the specific choice was made in (\ref{the specific choice of hermitian form}).

\section{Examples: simple supercuspidals}
\label{section Examples for simple supercuspidals}

In a series of papers \cite{Oi-SO-odd,Oi-Sp-and-SO-even,Oi-U-unram}, M. Oi computed the endoscopic liftings of simple supercuspidals using the endoscopic character identity. In this section, we compare our preliminary results in Proposition \ref{first general form of lifting} with his and show that both provide the same liftings. Some of these results are also required to further compute the endoscopic liftings of epipelagic representations. 

\subsection{Constructions}

We first recall the explicit construction of simple supercuspidals. Our treatment here is a slightly more general from Oi: there he chose in advance some convenient representatives of conjugacy classes of affine generic characters. We will consider these characters in general and show that our results are the same as those given by Oi. This generality has the advantage for switching across equivalent simple strata.

In subsequent subsections, we will provide the descriptions of simple supercuspidals for all types of quasi-split classical groups except the ramified unitary groups, and apply (\ref{criterion of isomorphic simple supercuspidal}) to provide conditions for two affine generic characters inducing isomorphic supercuspidals. We also modify the setting to accommodate the liftings for non-quasi-split pure inner forms, which are necessary for the calculations in Section \ref{subsection Reducibility results for different classical groups}.

\subsubsection{General linear groups} 
Let $\tilde G = \GL_{m}$, then the simple affine roots $\Delta_{\mathrm{aff}} = \{\alpha_i\}_{i=1}^{m}$ are $\alpha_i = e_i- e_{i+1}$ for  $1\leq i\leq m-1$ and $\alpha_0 = 1-(e_1-e_m)$. We express elements in the pro-p Iwahori subgroup as 
$$u =I+\antidiag( 
\diag( u_1, \dots, u_{m-1}), 
u_0\varpi )\in \mathcal{I}^+,$$ 
 and the affine generic character $\tilde{\lambda}$ restricts to each root space as $u_i\mapsto  \psi_{a_i}(u_i)$ for all $i$, for some $a_i\in \mathfrak{o}_F^\times$. Put $\tilde a = \prod_{i=0}^{m-1}a_i\bmod \mathcal{U}^1(F)$.

It is easy to show that  $N(\tilde{\lambda})= \Omega Z \mathcal{I}^+$. The generator of the cyclic group $\bar \Omega:=\Omega/(\Omega\cap Z \mathcal{I}^+)$ of order $m$ can be lifted to 
$$\tilde\varpi = \antidiag( 
\diag( a_1^{-1}, \dots, a_{m-1}^{-1}), 
a_0^{-1}\varpi),$$
so that $\tilde\varpi^m = (\varpi \tilde a^{-1})I\in Z$. Take a character $\phi$ of $\mu_F\subset Z$ and any $\xi \in \mathbb C^\times$ to define the extended character $\tilde{\boldsymbol{\lambda}} = \tilde{\boldsymbol{\lambda}}((a_{i})_{i=0}^{m-1},\phi,\xi)$ by 
$$\tilde{\boldsymbol{\lambda}}(\tilde\varpi^j z u) = \xi^j\phi(z)\prod_{i}\psi_{a_i} (u_i),\,\qquad \tilde\varpi^j zu \in N(\tilde\lambda) =  \left<\tilde\varpi\right>\mu_F \mathcal{I}^+.$$
The two simple supercuspidal representations $\tilde \pi_{\tilde{\boldsymbol{\lambda}}_j}$ for $j=1,2$, where $\tilde{\boldsymbol{\lambda}}_j =\tilde{\boldsymbol{\lambda}}((a_{i}^j)_{i=0}^{m-1},\phi_j,\xi_j)$, are isomorphic if and only if $(\tilde a^1,\phi_1 ,\xi_1 )= (\tilde a^2,\phi_2, \xi_2)$.

\begin{rmk}
\cite[2.2 Prop]{BK-epipelagic} asserts that $\tilde \pi_{\tilde{\boldsymbol{\lambda}}_1}\cong \tilde \pi_{\tilde{\boldsymbol{\lambda}}_2}$ if and only if $(\tilde a^1,\phi_1  )= (\tilde a^2,\phi_2)$ and $$\epsilon(\tilde \pi_{\tilde{\boldsymbol{\lambda}}_1},1/2,\psi) = \epsilon(\tilde \pi_{\tilde{\boldsymbol{\lambda}}_2},1/2,\psi),$$
 where $\epsilon(\tilde \pi,s,\psi)$, with $s\in \mathbb C$, is the Godement-Jacquet local constant of $\tilde\pi$. By \emph{loc. cit.} [2.2 Lem(1)] , the last condition is equivalent to $\xi_1=\xi_2$.
\qed\end{rmk}

An Oi's representative of $\tilde{\boldsymbol{\lambda}}$ takes the form $(\tilde a,\tilde\lambda|_{\mu_F},\tilde{\boldsymbol{\lambda}}(\tilde\varpi))\in \mathbb F^\times \times \hat{\mathbb F}^\times \times \mathbb C^\times$, where $\tilde a$ represents $(a_0,\dots,a_{m-1})  $ $= (\tilde a,1,\dots,1)$.

\subsubsection{Unramified unitary groups}  

Let $G = \mathrm{U}_{N,F/\Fo}$ where $F/\Fo$ is unramified, and take $\varpi\in \Fo$. We take $H = \diag(1,-1,1,-1,\dots)$, so that ${}^t\overline H = (-1)^{n-1} H$. Put $n=\lfloor N/2\rfloor$. The simple affine roots $\Delta_{\mathrm{aff}} = \{\alpha_i\}_{i=1}^{n}$ are $\alpha_i = e_{i+1}-e_i$ for $i<n$, 
$\alpha_n = e_n$ if $N$ is odd or $2e_n$ if $N$ is even, and $\alpha_n = 1-2\sum_{i=1}^n{\alpha_i}$. We express
\begin{equation*}
   \begin{split}
     u &=I+ 
     \antidiag( 
\diag( u_1, \dots, u_n, \overline{u_n}, \dots, \overline{u_1}), 
u_0\varpi ),\quad \text{$\overline{u_0} = -u_0$, when $N$ is odd, or }
\\
    & = I+ \antidiag( 
\diag( u_1, \dots,  u_{n-1}, u_n, \overline{u_{n-1}}, \dots, \overline{u_1}), 
u_0\varpi )\quad \text{$u_0,u_n\in \mathfrak o_{\Fo}\bmod \mathfrak p_{\Fo}$, when $N$ is even.}   \end{split}
\end{equation*}
The affine generic character ${\lambda}$ restricts to each root space as $u_i\mapsto  \psi_{a_i}(u_i)$ for all $i$, for some $a_i\in \mathfrak{o}_F^\times$, with extra conditions: $a_0\in \ker\tr_{F/\Fo}$ when $N$ is odd, and $a_0,a_n\in {\Fo}$ when $N$ is even.

 The group $\Omega$ is trivial, and $N_G(\lambda) = Z \mathcal{I}^+$.  Take a character $\phi$ of $\mu(F)_1\subset Z$ and define the extended character $\lambda = \lambda((a_{i})_{i=0}^{n},\phi)$ by 
$$\lambda( z u) = \phi(z)\prod_{i}\psi_{a_i} (u_i),\,\qquad zu \in N_G( \lambda) =  \mu(F)_1 \mathcal{I}^+.$$
 Put
  \begin{equation*}
   \begin{split}
  a&=a_{0}\left(\textstyle\prod_{i =1}^{n} N_{F/\Fo}a_{i}\right)
  \mod (F^\times)^2\mathcal{U}^1(F)
    \quad \text{when $N$ is odd, }
\\
&=a_{0}\left(\textstyle\prod_{i =1}^{n-1} N_{F/\Fo}a_{i}\right)a_{n} \mod \mathcal{U}^1(F)
    \quad \text{when $N$ is even.}   \end{split}
\end{equation*}
 Two representations $\pi_{\lambda_j}$ for $j=1,2$, where $\lambda_j = \lambda((a_{i}^j)_{i=0}^{n},\phi_j)$, are isomorphic if and only $(a^1 ,\phi_1)= (a^2, \phi_2)$.

From \cite{Oi-U-unram}, an Oi's representative of $\lambda$ takes the form $(a,\lambda|_{\mu(F)_1})$ where $a$ represents $(a_0,\dots,a_{n}) = (a,1,\dots,1)$ with $a\in \mathbb F_\bullet^\times$ (resp. $\ker\tr_{\mathbb F/\mathbb F_\bullet}\cap \mathbb F^\times$) if $N=2n$ is even (resp. $N=2n+1$ odd). Its lifting is given by 
$$\pi(a,\lambda|_{\mu(F)_1})\mapsto \tilde\pi( a , \tilde\lambda|_{\mu_F},(-1)^{n-1}\lambda(-1)),$$ where  $\tilde\lambda|_{\mu_F} = \lambda|_{\mu(F)_1}\circ(1-c)$.

\subsubsection{Odd special orthogonal groups.} 
Let $G $ be the split $\SO_{2n+1}$ with symmetric matrix $H = \antidiag(1,-1,1\dots,-1,1)$. The simple affine roots of $G$ are $\alpha_i = e_{i}-e_{i+1}$ for $1\leq i\leq n-1$, $\alpha_n = e_n$, and $\alpha_0 = 1-\alpha_l = 1-(e_1+e_2)$. Here the root subgroup of $\alpha_0$ occupies $(2n,1)$ and $(2n+1,2)$-entries, i.e., an element $u\in \mathcal{I}^+$ is of the form 
$$u = I+\begin{bmatrix}
&u_1&&
\\
&\cdot&\diag( u_2, \dots, u_{n},u_n,\dots,u_2)&
\\
u_0\varpi&\cdot&\cdot&u_1
\\
&u_0\varpi&&
\end{bmatrix}.
$$
 The character $\lambda$ maps $u_i\mapsto  \psi_{a_i}(u_i)$ for all $i\in \{0,1,\dots,n\}$, where $a_i\in \mathfrak{o}_F^\times$. The normalizer is $N(\lambda)= \Omega \mathcal{I}^+$. The order of $\bar\Omega$ is 2 and the non-trivial element in $\bar\Omega$ is represented by $\omega =  -\mathbbm p$, where
   $$ \mathbbm p= \antidiag ((a_0/a_1)\varpi^{-1},I_{2n-1},(a_1/a_0)\varpi) , $$
so that $\det \omega = 1$ and $\omega^2 =I$. Take a sign $\xi\in \{\pm 1\}$ and define the extended character
$$\lambda(\omega^j  u) = \xi^j\prod_{i}\psi_{a_i} (u_i),\,\qquad \text{ for }j\in \mathbb Z/2\text{ and } u\in \mathcal{I}^+.$$
Put $a = a_{0}a_{1}\prod_{i\geq 2} a_{i}^2\bmod \mathcal{U}^1(F)$. Two representations $\pi_{\lambda_1}$ and $\pi_{\lambda_2}$, where $\lambda_j = \lambda(\{a_{i}^j\}_i,\xi_j)$, are isomorphic if and only if $(a^1,\xi_1) = (a^2,\xi_2)$.

From  \cite{Oi-SO-odd}, a representative takes the form  $(a,\lambda(\omega))\in \mathbb F^\times \times \{\pm 1\} $, where $a$ represents $(a_0,\dots,a_{n}) = (a,1,\dots,1)$, and its lifting is 
\begin{equation}
\label{Oi's lifting odd orthogonal}
\pi(a,\lambda(\omega)) \mapsto \tilde \pi(2a,\mathbf 1_{\mu_F}, \lambda(\omega)).
\end{equation}

\subsubsection{\bf Symplectic groups} 
Let $G = \SP_{2n}$ be defined by the involution $g\mapsto H{}^tg^{-1}H^{-1} $, where $H = \antidiag(1,-1,1,-1,\dots)$. The simple affine roots are $\alpha_i = e_i - e_{i+1}$ for  $1\leq i\leq n-1$, $\alpha_n = 2e_n$, and $\alpha_0 = 1-\alpha_l = 1-2e_1$. Here the root subgroup of $\alpha_0$ occupies the $(2n,1)$-entry, and so 
$$u  =I+\antidiag( 
\diag( u_1, \dots,  u_{n-1}, u_n, u_{n-1}, \dots, {u_1}), 
u_0\varpi )\in \mathcal I^+.$$
 Since $G$ is simply connected, the group $\Omega$ is trivial, and so the normalizer $N(\lambda) $ is $Z \mathcal{I}^+ = \{\pm 1\} \mathcal{I}^+$. We take a sign $\xi\in \{\pm 1\}$
and define the extended character 
$$\lambda((-1)^j  u) = \xi^j\prod_{i}\psi_{a_i} (u_i),\,\qquad \text{ for }j\in \mathbb Z/2\text{ and } u\in \mathcal{I}^+.$$
Put $a = a_{0}\left(\prod_{i=1}^{n-1} a_{i}^2\right)a_{n} \mod \mathcal{U}^1(F)$. Two representations $\pi_{\lambda_1}$ and $\pi_{\lambda_2}$, where $\lambda_j = \lambda(\{a_{i}^j\}_i,\xi_j)$ are isomorphic if and only if $(a^1,\left(\frac{a_{n}^1}{\mu_F}\right),\xi_1) = (a^2,\left(\frac{a_{n}^2}{\mu_F}\right),\xi_2)$.

By \cite[Th 7.17]{Oi-Sp-and-SO-even}, a representative takes the form $(a,\kappa,\lambda(-1))\in \mathbb F^\times \times \mathbb F^\times / \mathbb F^{\times2} \times \{\pm 1\} $, where $(a,\kappa)$ represents $(a_0,\dots,a_{n}) = (a\kappa^{-1},1,\dots,1,\kappa)$. Note that we can take $\kappa\in \{1,\zeta\}$. Both  $\pi(a\kappa^{-1},\kappa,\lambda(-1))$ lift to 
\begin{equation}
\label{Oi's result for symplectic groups}
\tilde\pi (4a,\left(\tfrac{\cdot}{\mu_F}\right),\lambda(-1)\left(\tfrac{-1}{\mu_F}\right)\mathfrak n_\psi)\oplus \chi_{F[\sqrt{(-1)^{n-1}a\varpi}]/F},
\end{equation}
where, if $E/F$ is a quadratic extension, $\chi_{E/F}$ is the character of $F^\times$ whose kernel is the image of $N_{E/F}$.

\subsubsection{Ramified even orthogonal groups}

Take ramified $G = \SO^{}_{2n,F[\sqrt\varpi]/F}$ defined by 
\begin{equation}
\label{H matrix for ramified even orthogonal group G-plus}
H = \antidiag(1,\dots,1,\diag(-\varpi, 1),1,\dots,1).
\end{equation}
The simple affine roots
are
$\alpha_i = e_{i}-e_{i+1}$, for $1< i\leq n-2$, $\alpha_{n-1}=e_{n-1}$, and $\alpha_0 = \tfrac{1}{2}-e_1$ for $n>2$. If $n=2$, then $\Delta^{\text{aff}} = \{e_1, 1-e_1\}$. We express elements in $\mathcal I^+$ as $u = u(u_0,\dots,u_{n-1})$, where
$$u  =I+
\begin{bmatrix}
&\diag(u_1,\dots,u_{n-2})&&&&
\\
&&\cdot&u_{n-1}&\cdot&
\\
u_0&\cdot&\cdot&\cdot&\cdot&
\\
&&\cdot&\cdot&-u_{n-1}&
\\
&&\cdot&&&-\diag(u_{n-2},\dots,u_1)
\\
&&u_0\varpi&&&
\end{bmatrix}.
$$
The normalizer is $N(\lambda)  = Z \mathcal{I}^+ = \{\pm 1\} \mathcal{I}^+$. Take $\xi\in \{\pm 1\}$ and define
$$\lambda((-1)^k u) = \xi^k\prod_{i=0}^{n-1}\psi_{a_i} (u_i),\,\qquad \text{ for }k=1,2\text{ and } u\in \mathcal{I}^+.$$
Put $a=\prod_{i=0}^{n-1} a_{i} \bmod \mathcal{U}^1(F)$, then two representations $\pi_{\lambda_1}$ and $\pi_{\lambda_2}$, where $\lambda_j = \lambda(\{a_{i}^j\}_i,	\xi_j)$, are isomorphic if and only if $(a^1,\xi_1 )= (a^2,\xi_2)$.

An outer automorphism of $G$ is represented by  a non-trivial element in $G^\sharp = \mathrm O_{2n}$ not in $G$. For example, we take \begin{equation}
\label{p element in ramified even orthogonal groups}
     \mathbbm p  = \diag(I_{n-1},\diag(1,-1),I_{n-1}),
     \end{equation}
which conjugates $u(u_0,\dots,u_{n-2},u_{n-1})$ to $u(u_0,\dots,u_{n-2},-u_{n-1})$.

From \cite[Th 7.16]{Oi-Sp-and-SO-even}, a representative takes the form $(a, \lambda(-1))$ where  $a\in \mathbb F^\times$ 
represents $(a_0,\dots,a_{n}) = (a,1,\dots,1)$. Both $\pi(a,\lambda(-1))$ and $\pi(-a,\lambda(-1))$, related by the outer automorphism of $G$, lift to \begin{equation}
\label{Oi's lifting ramified even orthogonal}
\tilde\pi((-1)^{n-1}a^2, \left(\tfrac{\cdot}{\mu_F}\right), \mathfrak n_\psi \lambda(-1)).
\end{equation}

\subsubsection{Split and unramified even orthogonal groups}

 Let $G=\SO_{2n}$ or $G=\SO^{ur}_{2n}$ (note that Oi considered $\SO_{2n+2}$ and $\SO^{ur}_{2n+2}$ instead). The calculations for these two types of groups are similar (since they are isomorphic over $K$), so we put them in a single subsection.

Henceforth, $G$ is an even orthogonal group defined by the form 
$$H = \antidiag(1,\dots,1,H',1,\dots,1),$$
where 
$$H'=\antidiag(1,1) \quad \text{if $G$ is split, or}\quad \diag(-\zeta,1)\quad\text{if $G$ is unramified}.$$

We assume that $n\geq 3 $, then the simple affine roots are 
\begin{equation*}
   \begin{split}
    &\text{$\alpha_i = e_i - e_{i+1}$ for $1\leq i\leq n-1$, $\alpha_n = e_{n-1} +e_n  $, and $\alpha_0 = 1-\alpha_l = 1-(e_1+e_2)$}, 
    \text{ or }\\
       &\text{$\alpha_i = e_{i}-e_{i+1}$, for $1< i\leq n-2$, $\alpha_{n-1}=e_{n-1}$, and $\alpha_0 = 1-(e_1+e_2)$,}\end{split}
\end{equation*}
according to whether $G$ is split or unramified.

We express an element in $\mathcal I^+$ as
$$u  =I+
\begin{bmatrix}
&u_1&&
\\
&&\diag(u_2,\dots,u_{n-2}, u', -u_{n-2},\dots,-u_2)&
\\
u_0\varpi&\cdot&&-u_1
\\
&-u_0\varpi&&
\end{bmatrix}. 
$$
where 
$$u'=\begin{bmatrix}
u_{n-1}&u_{n}&
\\
&&-u_{n}
\\
&&-u_{n-1}\end{bmatrix} \quad \text{if $G$ is split, or}\quad \begin{bmatrix}
u_{n}\zeta&u_{n-1}&
\\
&&u_{n}
\\
&&-u_{n-1}
\end{bmatrix}\quad\text{if $G$ is unramified}. $$
The order of $\Omega$ is 2, whose non-trivial element can be represented by  
\begin{equation*}
   \omega = \antidiag(-\tfrac{a_0}{a_1}\varpi^{-1}, \diag(I_{n-2},\omega',I_{n-2}),-\tfrac{a_1}{a_0}\varpi)
  ,\end{equation*}
where 
\begin{equation*}
   \begin{split}
     \omega'&=\antidiag\left(\frac{a_{n-1}}{a_n},\frac{a_n}{a_{n-1}}\right) \quad \text{if $G$ is split, or}
     \\
     \quad& \frac{1}{a_n^2-a_{n-1}^2\zeta}\begin{bmatrix}
a_n^2+a_{n-1}^2\zeta& -2a_na_{n-1}
\\
2a_na_{n-1}\zeta & -(a_n^2+a_{n-1}^2\zeta)
\end{bmatrix}\quad\text{if $G$ is unramified},   \end{split}
\end{equation*}
 so that $\omega^2 =I$. Take two signs $\eta, \xi\in \{\pm 1\}$
and define
$$\lambda(\omega^j (-1)^k u) = \eta^j \xi^k\prod_{i}\psi_{a_i} (u_i),\,\qquad \text{ for }j,k\in \{1,2\}\text{ and } u\in \mathcal{I}^+,$$
i.e., the affine generic character can be extended to the full normalizer $N(\chi) = \Omega \{\pm 1\} \mathcal{I}^+$. Put 
\begin{equation}
\label{a element in split or unramified orthogonal groups}
 a=a_{0}a_{1}\left(\textstyle\prod_{i=2}^{n-2} a_{i}^2\right)N_{a_{n-1},a_n},
\end{equation}
where 
 $$N_{a_{n-1},a_n}=a_{n-1}a_n\quad \text{if $G$ is split, or}\quad a_{n-1}^2-a_{n}^2\zeta^{-1}\quad\text{if $G$ is unramified}.$$
Two representations $\pi_{\lambda_1}$ and $\pi_{\lambda_2}$, where $\lambda_j = \lambda(\{a_{i}^j\}_i,\xi_j,\eta_j)$, are isomorphic if and only if $(a^1,\left(\tfrac{N_{a^1_{n-1},a^1_n}}{\mu_F}\right)$, $\xi_1,\eta_1)=(a^2,\left(\tfrac{N_{a^2_{n-1},a^2_n}}{\mu_F}\right),\xi_2,\eta_2)$.

An non-trivial element in $G^\sharp = \mathrm O_{2n}$ not in $G$ can be represented by  
\begin{equation}
\label{p element in split or unramified orthogonal groups}
     \mathbbm p  = \diag(I_{n-1},\omega',I_{n-1}).
     \end{equation}
Whether $G$ is split or unramified,  $\mathbbm p$ stabilizes the character $\lambda$.

To express Oi's lifting result, we first put
 $$u_G=0 \quad \text{if $G$ is split, or}\quad 1\quad\text{if $G$ is unramified}.$$
Depending on whether $G$ is split or unramified, a representative takes the form 
\begin{equation*}
(a,\kappa,\lambda(-1),\lambda(\omega))\in 
   \begin{cases}
     \mathbb \mu_F \times \mu_F / \mu_F^{2} \times \{\pm 1\} \times \{\pm 1\} , \\
        \mathbb \mu_F \times \mu_{F[\sqrt\zeta]} / (\mu_F \ker N_{F[\sqrt\zeta]/F} )\times \{\pm 1\} \times \{\pm 1\} , \\
   \end{cases}
\end{equation*}
where $(a,\kappa)$ represents $(a_0,\dots,a_{n}) = (a\kappa^{-1},1,\dots,1,\kappa)$ or $(a(N_{F[\sqrt\zeta]/F}\kappa)^{-1},1,\dots,1,\kappa)$, and we can take $\kappa \in \{1,\zeta\}$ or $\in \{1,\zeta'\}$ where $N_{F[\sqrt\zeta]/F}\zeta' = \zeta$.

From  \cite[Th 8.7]{Oi-Sp-and-SO-even}, both $\pi( a,\kappa,\lambda(-1),\lambda(\omega))$, where $\kappa \in \{1,\zeta\}$ or $\in \{1,\zeta'\}$, lift to 
\begin{equation}
\label{Oi's lifting split or unramified even orthogonal}
\tilde\pi((-1)^{n}2^{2-u_G}a, \left(\tfrac{\cdot}{\mu_F}\right), 
(-1)^{u_G}\lambda(\omega)\left(\tfrac{-1}{\mu_F}\right)\mathfrak n_\psi )
\times 
\tilde\chi_1\times 
\tilde\chi_{2},
\end{equation}
where both $\tilde\chi_1$ and  $\tilde\chi_2$ are tamely ramified characters of $F^\times$, and
\begin{equation}
\label{Oi's result for split and unramified even orthogonal groups}
   \begin{split}
  \tilde \chi_1 \equiv \mathbf 1_{\mu_F}, \quad \tilde\chi_1(\varpi) = \lambda(\omega_o);\quad \tilde \chi_2 \equiv \left(\frac{\cdot}{\mu_F}\right), \quad 
     \tilde\chi_2(\varpi) = (-1)^{u_G}\left(\frac{-2^{u_G}a}{\mu_F}\right)\lambda(\omega_o).
         \end{split}
\end{equation}

\subsection{Comparison with our results}

We use our reducibility results in Propositions \ref{first general form of lifting, unitary group} and \ref{first general form of lifting} to determine the endoscopic liftings of simple supercuspidals, and show that our results are the same as Oi’s. The groups considered in Oi's are all quasi-split, so that $I_\zeta = \O$ in the simple supercuspidal case.

\subsubsection{Unramified unitary groups}
\label{subsubsection Unramified unitary groups}

Our result in Proposition \ref{first general form of lifting, unitary group} shows that $\tilde{\boldsymbol{\lambda}}(\varpi_i) = \lambda(-1) \mathfrak n(\mathbf q) $ where where $\mathbf q = \mathbf q_{z,s}$ in Section \ref{subsection Non-degeneracy of a quadratic form}, while Oi's is $\tilde{\boldsymbol{\lambda}}(\varpi_i) = (-1)^{n-1}\lambda(-1)$ where $n=\dim_FV$. Hence our aim is to show that 
$$\mathfrak n(\mathbf q) = (-1)^{n-1}.$$ 
 Let $\zeta\in \mu_{\Fo}$ be a generator, so that $F = \Fo[\sqrt\zeta]$ with $\overline{\sqrt\zeta} = -\sqrt\zeta$. As a quadratic form over $\mathbb F_\bullet$, the discriminant of $\mathbf q$ is equal to $(-\zeta)^{n-1}N_{\mathbb{F/F}_\bullet}(\det{}_{\mathbb F} \mathbf q)$, where $\det{}_{\mathbb F} \mathbf q$ is the discriminant of $\mathbf q$ over $\mathbb F$. Proposition \ref{properties of quadratic forms}\ref{Hermitian symmetry of the quadratic form} implies that 
 $\det{}_{\mathbb F} \mathbf q\in \mathbb F_\bullet$. Hence 
$$\mathfrak n(\mathbf q) =\left(\frac{(-\zeta)^{n-1}(\det{}_{\mathbb F} \mathbf q)^2}{\mu_\Fo}\right)\mathfrak n_\psi^{2n-2} = \left(\frac{\zeta}{\mu_\Fo}\right)^{n-1} = (-1)^{n-1},$$
which is the desired result. We remark that the result applies to both forms $G_+$ and $G_-$.

\subsubsection{Ramified even orthogonal groups}
\label{subsubsection Ramified even orthogonal groups}

The ramified even orthogonal case is much simpler than the split and unramified even orthogonal case, and so we beginning from there.

Let $G = G(V,h)$ be an even ramified even orthogonal group. Here the Hermitian matrix is $H = H_+$ defined in (\ref{H matrix for ramified even orthogonal group G-plus}). The lattice sequence determined by 
\begin{equation*}
   \begin{split}
     &\Lambda(\tfrac{1}{2n}-\tfrac{1}{2}) = [\mathfrak{o}_F^{\oplus 2n}],\quad \dots,
     \quad 
     \Lambda(0) = [\mathfrak{o}_F^{\oplus n+1},\mathfrak{p}_F^{\oplus n-1}],
     \quad
     \Lambda(\tfrac{1}{2n}) = [\mathfrak{o}_F^{\oplus n},\mathfrak{p}_F^{\oplus n}]
     \\
& \Lambda(\tfrac{1}{n}) = [\mathfrak{o}_F^{\oplus n-2 },\mathfrak{p}_F,\mathfrak{o}_F,\mathfrak{p}_F^{\oplus n }],\quad \dots,
\quad \Lambda(\tfrac{1}{2}) = [\mathfrak{p}_F^{\oplus n-1},\mathfrak{o}_F,\mathfrak{p}_F^{\oplus n}],
    \end{split}
\end{equation*}
is self-dual. The affine generic functional $\beta = \beta(a_0,\dots,a_{n-2}, a_{n-1}) $ in this case is represented by 
$$ \frac{1}{2}\begin{bmatrix}
&&a_0&&&
\\
\diag(a_1,\dots,a_{n-2})&&\cdot&&&
\\
&\cdot&\cdot&\cdot&\cdot&a_0\varpi^{-1}
\\
&a_{n-1}&\cdot&\cdot&&
\\
&\cdot&\cdot&-a_{n-1}&&
\\
&&&&\diag(-a_{n-2},\dots,-a_{1})&
\end{bmatrix},$$ 
then $[\Lambda_{},0,2\beta]$ is equivalent to $[ \Lambda_{},0,\tilde \beta]$ where 
\begin{equation}
\label{beta in ramified even orthogonal group}
\tilde\beta  = \varpi_i^{-1} = \antidiag((-1)^{n-1}a^2\varpi^{-1},I_{2n-1})
\end{equation}
and $a=\prod_{i=0}^{n-1} a_{i} \bmod \mathcal{U}^1(F)$. If we take $\mathbbm p = \diag(I_{n-1},1,-1,I_{n-1})$ as in (\ref{p element in ramified even orthogonal groups}), then 
\begin{equation}
\label{action of p on beta, ramified even orthogonal}
\Ad(\mathbbm p)\beta(a_0,\dots,a_{n-2}, a_{n-1}) = \beta(a_0,\dots,a_{n-2}, -a_{n-1}).
\end{equation}
By comparing Oi's result (\ref{Oi's lifting ramified even orthogonal}) with ours in Proposition \ref{first general form of lifting}, we have to show that 
$$\mathfrak n(\mathbf q) = \left(\frac{-2}{\mu_F}\right)\mathfrak n_\psi.$$
The bilinear form associated to $\mathbf q$ is 
$$\mathbf h:(X,X')\mapsto  \tfrac{1}{2}\tr_{\tilde{\mathfrak g}/F}(( X - \tilde \beta^{-1} X\tilde \beta){}^\alpha X'),$$ 
with $X=\diag(x_1,\dots,x_{2n})\in \tilde{\mathfrak g} $ mod center,
is expanded into 
\begin{equation}
\label{quadratic form in the even ramified orthogonal case}
   \begin{split}
    \mathbf h:(X,X')\mapsto &\tfrac{1}{2}(
(x_n-x_{2n})x_1'  +(x_{1}-x_{n})x'_{n}
+(x_{n+2}-x_{n+1})x'_{n+1}+
(x_{n+1}-x_{n-1})x'_{n+2}
\\
&+
\left(\sum_{i=2}^{n-1} +\sum_{i=n+3}^{2n}\right) (x_{2n+2-i}-x_{2n+1-i})x_i' ).
\end{split}
\end{equation}
Modulo the 1-dimensional radical, i.e., the center of $\tilde{\mathfrak g}$, it has discriminant $2^{2n-1}(-1)^{n}$. Hence 
$$\mathfrak n(\mathbf q) =  \left(\frac{2^{2n-1}(-1)^{n}}{\mu_F}\right)\mathfrak n_\psi^{2n-1}  = \left(\frac{-2}{\mu_F}\right)\mathfrak n_\psi$$
which is the desired result.

\begin{rmk}
The above result applies to $G_+$ defined by the Hermitian form (\ref{H matrix for ramified even orthogonal group G-plus}). The result for $G_-$  is just very similar, in which case we take $H = \antidiag(1,\dots,1,\zeta\diag(-\varpi, 1),1,\dots,1)$, so that an affine unipotent element is of the form
$$u  =I+
\begin{bmatrix}
&\diag(u_1,\dots,u_{n-2})&&&&
\\
&&\cdot&u_{n-1}&\cdot&
\\
u_0&\cdot&\cdot&\cdot&\cdot&
\\
&&\cdot&\cdot&-u_{n-1}\zeta^{-1}&
\\
&&\cdot&&&-\diag(u_{n-2},\dots,u_1)
\\
&&u_0\zeta\varpi&&&
\end{bmatrix},
$$
and the affine generic functional $\beta $ is represented by 
$$ \frac{1}{2}\begin{bmatrix}
&&a_0&&&
\\
\diag(a_1,\dots,a_{n-2})&&\cdot&&&
\\
&\cdot&\cdot&\cdot&\cdot&a_0\zeta^{-1}\varpi^{-1}
\\
&a_{n-1}&\cdot&\cdot&&
\\
&\cdot&\cdot&-a_{n-1}&&
\\
&&&&\diag(-a_{n-2},\dots,-a_{1})&
\end{bmatrix}.$$ 
The rest of the calculation is just similar to the case when $G=G_+$.
\qed\end{rmk}

\subsubsection{Odd orthogonal groups}
\label{subsubsection Odd orthogonal groups}

The affine generic functional $\beta$ defining a character of $\mathcal I^+$ is represented by $$\tfrac{1}{2}\begin{bmatrix}
&&a_0\varpi^{-1}&
\\
a_1&&\cdot &a_0\varpi^{-1}
\\
&\diag(a_2,\dots,a_{n},a_n,\dots,a_2)&&
\\
&&a_1&
\end{bmatrix}.$$
Put $a = a_0a_1\prod_{i\geq 2}a_i^2$ and define, with $I=\{i,o\}$, $$\beta_i = \varpi_i^{-1}  = \antidiag(a(2^{2n-3}\varpi)^{-1},I_{2n-1})\in \tilde {\mathfrak g}(V_{i}),$$ 
and take $\beta_o=0\in F$. The stratum $[\Lambda_{}, 0,\beta]$ is then equivalent to $[\Lambda_i\oplus \Lambda_o,0,\beta_i\oplus\beta_o]$ (note that  $\beta_i\oplus\beta_o$ has the same characteristic polynomial as $\beta$), and we take $\tilde \beta =2\beta_i$ to define $\tilde{\boldsymbol{\lambda}}_i|_{\mathcal U^{0_+}(\Lambda_i)}$.

To fit into the calculation (\ref{beta in ramified even orthogonal group}) of the even orthogonal group case, we modify our uniformizer by defining $\varpi'$ such that 
$$a(2^{2n-3}\varpi)^{-1}=(-1)^{n-1}b^2\varpi'^{-1} \mod \mathcal U^1(F) $$
 for some $b\in \mu_F$, and take \begin{equation*}
\label{odd orthogonal group, H1 and Ho}
H_i = \antidiag(1,\dots,1,\diag(-\varpi', 1),1,\dots,1),\quad H_o=[\varpi'].
\end{equation*}
Comparing with (\ref{Oi's lifting odd orthogonal}), our aim is to show that $\tilde{\boldsymbol{\lambda}}(\varpi_i) = \lambda (\omega)$, where $\omega\in G(V,h)$ corresponds to $\omega_i = \diag(-I_{2n},1)\in G(V,h_{H_i\oplus H_o})$. It suffices to prove that 
$$\mathfrak{n}(\mathbf q)=1,$$
where, with $H_i\oplus H_o$ defined as above, we have $\mathbf q = \mathbf q_{i}\oplus \mathbf q_{o}$, where $\mathbf q_{i}$ is defined as in (\ref{quadratic form in the even ramified orthogonal case}), and $\mathbf q_{o}=\tfrac{1}{2}x_0x_0'$. 
Here $X^i = \diag(x_1,\dots,x_{2n})\in\tilde {\mathfrak g}(V_{i})$ and $x_0\in \End(V_o,V_{i})$ is a column vector with almost all entries zero except possibly one $x_0$. The discriminant of $\mathbf q_{i}$ was computed in the even ramified orthogonal group case in the previous Subsection \ref{subsubsection Ramified even orthogonal groups}, which is $2^{2n-1}(-1)^{n}$ (with the modified uniformizer $\varpi'$), while that of $\mathbf q_{o}$ is $-2$. Hence the discriminant of $\mathbf q$ is $2^{2n-1}(-1)^{n}\cdot2 \equiv (-1)^{n}\bmod \mathbb F^{\times2}$, and so 
$$ \mathfrak{n}(\mathbf q) = \left(\frac{(-1)^{n}}{\mu_F}\right)\mathfrak n_\psi^{2n}=1$$
which is the desired result.

\subsubsection{Symplectic groups}
\label{subsubsection Symplectic groups comparison}

Let's look at $i\neq o$ first. The affine generic functional $\beta $ is represented by 
$$\tfrac{1}{2}\antidiag(2a_0\varpi^{-1}, \diag(a_1,\dots,a_{n-1},2a_n, a_{n-1},\dots,a_1)).$$
Put $a = a_0(\prod_{i=1}^{n-1}a_i^2)a_n$, then stratum $[\Lambda_{},0,2\beta]$ is hence equivalent to $[\Lambda,0,\tilde\beta]$ where 
$$\tilde\beta = \varpi_i ^{-1} = \antidiag(a(2^{2n-2}\varpi)^{-1},I_{2n-1}).$$

Our results in Proposition \ref{first general form of lifting} implies that the lifting is $\tilde\pi (4a,\left(\tfrac{\cdot}{\mu_F}\right), \left(\tfrac{-2}{\mu_F}\right)
\lambda (-1)\mathfrak{n}(\mathbf q))$. Compare this with (\ref{Oi's result for symplectic groups}), it suffices to show that 
$$\mathfrak{n}(\mathbf q) = \left(\frac{2}{\mu_F}\right)\mathfrak n_\psi.$$
The bilinear form associated to $\mathbf q$ is 
$$(X,X')\mapsto  \tfrac{1}{2}\tr_{\tilde{\mathfrak g}/F}(( X - \tilde \beta^{-1} X\tilde \beta){}^\alpha X').$$ 
With $X=\diag(x_1,\dots,x_{2n})\in \mathbb F^{\oplus 2n}\bmod \Delta\mathbb F$, it is expanded into
$$ \frac{1}{2}\left( (x_{1}-x_{2n})x_1' 
+\sum_{i=2}^{2n}(x_{2n+2-i}-x_{2n+1-i})x_{i}'\right).$$
Modulo the 1-dimensional radical, its discriminant is $2^{2n-1}(-1)^{n-1} = 2(-1)^{n-1}\bmod(\mathbb F^\times)^2$. The associated normalized Gauss sum is 
\begin{equation*}
\label{Gauss sum appearing in tripling method, SP-case}
\mathfrak{n}(\mathbf q)=\left(\frac{2(-1)^{n-1}}{\mu_F}\right)
\mathfrak{n}_\psi^{2n-1}  = \left(\frac{2}{\mu_F}\right)\mathfrak{n}_\psi
\end{equation*}
which is the desired result.

\subsubsection{Symplectic groups with $i=o$}
\label{subsubsection Symplectic groups comparison, i=o}

Now look at $i=o$.  Oi's suggested character is the central character $\tilde\chi$ of the above lifting, which is trivial on ${\mathcal U^1(F)}$, with 
$$\tilde\chi|_{\mu_F} = \left(\frac{\cdot}{\mu_F}\right)
\quad \text{and} \quad
\tilde\chi(\varpi) = \tilde{\boldsymbol{\lambda}}((4a)^{-1}\varpi_i^{2n}) =  \left(\frac{a}{\mu_F}\right)\mathfrak n_\psi^{2n} = \left(\frac{(-1)^na}{\mu_F}\right).
$$
This is the quadratic character associated to the ramified quadratic extension $F[\sqrt{(-1)^{n-1}a\varpi}]/F$.

Recall that our Hermitian space is defined by the matrix $H = \antidiag(1,-1,1,-1,\dots)$, and also in this case $\dim \tilde V  =1$. The operator $\alpha$ on $X$ is 
$$[x_1,\dots,x_{2n}] \mapsto {}^t[x_{2n},-x_{2n-1},\dots,x_{2},-x_1],$$
and on $Y$ is identity, so that the relation $X{}^\alpha X = Y- {}^\alpha Y$ becomes $0=0$, i.e., there is no relation between $X$ and $Y$.

The maximal lattice sequence $\Lambda$ in $V$ determined by 
$$\Lambda(0) = [\mathfrak{o}_F^{\oplus n},  \mathfrak{p}_F^{\oplus n} ],\quad  \Lambda(\tfrac{1}{2n})= [\mathfrak{o}_F^{\oplus n-1},  \mathfrak{p}_F^{\oplus n+1} ],\cdots$$
is self-dual, and we take 
$$\Lambda_- (0) = \mathfrak{o}_F
\quad\text{and}\quad
 \Lambda_- (0_+) = \Lambda_- (1) = \mathfrak{p}_F.$$
We take the null stratum $[\{\mathfrak{p}^k\}_k,0,0]$ in $F$. Then the group ${\mathcal J}_P$ can be expressed as
\begin{equation*}
\label{JP-group, with dim V=1, symplectic case}
\begin{bmatrix}
\mu_F+ \mathfrak{p}_F&[\mathfrak{p}_F^{\oplus n},  \mathfrak{o}_F^{\oplus n} ]&\mathfrak{o}_F
\\
{}^t [\mathfrak{p}_F^{\oplus (n+1)},  \mathfrak{o}_F^{\oplus (n-1)} ]&\mu_F + \mathcal{I}^+&{}^t[\mathfrak{o}_F^{\oplus n},  \mathfrak{p}_F^{\oplus n} ]
\\
\mathfrak{p}_F&[\mathfrak{p}_F^{\oplus (n+1)},  \mathfrak{o}_F^{\oplus (n-1)} ]&\mu_F + \mathfrak{p}_F
\end{bmatrix}.
\end{equation*}
We now compute the coefficients $b_w$. For $w=y$, we have
$$X\in[\mathfrak{p}_F^{\oplus n},  \mathfrak{o}_F^{\oplus n} ]/[\mathfrak{p}_F^{\oplus (n+1)},  \mathfrak{o}_F^{\oplus (n-1)} ]\cong \mathfrak{o}_F/\mathfrak{p}_F
\qquad\text{and}\qquad
Y\in \mathfrak{o}_F^\times/\mathfrak{p}_F.
$$
Here $X$ is represented by $[0^n,x_{n+1},0^{n-1}]$, then $-\beta{}^\alpha X Y^{-1} X = (-1)^na_nx_{n+1}^2Y^{-1}$, so that
$$b_y \tilde T_y(s_y)= \sum_{Y}\tilde{\lambda}(Y) \sum_{x_{n+1}}\psi((-1)^na_nx_{n+1}^2Y^{-1})=\left(\frac{(-1)^n a_n}{\mu_F}\right)q^{1/2}\mathfrak n_\psi(q-1)$$
when $\tilde\lambda |_{\mu_F} $ is quadratic, in which case \begin{equation*}
c_y = q^2, \quad r_y = 1
\quad\text{and}\quad 
\epsilon_y T_y(s_y)= \left(\frac{(-1)^n a_n}{\mu_F}\right)\mathfrak n_\psi.
\end{equation*}
The calculation for $w=z$ is similar; we have
\begin{equation*}
\label{definition of P_(F,Lambda)}
X\in \varpi^{-1}[\mathfrak{p}_F^{\oplus (n+1)},  \mathfrak{o}_F^{\oplus (n-1)} ]/[\mathfrak{p}_F^{\oplus n},  \mathfrak{o}_F^{\oplus n} ] \quad\text{and}\qquad
\varpi^{-1}Y_0\in \mathfrak{p}_F^{-1}/\mathfrak{o}_F,
\end{equation*}
with $X$ represented by $[x_1,\dots,x_{n},0,x_{n+2}\varpi^{-1},\dots,x_{2n}\varpi^{-1}]$, and 
$$\tr(-\beta{}^\alpha X Y^{-1} X )= (a_0x_1^2+2(-a_1x_2x_{2n}+\cdots+(-1)^{n-1}a_{n-1}x_nx_{n+2}) )Y_0^{-1},$$ so that
$$b_z \tilde T_z(s_z)= \sum_{Y_0}\tilde{\lambda}(Y_0) \sum_{X}\psi(\tr(-\beta{}^\alpha X Y^{-1} X ))=q^{(2n-1)/2}(q-1)\left(\frac{a/a_n}{\mu_F}\right)\mathfrak n_\psi .$$
again when $\tilde\lambda |_{\mu_F} $ is quadratic, in which case 
\begin{equation*}
c_z = q^{2n}, \quad r_z = 1
\quad\text{and}\quad 
\epsilon_z T_z(s_z)= \left(\frac{a/a_n}{\mu_F}\right)\mathfrak n_\psi.
\end{equation*}
Therefore, 
\begin{equation}
\label{symplectic groups, the character in the simple case, i=o}
\tilde\lambda|_{\mu_F}  = \left(\frac{\cdot}{\mu_F}\right)
\quad\text{and}\quad
\tilde{\boldsymbol{\lambda}}(\varpi)= \left(\frac{(-1)^{n}a_n(a/a_n)}{\mu_F}\right)
\mathfrak{n}_\psi^{2} =\left(\frac{(-1)^{n-1}a}{\mu_F}\right) = \left(\frac{\varpi\det \beta}{\mu_F}\right)
\end{equation}
which is the character appearing in (\ref{Oi's result for symplectic groups}).

\subsubsection{Split and unramified orthogonal groups}

The affine generic functional $\beta$ is represented by $$  \frac{1}{2}\begin{bmatrix}
&&a_0\varpi^{-1}&
\\
a_1&&\cdot&-a_0\varpi^{-1}
\\
&\diag(a_2,\dots,a_{n-2}, \beta', -a_{n-2},\dots,-a_{2})&&
\\
&&-a_1&
\end{bmatrix},$$
where 
$$\beta'=\begin{bmatrix}
a_{n-1}&&\\
a_{n}&&
\\
&-a_{n}&-a_{n-1}
\end{bmatrix} \quad \text{if $G$ is split, or}\quad \begin{bmatrix}
a_{n}&&
\\
a_{n-1}&&
\\
&a_{n}\zeta^{}&-a_{n-1}
\end{bmatrix}\quad\text{if $G$ is unramified}. $$
The outer automorphism  $\mathbbm p $ given in (\ref{p element in split or unramified orthogonal groups}) acts on $\beta$ by switching the entries $a_{n-1}$ and $a_n$ if $G$ is split 
and by conjugation $(a_{n-1},a_n) \mapsto (a_{n-1},-a_n) $ if $G$ is unramified. The characteristic polynomial of $\beta$ is 
$$T^{2n}-(-1)^{n}a(2^{2-u_G} \varpi)^{-1}T^2,$$
where $a$ is given in (\ref{a element in split or unramified orthogonal groups}). Hence $[\Lambda_{},0,2\beta]$ is equivalent to $[\Lambda_{i}\oplus \Lambda_o,0,\tilde\beta\oplus (0,  0)]$, where $\Lambda_{i}$ and $\Lambda_{o}$ are self-dual lattice sequences in $V_{i}$ and $V_{o}$ respectively with $\dim V_{i}=2n-2$ and $\dim V_o = 2$, and 
$$\tilde\beta = \varpi_i ^{-1} = \antidiag((-1)^{n}a(2^{2-u_G}\varpi)^{-1},I_{2n-3}).$$

As in the odd orthogonal group case in Subsection \ref{subsubsection Odd orthogonal groups}, to fit into the calculation (\ref{beta in ramified even orthogonal group}) of the even ramified orthogonal group case in Subsection \ref{subsubsection Ramified even orthogonal groups}, we modify our uniformizer by defining 
\begin{equation}
\label{Modified uniformizer}
(-1)^{n}a(2^{2-u_G} \varpi)^{-1}=(-1)^{n-1}b^2\varpi'^{-1}\mod \mathcal U^1(F) 
\end{equation}
for some $b\in \mu_F$, and take 
$$H_1 = \antidiag(1,\dots,1,\diag(-\varpi', 1),1,\dots,1),\quad H_o=\diag(-\zeta^{u_G}\varpi',1).$$
Here $H_o$ is put into the center of $H_1$, i.e., $H = \antidiag(1,\dots,1,\diag(-\varpi', -\zeta^{u_G}\varpi', 1,1),1,\dots,1)$.

For $i\neq o$, comparing with (\ref{Oi's lifting split or unramified even orthogonal}), it suffices to show that
$$
\mathfrak{n}(\mathbf q) = (-1)^{u_G}\left(\frac{2}{\mu_F}\right) \mathfrak n_\psi
.$$ 
We express $\mathbf q$ as the sum of two quadratic forms, $\mathbf q_{i}$ in the even ramified orthogonal case in (\ref{quadratic form in the even ramified orthogonal case}) but with $\dim_{\mathbb F} = 2n-3$, and another one $\mathbf q_{o}$ with $\dim_{\mathbb F} = 2$, i.e., the associated bilinear form is 
$(X,X')\mapsto \mathbf h_{i}(X,X')+ \mathbf h_{o}(X,X')$, where $\mathbf h_{i}$ as in (\ref{quadratic form in the even ramified orthogonal case}) and 
$\mathbf h_{o} = -\tfrac{1}{2}(\zeta^{u_G}x_0x_0'+y_0y_0')$, and so $$\mathrm{disc}(\mathbf q_{o}) = \zeta^{u_G}
\quad\Rightarrow\quad \mathfrak n(\mathbf q_{o}) = (-1)^{u_G}\left(\frac{-1}{\mu_F}\right). $$
Recall that $ \mathfrak n(\mathbf q_{i})= \left(\tfrac{-2}{\mu_F}\right) \mathfrak n_\psi$. Therefore, 
$$\mathfrak{n}(\mathbf q) =   \mathfrak n(\mathbf q_{i}) \mathfrak n(\mathbf q_{o})= (-1)^{u_G}\left(\frac{2}{\mu_F}\right) \mathfrak n_\psi$$
which is the desired result.

\subsubsection{Split and unramified orthogonal groups with $i=o$}

 We continue from the previous subsection, but with now $i=o$. The calculation is somewhat similar to Section \ref{subsubsection Example 2: ramified SO(2)}; although $\dim V_o =2$, we have to take a self-dual tamely ramified character $\tilde{\boldsymbol{\lambda}}$ of $\tilde G(F) = F^\times$, i.e., $\dim \tilde V =1$, and the null stratum $[\{\mathfrak{p}^k\}_k,0,0]$ in $F$ (but not any epipelagic stratum of degree 2).

The operator $\alpha$ on $X$ is 
\begin{equation*}
   \begin{split}
&[x_1,\dots,x_{n-2},
(x_{n-1}\varpi,x_{n}\varpi,x_{n+1},x_{n+2}),
x_{n+3},\dots,x_{2n}]
\\
&\mapsto {}^t[-x_{2n},,\dots,-x_{n+3},(\varpi'^{-1} x_{n-1},\zeta^{-u_G}\varpi'^{-1}x_{n},-x_{n+1},-x_{n+2}),-x_{n-2},\dots,-x_1],
   \end{split}
\end{equation*}
and on $Y$ is minus-identity. Here $\varpi'$ if the modified uniformizer in (\ref{Modified uniformizer}). We put 
$$\mathbbm p_o = \diag(I_{n-1},\diag(1,-1),I_{n-1})
\quad \text{and}\quad
\omega_o = \diag(I_{n-1},-I_2,I_{n-1}).$$
 When $w=y$, we represent $X = [0_{n-1},(0,x_{n+1}),0_{n-1}]$, and so $2Y = -x_{n+1}^2$. 
One can show that $I - {}^\alpha X Y^{-1}X $ is $G(F)$-conjugate to $\mathbbm p_o$, and 
$$b_y\tilde T_y(s_y) = \sum_{-2Y = x_{n+1}^2}\tilde\lambda (Y) \lambda( \mathbbm p_o^2) = \tilde\lambda(-2) (q-1),$$
which implies that
\begin{equation*}
c_y = q , \quad 
r_y = 1
\quad\text{and}\quad 
\epsilon_y T_y(s_y)= \tilde\lambda(-2).
\end{equation*}
When $w=z$, with representatives $X = [x_1,\dots,x_{n-1},
(x_{n},0),
x_{n+2}\varpi^{-1},\dots,x_{2n}\varpi^{-1}]$, we have 
\begin{equation*}
\label{relation in split or unramified even orthogonal groups, w=z, i=o}
2Y\varpi^{-1} = X{}^\alpha X =-\varpi^{-1}\zeta^{-u_G}x_n^2.
\end{equation*}
Regardless of the Hermitian form, one can show that indeed the trace of $\beta_j{}^\alpha X^jX^j$ is 0, which means that $\mathbf q_{j}$ is trivial. Now $I - {}^\alpha X Y^{-1}X $ is $G(F)$-conjugate to $\omega_o \mathbbm p_o$. Hence 
$$b_z\tilde T_z(s_z) = \sum_{2Y =-\zeta^{-u_G}x_n^2.}\tilde\lambda (Y) \lambda(\omega_o \mathbbm p_o^2) = \tilde\lambda(\zeta)^{u_G}\tilde\lambda(- 2)\lambda(\omega_o) (q-1),$$
which implies that
\begin{equation*}
c_z = q, \quad 
r_z = 1
\quad\text{and}\quad 
\epsilon_z T_z(s_z)= \tilde\lambda(\zeta)^{u_G}\tilde\lambda(- 2)\lambda(\omega_o).
\end{equation*}
We see that the deduction is independent of $\tilde\lambda|_{\mu_F}$, which gives two tamely ramified characters ${\tilde{\boldsymbol{\lambda}}}_1$ and ${\tilde{\boldsymbol{\lambda}}}_2$ for $\tilde\lambda$: 
 \begin{equation}
 \label{split or unramified orthogonal groups, the two charaacters}
   \begin{split}
      \tilde \lambda_1 \equiv \mathbf 1_{\mu_F}& \quad\Rightarrow \quad {\tilde{\boldsymbol{\lambda}}}_1(\varpi') = \lambda(\omega_o); \\
       \tilde \lambda_2 \equiv \left(\frac{\cdot}{\mu_F}\right)
       & \quad\Rightarrow \quad {\tilde{\boldsymbol{\lambda}}}_2(\varpi') = (-1)^{u_G}\lambda(\omega_o),
   \end{split}
\end{equation}
as the required results in (\ref{Oi's result for split and unramified even orthogonal groups}).

\subsection{Ramified unitary groups}
\label{subsection Ramified unitary groups}

The calculation applies equally well to ramified unitary groups, the type of groups that is not covered in Oi's series. Let $G=\mathrm U_{N}(F/\Fo)$ where $F = \Fo[\varpi]$ with $\varpi^2=-\varpi_\bullet$. 

\subsubsection{Odd case}

If $N=2n+1$, we take $H=\antidiag(1,-1,1,\dots,1)$, then an 
element $u\in \mathcal{I}^+$ is of the form 
$$u = I+\antidiag( \diag(u_1,\dots,u_{n},u_{n},\dots,u_1), u_0\varpi), \quad u_0,u_1,\dots,u_{n}\in \mu_F.$$
and the affine generic functional is 
$$\beta = \tfrac{1}{2}\antidiag(2a_0\varpi^{-1},\diag(a_1,\dots,a_{n},{a_{n}},\dots,{a_1})),\quad a_0,a_1,\dots,a_{n}\in \mu_F.$$
The normalizer is $N(\chi) =  \{\pm 1\} \mathcal{I}^+$. Take a sign $\xi\in \{\pm 1\}$
and define
$$\lambda( (-1)^k u) =  \xi^k\prod_{i}\psi_{a_i} (u_i),\,\qquad \text{ for }k=1,2\text{ and } u\in \mathcal{I}^+.$$
Put 
$ a=a_{0}(\textstyle\prod_{i=1}^{n} a_{i}^2)
$, then two representations $\pi_{\lambda_1}$ and $\pi_{\lambda_2}$, where $\lambda_j = \lambda(\{a_{i}^j\}_i,\xi_j)$ are isomorphic if and only if $(a^1,\xi_1)=(a^2,\xi_2)$.

The quadratic form $\mathbf q = \mathbf q_{z,\mathbf s}$ is on the space of dimension $2n$ and has discriminant $(-1)^n$. The normalized Gauss sum $\mathfrak n(\mathbf q)=\left(\tfrac{(-1)^n}{\mu_F}\right)\mathfrak n_\psi^{2n}=1$. Hence if we take $\varpi_i = \antidiag(I_{2n},2a\varpi^{-1})$, 
\begin{equation*}
\tilde \lambda|_{\tilde {\mathcal I}^+} =\psi_{2\beta} , \quad \tilde \lambda|_{\mu_{F}}= \mathbf 1_{\mu_F}
\quad\text{and}\quad
\tilde{\boldsymbol{\lambda}}(\varpi_i)= \tilde\lambda(-2)\lambda(-1),
\end{equation*}
then the lifting is given by 
$$\pi(a,\lambda(-1)) \mapsto \tilde \pi(2a,\mathbf 1_{\mu_F}, \lambda(-1)).$$ 

\subsubsection{Even case}

If $N=2n$ is even, suppose that $H=\antidiag(1,-1,1,\dots,-1)$, then an 
element $u\in \mathcal{I}^+$ is of the form 
$$u = I+\antidiag(\diag(u_1,\dots,u_{n},\dots,u_1),u_0\varpi_\bullet)
$$
and 
$$\beta = 
\tfrac{1}{2}\antidiag(2a\varpi_\bullet^{-1}, \diag(a_1,\dots,a_{n-1},2a_n,a_{n-1},\dots,a_1)).$$
The normalizer is $N(\chi) = \{\pm 1\}  \mathcal{I}^+$ Take a sign $\xi\in \{\pm 1\}$
and define
$$\lambda( (-1)^k u) =  \xi^k\prod_{i}\psi_{a_i} (u_i),\,\qquad \text{ for }k=1,2\text{ and } u\in \mathcal{I}^+.$$
Put 
$ a=a_{0}(\textstyle\prod_{i=1}^{n-1} a_{i}^2)a_n
$, then two representations $\pi_{\lambda_1}$ and $\pi_{\lambda_2}$, where $\lambda_j = \lambda(\{a_{i}^j\}_i,\xi_j)$ are isomorphic if and only if $(a^1,\xi_1)=(a^2,\xi_2)$.

The characteristic polynomial of $\beta$ is $T^{2n}-a(2^{2n-2}\varpi)^{-1}$. We hence take $I=\{j,o\}$ and the Hermitian forms $H_j=\varpi\antidiag(1,-1,1,\dots,1)$ and $H_o=[\varpi]$. Put $\beta_j = \tfrac{1}{2}\antidiag(a\varpi^{-1},I_{2n-2})$ and $\beta_o=0$, and define $\tilde\beta = 2[\beta_j,\beta_o]$.

If $i=j$, then $\mathbf q_{i}$ is the same as in the odd case (but with $n$ replaced by $n-1$) and $\mathbf q_{o}= -x_o^2$. We hence have $\mathfrak n(\mathbf q)=\left(\tfrac{-1}{\mu_F}\right)\mathfrak n_\psi^{}$, and define
\begin{equation*}
\tilde \lambda|_{\tilde {\mathcal I}^+} =\psi_{2\beta_j} , \quad, \tilde \lambda_i|_{\mu_{F}}= \left(\frac{\cdot}{\mu_F}\right)
\quad\text{and}\quad
\tilde{\boldsymbol{\lambda}}_i(\varpi_i)= \lambda(\omega_i)\left(\frac{2}{\mu_F}\right)\mathfrak n_\psi. 
\end{equation*}
If $i=o$, then $\mathbf q_{o}=1$ as $\mathfrak W_{z,o}$ is trivial, while  
$$\mathbf q_{j}=-\varpi X^j \beta_j {}^\alpha X^j = (\det \beta_j) x_1^2, \quad 
X^j\in \mathfrak W_{z,j}.$$
Therefore, we define $\tilde{\boldsymbol{\lambda}}_o$ to be the tamely ramified character of $F^\times$ with
\begin{equation*}
\tilde \lambda_o|_{\mu_{F}}= \left(\frac{\cdot}{\mu_F}\right)
\quad\text{and}\quad
\tilde{\boldsymbol{\lambda}}_o(\varpi)= \lambda(\omega_o)\left(\frac{-2\det \beta_j}{\mu_F}\right)\mathfrak n_\psi.
\end{equation*}
The lift of $\pi_\lambda$ is the parabolically induced representation 
$$\tilde\pi(a,\left(\tfrac{\cdot}{\mu_F}\right), \lambda(\omega_i)\left(\tfrac{2}{\mu_F}\right)\mathfrak n_\psi)\times \tilde{\boldsymbol{\lambda}}_o$$ 
of $\GL_{2n}(F)$.

\section{Main results on epipelagic representations}
\label{subsection Reducibility results for different classical groups}

Generalizing from the simple case, we now give the calculations on the the Gauss sums $\mathfrak{n}_z(\varpi_i,\mathbf{s},\psi,h_{})$ appearing 
in Propositions \ref{first general form of lifting, unitary group} and \ref{first general form of lifting}, leading to simpler expressions of the liftings of epipelagic representations.

We skip the discussion for unramified unitary groups, since all epipelagic representations are simple supercuspidal in this case, and so is covered in Section \ref{subsubsection Unramified unitary groups}.

\subsection{Endoscopic liftings for classical groups}
\label{subsection Liftings of epipelagic representations for classical groups}

Again, $G=G(V,h)$ is a connected classical group, and is assumed to be not of unramified unitary type. Let $\pi = \cInd_{\mathcal J}^{G(F)}\lambda$ be an epipelagic representation constructed from an epipelagic semi-simple stratum $\mathbf s = [\Lambda,0,\beta]$. Each component $\mathbf s_i =[\Lambda_i,0,\beta_i]$ is self-dual with respect to $(V_i,h_i)$, such that the orthogonal sum $\oplus_{i\in I}(V_i,h_i)$ is isometric to $(V,h)$. By choosing a basis of $V_i$, we suppose that $h_i$ is represented by a Hermitian matrix $H_i$ which is of the form $H_+$ or $H_-$ in Section \ref{subsection Cohomological classification}.

For all $j\in I\smallsetminus \{i\}$, put 
$$\gamma_j=1-\det \beta_j \det\beta_i^{-1}
\quad\text{ and }\quad \theta_i = \det H_j\det H_i^{-1},$$ 
as in (\ref{discriminants for j neq i}). Clearly, $\gamma_j\in \mu_F$ if $i\neq j$.  Also, our choices of $H_j$, for all $j\in I$, imply that $\theta_j\in \mu_F$ too. We have
$$\mathfrak n(\mathbf q_{j}) = \left(\frac{\gamma_j\theta_j}{\mu_F}\right)\mathfrak n_\psi^{2e} = \left(\frac{(-1)^e\gamma_j\theta_j}{\mu_F}\right).$$
Here $2e $ is the common degree $[E_i:F]$ for $i\in I\smallsetminus \{o\}$.

We hence obtain,   
$$\mathfrak{n}_z(\varpi_i,\mathbf{s},\psi,h_{}) =  \mathfrak{n}(\mathbf q_{i}) \mathfrak{n}(\mathbf q_{o})\prod_{j\in I\smallsetminus\{i,o\}}\mathfrak{n}(\mathbf q_{j}),\quad \text{for $i\neq o$}. $$
Here we put $\mathfrak{n}(\mathbf q_{o})=1$ if $o\notin I$, and the values of $\mathfrak{n}(\mathbf q_{i}) $ and $\mathfrak{n}(\mathbf q_{o})$ when $G$ is orthogonal or ramified unitary and $o\in I$ are computed in the simple supercuspidal case.

Finally, we put
$$\kappa_i =  \left(\frac{(-1)^{n-e}\textstyle\prod_{j\in I\smallsetminus\{ i,o\}}\gamma_j\theta_j}{\mu_F}\right),\quad \text{for $i\neq o$}.$$

Substitute these values into Proposition \ref{first general form of lifting}, we obtain the following explicit values of epipelagic characters.

\begin{enumerate}[(i)]

\item If $G$ is odd orthogonal, then \begin{equation*}
\tilde\lambda_i|_{\mu_F} = \mathbf 1 \quad\text{and}\quad{\tilde{\boldsymbol{\lambda}}}_i(\varpi_{i}) =
\lambda(\omega_i)\kappa_i .
\end{equation*}

\item If $G$ is symplectic and $i\neq o$, then 
\begin{equation*}
\tilde\lambda_i|_{\mu_F} = \left(\frac{\cdot}{\mu_F}\right) \quad\text{and}\quad{\tilde{\boldsymbol{\lambda}}}_i(\varpi_{i}) =
\lambda(\omega_i)\left(\frac{-1}{\mu_F}\right)\kappa_i \mathfrak n_\psi;
\end{equation*}
while if $i=o$, then 
$$ \tilde\lambda_o|_{\mu_F} = \left(\frac{\cdot}{\mu_F}\right)^{\#I } 
\quad\text{and}\quad
{\tilde{\boldsymbol{\lambda}}_o}(\varpi_{}) = \prod_{i\in I}\left(\frac{\varpi\det\beta_i}{\mu_F}\right),$$
as generalizing the character in the simple case in (\ref{symplectic groups, the character in the simple case, i=o}). 
  \end{enumerate}

\begin{rmk}Note that the choice of the uniformizer $\varpi$ in the last equality is unimportant, since if $\#I$ is odd and we change $\varpi$ to $\zeta\varpi$, then both sides of the last equality produce the same sign $\left(\frac{\zeta}{\mu_F}\right)^{\#I } $.
\qed\end{rmk}

   \begin{enumerate}[resume*]
\item If $G$ is $\SO_{2n}$, then $\#I$ is odd if and only if $G$ is ramified (by considering discriminant). If $o\notin I$,
 \begin{equation*}
\tilde\lambda_i|_{\mu_F} = \left(\frac{\cdot}{\mu_F}\right)^{} \quad\text{and}\quad{\tilde{\boldsymbol{\lambda}}}_i(\varpi_{i}) =
\lambda(\omega_i)
\kappa_i
 \mathfrak n_\psi^{}.
\end{equation*}

If $i=o$, then there are two extra characters ${\tilde{\boldsymbol{\lambda}}}_1$ and ${\tilde{\boldsymbol{\lambda}}}_2$ as in (\ref{split or unramified orthogonal groups, the two charaacters}), except with a modification on 
 \begin{equation*}
{\tilde{\boldsymbol{\lambda}}}_2(\varpi') = \kappa_o\lambda(\omega_o),
\quad \text{where}\quad
\kappa_o=\left(\frac{-\varpi'\det H_o^{-1}}{\mu_F}\right).
\end{equation*}

\item 

If $G$ is ramified unitary over $F/\Fo$, where $F = \Fo[\varpi]$ with $\varpi^2=-\varpi_\bullet$, then $\# I \equiv N\bmod 2$. Put $ \delta_o =1 $ if $o\in I$, and $=0$ otherwise. We computed in Section \ref{subsection Ramified unitary groups} that, if $i\neq o$, then we have
$\mathfrak n(\mathbf q_i) = 1$, as well as the discriminants 
\begin{equation*}
 \text{ $\disc(\mathbf q_o) = (-1)^{\delta_o}$\quad  and \quad $\disc(\mathbf q_j) = \gamma_j\theta_j$,\quad  for all $j\in I\smallsetminus\{i,o\}$;}
 \end{equation*}
while if $i=o$, then 
$$   
\text{$\disc(\mathbf q_j) = \det \beta_j$ \quad for all $j\in I\smallsetminus\{o\}$.
}
$$
If $i\neq o$,
\begin{equation*}
\tilde\lambda_i|_{\mu_F} =  \left(\frac{\cdot}{\mu_F}\right)^{N- 1 } 
\quad\text{and}\quad
{\tilde{\boldsymbol{\lambda}}}_i(\varpi_{i}) = \lambda(\omega_i)
\kappa_i
\mathfrak n_\psi^{N-1},
\end{equation*}
where
\begin{equation*}
\kappa_i := \left(\frac{(-2)^{N- 1 }(-1)^{\delta_o} \textstyle\prod_{j\in I\smallsetminus\{i\}}(\gamma_j\theta_j) }{\mu_F}\right).
\end{equation*}
If $i=o$, then 
\begin{equation*}
\tilde\lambda_o|_{\mu_F} =\left(\frac{\cdot}{\mu_F}\right)^{N-1 } 
\quad\text{and}\quad{\tilde{\boldsymbol{\lambda}}}_o(\varpi) =\lambda(\omega_o)
\kappa_o
\mathfrak n_\psi^{N-1},
\end{equation*}
where
\begin{equation*}
\kappa_o := \left(\frac{ (-2)^{N-1 } \textstyle\prod_{j\in I\smallsetminus\{o\}}\det\beta_j }{\mu_F}\right)
.
\end{equation*}

\end{enumerate}

\subsection{L-packets}

Finally, we provide a description of L-packets for epipelagic representations of classical groups. The calculation in each case is simply an inversion process of the previous section. Again, we skip the discussion for unramified unitary groups, since the L-packets of epipelagic supercuspidals are just singleton.

\subsubsection{Odd orthogonal groups}

For $G=\SO_{2n+1}$, given
\begin{enumerate}[(i)]

\item  a stable functional $\beta= (\beta_i)_{i\in I}$ in $G(F)$, 

\item a partition $I_\zeta$ of $I$, and

\item a tuple of signs $\delta = (\delta_i)_{i\in I\smallsetminus\{o\}}\in \{\pm 1\}^{\#I\smallsetminus\{o\}}$,

\end{enumerate}
define a character $\lambda_{(\beta,I_\zeta,\delta)}$ by 
\begin{equation*}
 \lambda_{(\beta,I_\zeta,\delta)} = \psi_{m_{I_\zeta}(\beta)}\quad \text{and}\quad
\lambda_{(I_\zeta,\beta,\delta)}(\omega_i) = \delta_i\kappa_i,\quad\text{ for all }i\in I\smallsetminus\{o\}. 
\end{equation*}
Then $\lambda_{(\beta,I_\zeta,\delta)}$ is an epipelagic type. With fixed $(\beta,\delta)$, let  $\pi_{I_\zeta}:=\pi_{(\beta,I_\zeta,\delta)}$ be the induced epipelagic supercuspidal representation of $G(F)$.

We also define the following representation of $\GL_{2n}(F)$,
$$\tilde\pi_{(\beta,\delta)} = \prod_{i\in I\smallsetminus\{o\}}\tilde\pi(2\beta_i,\mathbf 1_{\mu_F}, \delta_i), $$
then $\tilde\pi_{(\beta,\delta)}$ is the lifting of the representations in the L-packet
$$\tilde\Pi_{(\beta,\delta)} := \{\pi_{I_\zeta}\}_{I_\zeta},$$
which contains $2^{\#I\smallsetminus\{o\}}$ representations: half of them belongs to $G_+$ and another half to $G_-$. Indeed if $\pi_{{\O}}$ belongs to $G_+$, then 
$$\text{$\pi_{I_\zeta}$ belongs to $G_\epsilon$ \quad$\Leftrightarrow$\quad $\#I_\zeta \equiv \tfrac{1}{2}(\epsilon-1)\bmod 2$,\quad for $\epsilon\in \{+,-\}$.}
$$

\subsubsection{Symplectic groups}

For $G=\SP_{2n}$, given 
\begin{enumerate}[(i)]

\item  a stable functional $\beta= (\beta_i)_{i\in I}$ in $G(F)$, 

\item a partition $I_\zeta$ of $I$, and

\item a tuple of signs $\delta = (\delta_i)_{i\in I}\in \{\pm 1\}^{\#I}$,

\end{enumerate}
define a character $\lambda_{(\beta,I_\zeta,\delta)}$ by 
\begin{equation*}
\lambda_{(\beta,I_\zeta,\delta)} = \psi_{m_{I_\zeta}(\beta)},
\quad\text{and}\quad 
\lambda_{(I_\zeta,\beta,\delta)}(\omega_i) = \delta_i\left(\tfrac{-1}{\mu_F}\right)\kappa_i,\quad\text{ for all }i\in I. 
\end{equation*}
Then $\lambda_{(\beta,I_\zeta,\delta)}$ is an epipelagic type. With fixed $(\beta,\delta)$, let  $\pi_{I_\zeta}:=\pi_{(\beta,I_\zeta,\delta)}$ be the induced epipelagic supercuspidal representation of $G(F)$.

Define a character $\tilde{\boldsymbol{\lambda}}_{\delta_o}$ of $F^\times$ by 
$$ \tilde\lambda_{\delta_o}|_{\mu_F} = \left(\frac{\cdot}{\mu_F}\right)^{\#I } 
\quad\text{and}\quad
\tilde{\boldsymbol{\lambda}}_{\delta_o}(\varpi_{}) = \prod_{i\in I}\delta_i\left(\frac{\varpi\det\beta_i}{\mu_F}\right),$$
as a slight modification of $\tilde{\boldsymbol{\lambda}}_{o}$ in  Section \ref{subsection Liftings of epipelagic representations for classical groups}, and then define an irreducible parabolically induced representation of $\GL_{2n+1}(F)$:
$$\tilde\pi_{(\beta,\delta)} = \left( \prod_{i\in I}\tilde\pi(2\beta_i,\left(\tfrac{\cdot}{\mu_F}\right), \delta_i\mathfrak n_\psi) \right)\times \tilde{\boldsymbol{\lambda}}_{\delta_o}.$$ 
then $\tilde\pi_{(\beta,\delta)}$ is the lifting of the representations in the L-packet
$$\tilde\Pi_{(\beta,\delta)} := \{\pi_{I_\zeta}\}_{I_\zeta},$$which contains $2^{\#I}$ representations.

\subsubsection{Even orthogonal groups}

For $G=\SO^{}_{2n}$, given 
\begin{enumerate}[(i)]

\item  a stable functional $\beta= (\beta_i)_{i\in I}$ in $G(F)$, 

\item a partition $I_\zeta$ of $I$,

\item a tuple of signs $\delta = (\delta_i)_{i\in I}\in \{\pm 1\}^{\#I}$, and

\item an extra sign $\xi\in \{\pm 1\}$,

\end{enumerate}
if $\beta_i = \beta_i(a_1,\dots,a_{k})$ for some $(a_1,\dots,a_{k})$, put 
$$\beta_i(+) = \beta_i
\quad\text{and}\quad\beta_i(-) =\beta_i(a_1,\dots,a_{k-1},-a_{k})$$ 
(see (\ref{action of p on beta, ramified even orthogonal})), and for any tuple $(\xi_i)_i\in \{\pm 1\}^{\#I}$ such that $\prod_{i}\xi_i = \xi$, define 
$$\beta_{\xi} = \bigoplus_{\prod_i\xi_i =\xi}\beta_i(\xi_i ),\quad\xi\in \{+,-\}.$$
Define a character $\lambda_{(\beta,I_\zeta,\delta,\xi)}$ by 
\begin{equation*}
 \lambda_{(\beta,I_\zeta,\delta,\xi)} = \psi_{m_{I_\zeta}(\beta_\xi)}\quad\text{and}\quad 
\lambda_{(I_\zeta,\beta,\delta)}(\omega_i) = \delta_i\kappa_i,\quad\text{ for all }i\in I.  
\end{equation*}
Then $ \lambda_{(\beta,I_\zeta,\delta,\xi)}$ is an epipelagic type. With fixed $(\beta,\delta)$, let  $\pi_{(I_\zeta,\xi)}:=\pi_{(\beta,I_\zeta,\delta,\xi)}$ be the induced epipelagic supercuspidal representation of $G(F)$.

If $o\notin I$, we define an irreducible parabolically induced representation of $\GL_{2n}(F)$, 
$$\tilde\pi_{(\beta,\delta,\xi)} =  \prod_{i\in I}\tilde\pi(\beta_i(\xi_i),\left(\tfrac{\cdot}{\mu_F}\right), \delta_i\mathfrak n_\psi),
$$
then the lifting of the packet
$$\tilde\Pi_{(\beta,\delta)} = \{\pi_{(I_\zeta,\xi)}
\}_{I_\zeta,\xi}$$
lies in $\{\tilde\pi_{(\beta,\delta,\xi)} \}_\xi$. The packet $\tilde\Pi_{(\beta,\delta)}$ contains $2^{\#I+1}$ representations: half of them belongs to $G_+$ and another half to $G_-$. Indeed, for fixed $\xi$, if $\pi_{(\O,\xi)}$ belongs to $G_+$, then 
$$\text{$\pi_{(I_\zeta,\xi)}$ belongs to $G_\epsilon$ \quad$\Leftrightarrow$\quad $\#I_\zeta \equiv \tfrac{1}{2}(\epsilon-1)\bmod 2$,\quad for $\epsilon\in \{+,-\}$.}
$$
Hence each group has a packet of cardinality $2^{\#I}$, which is a union of two L-packets $\tilde\Pi_{(\beta,\delta,+)} \sqcup \tilde\Pi_{(\beta,\delta,-)}$ according to Langlands' philosophy: 
$$\tilde\Pi_{(\beta,\delta,\xi)}
\quad\text{ lifts to } \quad
\tilde\pi_{(\beta,\delta,\xi)},
\quad\text{for }
\xi\in \{\pm \}.$$
Fix an embedding $I_\zeta$, then both $\pi_{(\beta,I_\zeta,\delta,+)}$ and $\pi_{(\beta,I_\zeta,\delta,-)}$ induce to isomorphic representation of $G^\sharp(F).$

If $o\in I$, we define two extra characters ${\tilde{\boldsymbol{\lambda}}}_{\delta_o}$ and ${\tilde{\boldsymbol{\lambda}}}_\delta$ of $F^\times$, tamely ramified and such that 
 \begin{equation*}
   \begin{split}
      \tilde\lambda_{\delta_o} |_{\mu_F}\equiv \mathbf 1_{\mu_F},{\tilde{\boldsymbol{\lambda}}}_{\delta_o}(\varpi) = {\delta_o}
      ;\quad 
     \tilde\lambda_{\delta} |_{\mu_F}\equiv \left(\frac{\cdot}{\mu_F}\right), {\tilde{\boldsymbol{\lambda}}}_{\delta}(\varpi) = \prod_{i\in I}\delta_i,
   \end{split}
\end{equation*}
as slight modifications of ${\tilde{\boldsymbol{\lambda}}}_{1}$ and ${\tilde{\boldsymbol{\lambda}}}_{2}$ in (\ref{split or unramified orthogonal groups, the two charaacters}), define an irreducible parabolically induced representation of $\GL_{2n}(F)$, 
$$\tilde\pi_{(\beta,\xi,\delta)} =  \prod_{i\in I\smallsetminus \{o\}}\tilde\pi(\beta_i(\xi_i),\left(\tfrac{\cdot}{\mu_F}\right), \delta_i\mathfrak n_\psi) 
\times   \tilde{\boldsymbol{\lambda}}_{\delta_o}
\times   \tilde{\boldsymbol{\lambda}}_{\delta}$$  as an irreducible parabolically induced representation of $\GL_{2n}(F)$, which is the lifting of the packet
$$\tilde\Pi_{(\beta,\delta)} = \{\pi_{(I_\zeta,\xi)}
\}_{I_\zeta,\xi}$$
containing $2^{\#I}$ representations: if $\pi_{{(\xi, {\O})}}$ belongs to $G_+$, then 
$$\text{$\pi_{(\xi,I_\zeta )}$ belongs to $G_\epsilon$ \quad$\Leftrightarrow$\quad $\#I_\zeta \equiv \tfrac{1}{2}(\epsilon-1)\bmod 2$,\quad for $\epsilon\in \{+,-\}$.}
$$
Hence each group has a packet of cardinality $2^{\#I-1}$, and is an L-packet. In contrast to the case $o\notin I$, this time all $\pi_{(I_\zeta,\xi)}$ are invariant by $G^\sharp(F)$.

\subsubsection{Ramified unitary groups}

Let $G=\mathrm U_{N}(F/\Fo)$ where $F = \Fo[\varpi]$ with $\varpi^2=-\varpi_\bullet$. Given 
\begin{enumerate}[(i)]

\item  a stable functional $\beta= (\beta_i)_{i\in I}$ in $G(F)$, 

\item a partition $I_\zeta$ of $I$, and

\item a tuple of signs $\delta = (\delta_i)_{i\in I}\in \{\pm 1\}^{\#I}$,

\end{enumerate}
define a character $\lambda_{(\beta,I_\zeta,\delta)}$ by 
\begin{equation*}
\lambda_{(\beta,I_\zeta,\delta)} = \psi_{m_{I_\zeta}(\beta)},
\quad\text{and}\quad 
\lambda_{(I_\zeta,\beta,\delta)}(\omega_i) = \delta_i\kappa_i,\quad\text{ for all }i\in I. 
\end{equation*}
Then $\lambda_{(\beta,I_\zeta,\delta)}$ is an epipelagic type, and let  $\pi_{I_\zeta}:=\pi_{(\beta,I_\zeta,\delta)}$ be the induced epipelagic supercuspidal representation of $G(F)$.

Define
$$\tilde\pi_{(\beta,\delta)} = \prod_{i\in I}\tilde\pi(\beta_i,\left(\tfrac{\cdot}{\mu_F}\right)^{N-1}, \delta_i\mathfrak n_\psi^{N-1}) $$ 
as a representation of $\GL_{N}(F)$. The L-packet of $\tilde\pi_{(\beta,\delta)}$ is then 
$$\tilde\Pi_{(\beta,\delta)} := \{\pi_{I_\zeta}\}_{I_\zeta},$$
which contains $2^{\#I}$ representations: half of them belongs to $G_+$ and another half to $G_-$. Indeed if $\pi_{{\O}}$ belongs to $G_+$, then 
$$\text{$\pi_{I_\zeta}$ belongs to $G_\epsilon$ \quad$\Leftrightarrow$\quad $\#I_\zeta \equiv  \epsilon \bmod 2$,\quad for $\epsilon\in \{+,-\}$.}
$$

\addcontentsline{toc}{section}{References} 
\bibliographystyle{alpha}
\bibliography{epipelagic}

\end{document}